\documentclass[10pt]{amsart}

\usepackage{amsmath, amssymb, amsfonts, amsbsy, color, graphicx}
\usepackage[mathscr]{euscript}
\usepackage{pst-all}

\def\resp{\mbox{\it resp}.\ }

\newcounter{ITEM}
\newcommand\ITEM[1]{\setcounter{ITEM}{#1}\leavevmode\hbox{\rm(\roman{ITEM})}}

\renewcommand\aa{a}
\renewcommand\AA{A}
\newcommand\act{\mathbin{\raisebox{1pt}{$\scriptscriptstyle\bullet$}}}
\newcommand\alter[1]{\hat{#1}}
\newcommand\bb{b}
\newcommand\BB{B}

\newcommand\BDD[1]{B_{#1}^{\oplus}}%OR Dubrovin
\newcommand\BP[1]{B^{\scriptscriptstyle+}_{#1}}% Braids
\newcommand\BR[1]{B_{#1}}%Braids
\newcommand\CAN[1]{\ell_{#1}}
\newcommand\card{\mathtt{\#}}
\newcommand\Cat[1]{\mathcal{C}\hspace{-0.3ex}a\hspace{-0.1ex}t(#1)}%GE Category associated with a germ
\newcommand\cc{c}

\newcommand\CCC{\mathcal{C}}

\newcommand\CCCi{\Isom\CCC}
\newcommand\CCCni{\CCC \setminus \CCCi}%Non-invertible elements
\newcommand\CCCu{\CCC_1}%PR Subcategory
\newcommand\CCCui{\Isom\CCCu}%PR Subcategory

\def\CELL[#1]{[\HS{-0.2}[#1]\HS{-0.2}]}%NF cell
\def\CELLbis[#1]{\langle\HS{-0.5}\langle#1\rangle\HS{-0.5}\rangle}%NF cell
\newcommand\CFMX{\CCC_{\HS{-0.3}\FF\HS{-0.2}}(\MM\HS{-0.2},\HS{-0.3}\XX)}%LD partial action

\newcommand\cl[1]{[#1]}%Equivalence class
\newcommand\clp[1]{[#1]^{\!\scriptscriptstyle+}}%PR positive class
\newcommand\conjt[3]{{ #1\xrightarrow{#2}#3}}%CO, DL triple
\newcommand\dd{d}

\newcommand\DELTA[1]{\Delta^{\![#1]}}%Garside maps
\newcommand\DELTAsym[1]{\symm\Delta_{[#1]}}%Garside maps
\newcommand\Der[2]{#1^{#2}}%GE: Derived germ
\newcommand\dist{\mathrm{dist}}% Combinatorial distance
\renewcommand\div{\divLloc{}}
\newcommand\Div{\mathrm{Div}}
\newcommand\Divl[1]{\mathrm{Div}_{#1}}%Left-divisors, local version
\newcommand\DIV[1]{\mathcal{D}iv(#1)}%Divisors
\newcommand\dive{\mathbin{\preccurlyeq}}%PR

\newcommand\diveLloc[1]{\mathbin{\preccurlyeq_{#1}\nobreak}}%Local version
\newcommand\diveR{\diveRloc{}}
\newcommand\diveRloc[1]{\mathbin{\widetilde{\smash{\preccurlyeq}}_{#1}\nobreak}}
\newcommand\diveRS{\diveRloc\SSSS}
\newcommand\diveRSs{\diveRS^{\HS{-1.2}*}}%Right divisib
\newcommand\diveS{\diveLloc\SSSS}

\newcommand\diveSu{\diveLloc{\SSSSu}}

\newcommand\divLloc[1]{\mathrel{\prec_{#1}\nobreak}}
\newcommand\divR{\divRloc{}}

\newcommand\divRloc[1]{\mathrel{\widetilde{\smash{\prec}}_{#1}\nobreak}}
\newcommand\divRS{\divRloc{\SSSS}}
\newcommand\divRSs{\mathrel{{\divRS\VR(2.2,0)}^{\HS{-1.3}*}\HS{0.5}}}%Right divisib
\newcommand\divS{\divLloc\SSSS}

\newcommand\DIVR[1]{\symm{\VR(2.2,0)\smash{\mathcal{D}iv}}(#1)}%Divisors
\newcommand\Divsym{\widetilde{\vrule width0pt height6pt\smash{\mathrm{Div}}}}%Divisors
\newcommand\double{\psi}%YB
\newcommand\Double{\Psi}%YB
\newcommand\dual{\partial}%Duality
\newcommand\DUAL[1]{\partial^{[#1]}}%Duality
\newcommand\dualD{\partial_{\HS{-0.4}_\Delta\HS{-0.1}}}%DE duality
\newcommand\dualsym{\symm\partial}%DE Duality obsolete
\newcommand\DUALsym[1]{\symm\partial_{[#1]}}%Duality
\newcommand\EE{E}

\newcommand\Env[1]{\mathcal{E}\!nv(#1)}%PR enveloping groupoid
\newcommand\eps{\varepsilon}
\newcommand\eqD{=_{\HS{-0.3}\scriptscriptstyle\Delta}}%CO
\newcommand\eqil{\mathrel{{}^{\scriptscriptstyle\times}\!\!=}}
\newcommand\eqilr{\mathrel{{}^{\scriptscriptstyle\times}{\!=}^{\!\scriptscriptstyle\times}}}
\newcommand\eqir{\mathbin{=^{\!\scriptscriptstyle\times}}}
\newcommand\eqirS{\mathrel{=^{\!\scriptscriptstyle\times}_{\!\SSSS}}}
\newcommand\EQU{\bowtie}%PR Ore's theorem
\newcommand\eqp{\equiv^{\scriptscriptstyle+}}

\newcommand\eqpRR{\equiv_{\HS{-0.1}\RRR}^{\scriptscriptstyle+}\nobreak}

\newcommand\etc{\pdots}

\newcommand\ew{\varepsilon}
\newcommand\fct{\phi}%SU Generic functor
\newcommand\ff{f}
\newcommand\fD{\phi_\Delta}%Garside functor
\newcommand\FF{F}

\newcommand\FFF{\mathcal{F}}%SU family (of subcategories)

\newcommand\fft{\alter{f}}
\renewcommand\ge{\geqslant\nobreak}

\newcommand\Gen{\Sigma}%BR Coxeter generators
\renewcommand\gg{g}
\newcommand\GG{G}
\newcommand\GGG{\mathcal{G}}
\newcommand\ggt{\alter{g}}
\newcommand\hh{h}
\newcommand\HH{H}% Head function
\newcommand\HHH{{\mathcal H}}

\newcommand\hht{\alter{h}}
\newcommand\Ht{\mathrm{ht}}%PR Height
\newcommand\HS[1]{\hspace{#1ex}}
\newcommand\id[1]{1_{\HS{-0.2}{#1}}}
\newcommand\Id[1]{\boldsymbol{1}_{\!#1}}
\newcommand\ie{\epsilon}%PR invertible\newcommand\ii{i}
\newcommand\ii{i}
\newcommand\II{I}
\newcommand\III{\mathscr{I}}
\newcommand\inc{\subset}
\newcommand\ince{\subseteq}

\newcommand\incne{\mathrel{\raise2pt\hbox{$\inc$}\kern-7pt\raise-2pt\hbox{\hskip0.5pt$\scriptstyle\not=$}}}
\newcommand\inv{^{-1}}
\newcommand\INV[1]{\overline{\vrule height 6pt width 0pt #1}}%Reversing
\newcommand\Isom[1]{#1^{\!\times}\!}

\newcommand\jj{j}
\newcommand\JJ{J}
\newcommand\JJJ{\mathscr{J}}

\newcommand\kk{k}
\newcommand\LDAct{F}%LD Action
\newcommand\LDCAT[2]{\CCC(#1, #2)}%LD 
\newcommand\LDDef[1]{\mathrm{Def}(#1)}%LD 
\newcommand\LDM{M_{\LDsmall}}%LD Monoid of LD
\newcommand\LDO[1]{\Sigma_{#1}}%LD Operator
\newcommand\LDpibis{\pi_{\HS{-0.3}\scriptscriptstyle\bullet}\HS{-0.3}}%LD Projection
\newcommand\LDsec{\widetilde\pi}%LD Section
\newcommand\LDsmall{{\HS{-0.3}\scriptscriptstyle L\!D}}%LD 
\newcommand\lD{<_{\HS{-0.2}\scriptscriptstyle\mathrm{D}}}%OR Dehornoy braid ordering
\newcommand\lDe{{\le_{\HS{-0.8}\scriptscriptstyle\mathrm{D}}\,}}%OR Dehornoy braid ordering

\renewcommand\le{\leqslant\nobreak}
\newcommand\length{\mathrm{lg}}
\newcommand\lexp[2]{\kern\scriptspace\vphantom{#2}^{#1}\kern-\scriptspace#2}
\renewcommand\lg[1]{\length(#1)}
\newcommand\LGG[2]{\Vert#2\Vert_{#1}}

\newcommand\LO[1]{\mathrm{LO}(#1)}%OR Space
\newcommand\LT[2]{\Vert#1\Vert_{#2}}%GE length
\newcommand\mm{m}
\newcommand\MM{M}

\newcommand\MMu{\MM_1}%SU submonoid

\newcommand\Mon[1]{\mathcal{M}\hspace{-0.3ex}o\hspace{-0.1ex}n(#1)}%GE Monoid associated with a germ
\newcommand\mult{\succ}

\newcommand\multeR{\mathrel{\widetilde{\smash{\succcurlyeq}}}}

\newcommand\multR{\multRloc{}}
\newcommand\multRloc[1]{\mathrel{\widetilde{\smash{\succ}}_{#1}\nobreak}}%Local version
\newcommand\NF[1]{\mathtt{NF}_{\!#1}}%Family of all normal sequences
\newcommand\nn{n}
\newcommand\NN{N}
\newcommand\NNNN{\mathbb{N}}

\def\NT(#1,#2,#3){\conjt{#1}{#2}{#3}}%LD

\newcommand\Obj{\mathcal{O}\HS{-0.15}b\HS{-0.25}j}
\newcommand\OP{\,{\scriptstyle\bullet}\,}%GE basic operation
\newcommand\OPDer{\OP}%GE operation for derived germ
\newcommand\OPu{\mathbin{{\scriptstyle\bullet}_1}}%GE operation with reference germ 
\newcommand\opR{\star}%YB cyclic operation
\newcommand\opRu{\underline\opR}%YB cyclic operation
\newcommand\pdots{\HS{0.2}{\cdot}{\cdot}{\cdot}\HS{0.2}}

\newcommand\Pow[2]{#1^{#2}}%Power
\newcommand\pp{p}
\newcommand\PP{P}
\newcommand\Pref[1]{\mathrel{{\le}_{#1}}\nobreak}%GE prefix
\newcommand\Pres[2]{(#1, #2)}
\newcommand\PRESp[2]{\langle#1\,\vert\, #2\rangle^{\scriptscriptstyle\!+}\!}% Presented monoid
\newcommand\Prod[3]{(#1\sep #2)^{[#3]}}%BR porduct
\newcommand\qq{q}

\newcommand\RC{\theta}%PR right-complement
\newcommand\RCt{\theta^*_3}

\newcommand\rep[1]{#1{}^0}%CO representative
\newcommand\rev{\curvearrowright}

\newcommand\revRh{\mathrel{\curvearrowright_{\!\RRh}}}

\newcommand\rn{r}%OR "r number" (avoid)
\newcommand\rr{r}
\newcommand\RR{R}
\newcommand\RRh{\widehat\RR}
\newcommand\RRR{{\mathcal{R}}}

\newcommand\sep{\HS{0.05}\vert\HS{0.05}}%PR path separator
\newcommand\sepp{\relax}%optional path separator
\newcommand\seq[1]{#1}
\newcommand\seqq[2]{#1\sep#2}
\newcommand\seqqq[3]{#1\sep#2\sep#3}
\newcommand\seqqqq[4]{#1\sep#2\sep#3\sep#4}

\newcommand\seqqqqqq[6]{#1\sep#2\sep#3\sep#4\sep#5\sep#6}

\newcommand\Seq[2]{\VR(3.5,0)\smash{#1^{[#2]}}}%BR Artin generators
\newcommand\sh[1]{{#1^{\scriptstyle\sharp}}}% Closure under right-multiplication by invertible elements
\newcommand\shift{\mathrm{sh}}%OR shift
\newcommand\sig[1]{\sigma_{\hspace{-0.2ex}#1}^{\null}}

\newcommand\sigg[2]{\sigma_{#1}^{#2}}
\newcommand\siginv[1]{\sigma_{#1}^{-1}}

\newcommand\sol{\rho}%YB solution YBE
\newcommand{\sss}{s}
\newcommand{\SSS}{S}

\newcommand\SSSS{\mathcal{S}}

\newcommand\SSSSg{\underline\SSSS}%GE : germ S
\newcommand\SSSSgu{\SSSSg_1}%GE : germ S
\newcommand\SSSs{\sh{\SSS}}
\newcommand\SSSSs{{\sh\SSSS}}
\newcommand\SSSSi{\Isom\SSSS}
\newcommand\SSSSu{\SSSS_1}%PR Subfamily
\newcommand\SSSSus{\sh{\SSSSu}}%SU subfamily
\newcommand\Sub[1]{\mathcal{S}\!ub(#1)}%GE subcategory
\newcommand\Suff[1]{\mathrel{\widetilde{\VR(1.8,0)\smash\le}_{#1}}\nobreak}%GE suffix
\newcommand\SUP[1]{\mathrm{sup}_{#1}}% supremum
\newcommand\symm{\widetilde}%symbol for right counterparts

\newcommand\term{T}%LD term
\renewcommand\tt{t}
\newcommand\tta{\mathtt{a}}

\newcommand\ttb{\mathtt{b}}

\newcommand\ttc{\mathtt{c}}

\newcommand\ttd{\mathtt{d}}
\newcommand\tte{\mathtt{e}}

\newcommand\ttf{\mathtt{f}}
\newcommand\ttg{\mathtt{g}}
\newcommand\tth{\mathtt{h}}
\newcommand\tti{\mathtt{i}}
\newcommand\ttj{\mathtt{j}}
\newcommand\ttk{\mathtt{k}}

\newcommand\ttm{\mathtt{m}}
\newcommand\ttn{\mathtt{n}}
\newcommand\ud{\hbox{-}}%for undefined variable
\newcommand\und[1]{\underline#1}%LD
\newcommand\under{\backslash}
\newcommand\uu{u}

\newcommand\uut{\alter{\uu}}

\newcommand\vs{{\it vs.}\ }
\newcommand\vv{v}

\def\VR(#1,#2){\vrule width0pt height#1mm depth#2mm}
\newcommand\vvt{\alter{\vv}}

\newcommand\wAt{\widetilde{\VR(2,0)\smash{\mathrm{A}}}_2}%NF affine braids
\newcommand\wdots{, ...\hspace{0.2ex},}
\newcommand\wit{\lambda}%Noetherianity witness
\newcommand\ww{w}
\newcommand\WW{W}
\newcommand\xx{x}
\newcommand\XX{X}

\newcommand\XXX{\mathcal{X}}%GA Noetherianity
\newcommand\yy{y}
\newcommand\YY{Y}

\newcommand\zz{z}
\newcommand\ZZ{Z}
\newcommand\ZZZZ{{\mathbb Z}}

\psset{unit=1mm}
\psset{fillcolor=white}
\psset{dotsep=1.5pt}
\psset{dash=1.5pt 1.5pt}
\psset{linewidth=0.8pt}
\psset{arrowsize=3.5pt}
\psset{doublesep=1pt}
\psset{coilwidth=1.2mm}
\psset{coilheight=2}
\psset{coilaspect=20}
\psset{coilarm=2}
\newpsstyle{eps}{linewidth=0.5pt, doubleline=true}
\newpsstyle{double}{linewidth=0.5pt, doubleline=true}
\newpsstyle{etc}{linestyle=dotted}
\newpsstyle{exist}{linestyle=dashed}
\newpsstyle{thick}{linewidth=1.5pt}
\newpsstyle{thickexist}{linewidth=1.5pt, linestyle=dashed}
\newpsstyle{thin}{linewidth=0.5pt}
\newpsstyle{thinexist}{linewidth=0.5pt, linestyle=dashed}

\def\BPoint(#1,#2,#3){\cnode[style=thin,fillcolor=black,fillstyle=solid](#1,#2){0.5}{#3}}
\def\GPoint(#1,#2,#3){\cnode[style=thin,fillcolor=lightgray,fillstyle=solid](#1,#2){0.5}{#3}}
\def\WPoint(#1,#2,#3){\cnode[style=thin,fillcolor=white,fillstyle=solid](#1,#2){0.5}{#3}}
\def\AArrow(#1,#2){\ncline[nodesep=0.3mm, linewidth=0.8pt, border=1.2pt]{->}{#1}{#2}}
\def\BArrow(#1,#2){\ncline[linecolor=blue, linewidth=1.2pt, linestyle=dotted, nodesep=0.3mm, border=1.2pt]{->}{#1}{#2}}
\def\CArrow(#1,#2){\ncline[linecolor=red, nodesep=0.3mm, linewidth=0.8pt, border=1.2pt]{->}{#1}{#2}}
\def\DArrow(#1,#2){\ncline[linecolor=red,style=thickexist,nodesep=0.5mm,border=2pt]{->}{#1}{#2}}

\def\Arrow(#1,#2){\ncline[nodesep=1mm,border=2pt]{->}{#1}{#2}}
\def\ArrowE(#1,#2){\ncline[nodesep=1mm,border=2pt,style=exist]{->}{#1}{#2}}
\def\ETC(#1,#2){\ncline[nodesep=2mm,border=2pt,style=etc]{#1}{#2}}

\def\Point(#1,#2,#3){\pnode(#1,#2){#3}}
\def\PointN(#1,#2,#3,#4){\pnode(#1,#2){#3}{#4}}
\def\SetGrid{\Point(0,0,00)\Point(15,0,10)\Point(30,0,20)\Point(45,0,30)\Point(60,0,40)\Point(75,0,50)\Point(90,0,60)\Point(105,0,70)\Point(120,0,80)
\Point(0,12,01)\Point(15,12,11)\Point(30,12,21)\Point(45,12,31)\Point(60,12,41)\Point(75,12,51)\Point(90,12,61)\Point(105,12,71)\Point(120,12,81)
\Point(0,24,02)\Point(15,24,12)\Point(30,24,22)\Point(45,24,32)\Point(60,24,42)\Point(75,24,52)\Point(90,24,62)\Point(105,24,72)\Point(120,24,82)
\Point(0,36,03)\Point(15,36,13)\Point(30,36,23)\Point(45,36,33)\Point(60,36,43)\Point(75,36,53)\Point(90,36,63)\Point(105,36,73)\Point(120,36,83)
\Point(0,48,04)\Point(15,48,14)\Point(30,48,24)\Point(45,48,34)\Point(60,48,44)\Point(75,48,54)\Point(90,48,64)\Point(105,48,74)\Point(120,48,84)
\Point(0,60,05)\Point(15,60,15)\Point(30,60,25)\Point(45,60,35)\Point(60,60,45)\Point(75,60,55)\Point(90,60,65)\Point(105,60,75)\Point(120,60,85)}

\allowdisplaybreaks

\title{Addenda to ``Foundations of Garside Theory''}

\author{Patrick DEHORNOY}
\address{Laboratoire de Math\'ematiques Nicolas Oresme,
CNRS UMR 6139, Universit\'e de Caen, 14032 Caen, France}
\email{patrick.dehornoy@unicaen.fr}
\urladdr{www.math.unicaen.fr/\!\hbox{$\sim$}dehornoy}
\author{\\\lowercase{with}\\ Fran\c cois DIGNE}
\address{Laboratoire Ami\'enois de Math\'ematique Fondamentale et Appliqu\'ee,
CNRS UMR 7352, Universit\'e de Picardie
Jules-Verne, 80039 Amiens, France}
\email{digne@u-picardie.fr}
\urladdr{www.mathinfo.u-picardie.fr/digne/}
\author{Eddy GODELLE}
\address{Laboratoire de Math\'ematiques Nicolas Oresme,
CNRS UMR 6139, Universit\'e de Caen, 14032 Caen, France}
\email{godelle@math.unicaen.fr}
\urladdr{www.math.unicaen.fr/\!\hbox{$\sim$}godelle}
\author{Daan KRAMMER}
\address{Mathematics Institute, University of Warwick, Coventry CV4
7AL, United Kingdom}
\email{D.Krammer@warwick.ac.uk}
\urladdr{www.warwick.ac.uk/\!\hbox{$\sim$}masbal/}
\author{Jean MICHEL}
\address{Institut de Math\'ematiques de Jussieu, CNRS UMR 7586,
Universit\'e Denis Diderot Paris 7, 75205 Paris 13, France}
\email{jmichel@math.jussieu.fr}
\urladdr{www.math.jussieu.fr/\!\hbox{$\sim$}jmichel/}

\makeatletter
\renewcommand{\ps@plain}{%
\renewcommand{\@oddhead}{Addenda to ``Foundations of Garside Theory'', by Dehornoy {\it et al.}\hfill\thepage}%
\renewcommand{\@evenhead}{\thepage\hfill Addenda to ``Foundations of Garside Theory'', by Dehornoy {\it et al.}}%
\renewcommand{\@oddfoot}{}%
\renewcommand{\@evenfoot}{}}
\makeatother
\newcommand\chapter[1]{\vskip 15mm plus 5mm minus 5mm\noindent{\LARGE\bf #1}\vskip 5mm}

\newenvironment{ADprop}[3]
{\bigskip\noindent{\bf Proposition}~{#1} ({\bf#2}).---\ \it}
{}

\newenvironment{ADlemm}[2]
{\bigskip\noindent{\bf Lemma}~{#1}.---\ \it}
{}

\newcommand\ADlabel[1]{\tag{{#1}}}

\newenvironment{exer}[3]
{\bigskip\noindent{\bf Exercice}~{\bf #1} ({\bf#2}).---\ \it}
{}

\newenvironment{solu}
{\medskip\noindent{\it Solution.}\ }
{}

\begin{document}

\maketitle

\begin{abstract}
This text consists of additions to the book ``Foundations of Garside Theory'', EMS Tracts in Mathematics, vol. 22 (2015)---see introduction and table of contents in arXiv:1309.0796---namely skipped proofs and solutions to selected exercises.
\end{abstract}

%%%%
\chapter{Chapter~I: Examples}

%%%%
\section*{Skipped proofs}

\noindent(none)

%%%%
\section*{Solution to selected exercises}

\noindent(none)

%%%%
\chapter{Chapter~II: Preliminaries}

%%%%
\section*{Skipped proofs}

\begin{ADprop}{II.1.18}{collapsing invertible elements}{}
For $\CCC$ a left-cancellative category, the following conditions  are equivalent:

\ITEM1 The equivalence relation~$\eqir$ is compatible with composition, in the sense that, if $\gg_1 \gg_2$ is defined and $\gg'_\ii \eqir \gg_\ii$ holds for $\ii = 1,2$, then $\gg'_1 \gg'_2$ is defined and $\gg'_1 \gg'_2 \eqir \gg_1 \gg_2$ holds;

\ITEM2 The family $\CCCi(\xx , \yy)$ is empty for $\xx \not= \yy$ and, for all $\gg, \gg'$ sharing the same source, $\gg \eqilr \gg'$ implies $\gg \eqir \gg'$;

\ITEM3 The family $\CCCi(\xx , \yy)$ is empty for $\xx \not= \yy$ and we have 
\begin{equation}
\ADlabel{II.119}
\forall\xx,\yy \in \Obj(\CCC) \ \forall\gg\in \CCC(\xx, \yy)\ \forall\ie \in \CCCi(\xx, \xx)\ \exists\ie' \in \CCCi(\yy, \yy)\ (\ie \gg = \gg \ie').
\end{equation}
When the above conditions are satisfied, the equivalence relation~$\eqilr$ is compatible with composition, and the quotient-category~$\CCC{/}{\eqilr}$ has no nontrivial invertible element.
\end{ADprop}

\begin{proof}
Assume \ITEM1. Let $\ie \in \CCCi(\xx, \yy)$. Then $\ie \eqir \id\xx$ and $\id\yy \eqir \id\yy$ are satisfied, and $\ie \, \id\yy$ is defined. By~\ITEM1, $\id\xx \, \id\yy$ must be defined as well, which is possible only for $\xx = \yy$. Let now $\gg \in \CCC(\xx, \yy)$ and $\ie \in \CCCi(\xx, \xx)$. Then $\ie \eqir \id\xx$ and $\gg \eqir \gg$ are satisfied, and $\ie \, \gg$ is defined. By~\ITEM1, we must have $\gg \eqir \ie \gg$, that is, there must exist~$\ie'$ in~$\CCC(\yy, \yy)$ satisfying $\ie \gg = \gg \ie'$. So \ITEM1 implies~\ITEM3.

Assume now \ITEM2. Let $\gg \in \CCC(\xx, \yy)$ and $\ie \in \CCCi(\xx, \xx)$. Then $\gg \eqilr \ie \gg$ holds, as we can write $\ie \, \gg = \ie\gg \, \id\yy$. By~\ITEM2, we deduce $\gg \eqir \ie\gg$, so, as above, there must exist~$\ie'$ in~$\CCC(\yy, \yy)$ satisfying $\ie \gg = \gg \ie'$. So \ITEM2 implies~\ITEM3.

Assume now \ITEM3. Let $\gg_1 \in \CCC(\xx, \yy)$, $\gg_2 \in \CCC(\yy, \zz)$, and assume $\gg'_1 \eqir \gg_1$ and $\gg'_2 \eqir \gg_2$. By assumption, there exists $\ie_\ii$ in~$\CCCi$ satisfying $\gg'_\ii = \gg_\ii \ie_\ii$ for $\ii = 1,2$. Applying~(II.1.19) to~$\gg_2$ and~$\ie_1$, we deduce that there exists~$\ie_1'$ in~$\CCCi(\zz, \zz)$ satisfying $\ie_1 \gg_2 = \gg_2 \ie_1'$. We deduce $\gg'_1 \gg'_2 = \gg_1 \gg_2 \ie_1' \ie_2$, whence $\gg'_1 \gg'_2 \eqir \gg_1 \gg_2$. So \ITEM3 implies~\ITEM1.

Next, assume \ITEM3 again, and $\gg' \eqilr \gg$. By definition, there exist $\ie, \ie'$ satisfying $\ie \gg' = \gg \ie'$. By~\ITEM3, there exists an invertible element~$\ie''$ satisfying $\ie \gg' = \gg' \ie''$, and we deduce $\gg' = \gg \, \ie' \ie''{}\inv$, whence $\gg' \eqir \gg$. So \ITEM3 implies~\ITEM2.

Finally, assume that \ITEM1, \ITEM2, and~\ITEM3 are satisfied, we have $\gg_1 \eqilr \gg'_1$ and $\gg_2 \eqilr\nobreak \gg'_2$, and $\gg_1 \gg_2$ and $\gg'_1 \gg'_2$ exist. By assumption, there exist for $\ii = 1,2$ invertible elements $\ie_\ii, \ie'_\ii$ satisfying $\ie'_\ii \gg_\ii = \gg'_\ii \ie_\ii$. Let $\xx$ and $\yy$ be the source and target of~$\gg_2$. By construction, $\ie_1\inv \ie'_2$ belongs to~$\CCCi(\xx, \xx)$. Applying~\ITEM3, we deduce the existence of~$\ie$ in~$\CCCi(\yy, \yy)$ satisfying $\ie_1\inv \ie'_2 \gg_2 = \gg_2 \ie$. Then we obtain
$$\ie'_1 \, \gg_1 \gg_2 = \ie'_1 \, \gg_1 \gg_2 \, \ie = \ie'_1 \gg_1 
\ie_1\inv \, \ie'_2 \gg_2 = \gg'_1 \gg'_2 \,\ie_2,$$
which shows that $\gg_1 \gg_2 \eqilr \gg'_1 \gg'_2$ is true. So $\eqilr$ is a congruence, and there exists a well defined quotient category~$\CCC/{\eqilr}$, obtained from~$\CCC$ by identifying elements that are $\eqilr$-equivalent. In the current case, according to~\ITEM2, distinct objects of~$\CCC$ are never connected by invertible elements, so the collapsing only involves the elements. Then, by construction, the category $\CCC/{\eqilr}$ has no nontrivial invertible element.
\end{proof}

%%%%
\section*{Solution to selected exercises}

%%%%
\begin{exer}{4}{atom}{}
Assume that $\CCC$ is a left-cancellative category. Show that,  for $\nn \ge 1$,  every element~$\gg$ of~$\CCC$ satisfying $\Ht(\gg) = \nn$ admits a decomposition into a product of $\nn$~atoms. 
\end{exer}

\begin{solu}
By Proposition~II.2.48, there exists a decomposition $(\gg_1 \wdots \gg_\nn)$ of~$\gg$ consisting of $\nn$~non-invertible entries. The assumption that $\gg_\ii$ is non-invertible implies $\Ht(\gg_\ii) \ge 1$ for every~$\ii$. Hence we must have $\Ht(\gg_\ii) = 1$ for each~$\ii$, so each factor~$\gg_\ii$ is an atom. Hence $\gg$ admits a decomposition as a product of $\nn$~atoms. 
\end{solu}

%%%%
\begin{exer}{5}{unique right-mcm}{Z:PRRightMcm}
Assume that $\CCC$ is a left-cancellative category that admits right-mcms. Assume that$\ff, \gg$ are elements of~$\CCC$ that admit a common right-multiple and any two right-mcms of~$\ff$ and~$\gg$ are $\eqir$-equivalent. Show that every right-mcm of~$\ff$ and~$\gg$ is a right-lcm of~$\ff$ and~$\gg$. \end{exer} 

\begin{solu}
Let $\hh$ be a right-mcm of~$\ff$ and~$\gg$, and $\hht$ be a common right-multiple of~$\ff$ and~$\gg$. Then $\hht$ is a right-multiple of some right-mcm~$\hh'$ of~$\ff$ and~$\gg$. By assumption, $\hh' \eqir \hh$ holds, hence $\hht$ is a right-multiple of~$\hh$ as well. Hence $\hh$ is a right-lcm of~$\ff$ and~$\gg$.
\end{solu}

%%%%
\begin{exer}{6}{right-gcd to right-mcm}{Z:PRRightGcdRightMcm}
Assume that $\CCC$ is a cancellative category that admits right-gcds, $\ff, \gg$ are elements of~$\CCC$ and $\hh$ is a common right-multiple of~$\ff$ and~$\gg$. Show that there exists a right-mcm~$\hh_0$ of~$\ff$ and~$\gg$ such that every common right-multiple of~$\ff$ and~$\gg$ that left-divide~$\hh$ is a right-multiple of~$\hh_0$.
\end{exer}

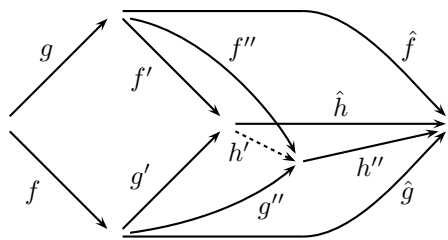
\begin{figure}[htb]
\begin{picture}(60,32)(0,0)
\pcline{->}(1,14)(14,1)
\put(3,5){$\ff$}
\pcline{->}(1,16)(14,29)
\put(5,24){$\gg$}
\pcline{->}(16,1)(29,14)
\put(17,20){$\ff'$}
\pcline{->}(16,29)(29,16)
\put(17,7){$\gg'$}
\psline(16,0)(40,0)
\psbezier{->}(40,0)(47,0)(54,9)(59,14)
\put(53,24){$\fft$}
\psline(16,30)(40,30)
\psbezier{->}(40,30)(47,30)(54,21)(59,16)
\put(53,5){$\ggt$}
\psline[style=exist]{->}(31,14)(39,10)
\psbezier{->}(17,29)(28,27)(37,15)(39,11)
\psbezier{->}(17,0.5)(28,2)(35,6)(39,9.5)
\put(30,24){$\ff''$}
\put(34,3){$\gg''$}
\pcline{->}(31,15)(59,15)
\taput{$\hht$}
\pcline{->}(40,10)(58,14)
\put(47,8){$\hh''$}
\put(30,10){$\hh'$}
\end{picture}
\caption[]{\sf\small Solution to Exercise~6}
\label{F:PRRightMcm}
\end{figure}

\begin{solu}
(See Figure~\ref{F:PRRightMcm}.) Write $\hh = \ff \ggt = \gg \fft$, and let $\hht$ be a right-gcd of~$\fft$ and~$\ggt$. By definition $\hht$ right-divides~$\fft$ and~$\ggt$, so there exist~$\ff', \gg'$ satisfying $\fft = \ff' \hht$ and $\ggt = \gg' \hht$. Then we have $\ff \gg' \hht = \ff \ggt = \gg \fft = \gg \ff' \hht$, whence $\ff \gg' = \gg \ff'$ by right-cancelling~$\hht$. Assume $\ff \gg'' = \gg \ff'' \dive \hh$, say $\hh = \ff \gg'' \hh''$. By left-cancelling~$\ff$, we deduce $\ggt = \gg'' \hh''$ and, similarly, $\fft = \ff'' \hh''$. So $\hh''$ is a common right-divisor of~$\fft$ and~$\ggt$, hence it is a right-divisor of~$\hht$, that is, there exists~$\hh'$ satisfying $\hht = \hh' \hh''$. This implies $\ff \gg'' \hh'' = \ff \ggt$, whence $\gg'' \hh'' = \ggt$ by left-cancelling~$\ff$, and, finally, $\hh'' = \ff \gg'' \dive \ff \ggt = \hh$. So every common right-multiple of~$\ff$ and~$\gg$ that left-divides~$\hh$ is a right-multiple of~$\ff \gg'$. 

Now assume $\ff \gg'' = \gg \ff'' \dive \ff \gg'$. A fortiori, we have $\ff \gg'' = \gg \ff'' \dive \hh$, so the above result implies $\ff \gg' \dive \ff \gg''$, whence $\ff \gg'' \eqir \ff \gg'$. So $\ff \gg'$ is a right-mcm of~$\ff$ and~$\gg$.
\end{solu}

%%%%
\begin{exer}{8}{conditional right-lcm}{Z:PRLocalRightLcm}
Assume that $\CCC$ is a left-cancellative category. \ITEM1 Show that every left-gcd of~$\ff\gg_1$ and~$\ff\gg_2$ (if any) is of the form~$\ff\gg$ where $\gg$ is a left-gcd of~$\gg_1$ and~$\gg_2$. \ITEM2 Assume moreover that $\CCC$ admits conditional right-lcms. Show that, if $\gg$ is a left-gcd of~$\gg_1$ and~$\gg_2$ and $\ff\gg$ is defined, then $\ff\gg$ is a left-gcd of~$\ff \gg_1$ and~$\ff \gg_2$.
\end{exer}

\begin{solu}
\ITEM1 Since $\ff$ left-divides~$\ff\gg_1$ and~$\ff\gg_2$, it left-divides every left-gcd of~$\ff\gg_1$ and~$\ff\gg_2$, so the latter can be written~$\ff\gg$. Now assume that $\ff\gg$ is a left-gcd of~$\ff\gg_1$ and~$\ff\gg_2$. As $\CCC$ is left-cancellative, $\ff\gg \dive \ff\gg_\ii$ implies $\gg \dive \gg_\ii$. Next, assume that $\hh$ left-divides~$\gg_1$ and~$\gg_2$. Then $\ff\hh$ left-divides~$\ff\gg_1$ and~$\ff\gg_2$, implying $\ff\hh \dive \ff\gg$, whence $\hh \dive \gg$. So $\gg$ is a left-gcd of~$\gg_1$ and~$\gg_2$.

\ITEM2 It is clear that $\ff \gg$ is a common left-divisor of~$\ff \gg_1$ and~$\ff\gg_2$. Conversely, assume that $\hh$ is a common left-divisor of~$\ff \gg_1$ and~$\ff \gg_2$. By assumption, $\ff$ and~$\hh$ admit~$\ff \gg_1$ as a common right-multiple, so they admit a right-lcm, say~$\ff \hh'$. By assumption, we have $\ff \hh' \dive \ff \gg_1$ and $\ff \hh' \dive \ff \gg_2$, whence $\hh' \dive \gg_1$ and~$\hh' \dive \gg_2$ by left-cancelling~$\ff$. This in turn implies $\hh' \dive \gg$ since $\gg$ is a left-gcd of~$\gg_1$ and~$\gg_2$. Hence we deduce $\hh \dive \ff \hh' \dive \ff \gg$, which shows that $\ff \gg$ is a left-gcd of~$\ff\gg_1$ and~$\ff\gg_2$.
\end{solu}

%%%%
\begin{exer}{9}{left-coprime}{Z:PRCoprime}
Assume that $\CCC$ is a left-cancellative category. Say that two elements~$\ff, \gg$ of~$\CCC$ sharing the same source~$\xx$ are \emph{left-coprime} if $\id\xx$ is a left-gcd of~$\ff$ and~$\gg$. Assume that $\gg_1, \gg_2 $ are elements of~$\CCC$ sharing the same source, and $\ff \gg_1$ and~$\ff \gg_2$ are defined. Consider the properties \ITEM1 The elements~$\gg_1$ and~$\gg_2$ are left-coprime; \ITEM2 The element~$\ff$ is a left-gcd for~$\ff \gg_1$ and~$\ff \gg_2$. Show that \ITEM2 implies~\ITEM1 and that, if $\CCC$ admits conditional right-lcms, \ITEM1 implies~\ITEM2. [Hint: Use Exercise~8.]
\end{exer}

\begin{solu}
If $\gg$ is a non-invertible common left-divisor of~$\gg_1$ and~$\gg_2$, then $\ff \gg$ is a non-invertible common left-divisor of~$\ff \gg_1$ and~$\ff \gg_2$, so clearly \ITEM2 implies~\ITEM1.

Conversely, assume that $\CCC$ admits conditional right-lcms and $\gg_1, \gg_2$ are left-coprime. By definition, $\id\xx$ is a left-gcd of~$\gg_1$ and~$\gg_2$. Hence, by Exercise~8, $\ff$ is a left-gcd of~$\ff \gg_1$ and~$\ff \gg_2$. So \ITEM1 implies~\ITEM2 in this case.
\end{solu}

%%%%
\begin{exer}{10}{subgroupoid}{Z:PRSubgroupoid}
Let $\MM = \PRESp{\tta, \ttb, \ttc, \ttd}{\tta\ttb = \ttb\ttc = \ttc\ttd = \ttd\tta}$ and $\Delta = \tta\ttb$. \ITEM1 Check that $\MM$ is a Garside monoid with Garside element $\Delta$. \ITEM2 Let $\MM_1$ (\resp $\MM_2$) be the submonoid of~$\MM$ generated by~$\tta$ and~$\ttc$ (\resp $\ttb$ and~$\ttc$). Show that $\MM_1$ and $\MM_2$ are free monoids of rank~$2$ with intersection reduced to~$\{1\}$. \ITEM3 Let $\GG$ be the group of fractions of~$\MM$. Show that the intersection of the subgroups of~$\GG$ generated by~$\MM_1$ and~$\MM_2$ is not~$\{1\}$.
\end{exer}

\begin{solu}
\ITEM2 No word in~$\{\tta, \ttc\}^*$ is eligible for any of the defining relations of~$\MM$, so two distinct such words represent distinct elements of~$\MM_1$. \ITEM3 In~$\GG$, we have $\ttc\inv \tta = \ttd \ttb\inv$ but $\ttd \ttb\inv$ cannot be expressed as $\ff\inv \gg$ with~$\ff, \gg$ in~$\MM_2$ since, otherwise, $\MM_2$ would not be free.
\end{solu}

%%%%
\begin{exer}{11}{weakly right-cancellative}{Z:PRWeaklyRC}
Say that a category~$\CCC$ is weakly right-cancellat\-ive if $\gg \hh = \hh$ implies that $\gg$ is invertible. \ITEM1 Observe that a right-cancellative category is weakly right-cancellative; \ITEM2 Assume that $\CCC$ is a left-cancel\-lative category. Show that $\CCC$ is weakly right-cancellative if and only if, for all~$\ff, \gg$ in~$\CCC$, the relation $\ff \divR \gg$ is equivalent to the conjunction of $\ff \diveR \gg$ and ``$\gg = \gg'\ff$ holds for no~$\gg'$ in~$\CCCi$''.
\end{exer}

\begin{solu}
\ITEM2 The conjunction of $\ff \diveR \gg$ and ``$\gg = \gg'\ff$ holds for no~$\gg'$ in~$\CCCi$'' always implies $\ff \divR \gg$. Assume that $\CCC$ is weakly right-cancellative and $\gg \multR \ff$ holds. Then we have $\gg = \gg' \ff$ for some~$\gg' \notin\CCCi$. Assume $\gg = \gg'' \ff$: if $\gg''$ is invertible, we deduce $\ff = \gg''{}\inv \gg = \gg''{}\inv \gg' \gg$. The assumption that $\CCC$ is weakly right-cancellative implies that $\gg''{}\inv \gg'$ is invertible, hence that $\gg'$ is invertible. So $\ff \divR \gg$ implies ``$\gg = \gg'\ff$ holds for no~$\gg'$ in~$\CCCi$''.
Conversely, assume that $\ff \divR \gg$ is equivalent to the conjunction of $\ff \diveR \gg$ and ``$\gg = \gg'\ff$ holds for no~$\gg'$ in~$\CCCi$''. Assume $\gg \hh = \hh$. If $\xx$ is the source of~$\hh$, we have $\id\xx \hh = \hh$ and $\id\xx$ is invertible. Then the assumption implies that $\hh \divR \hh$ fails, which, as $\hh \dive \hh$ is true, means that $\hh = \gg\hh$ holds for no non-invertible~$\gg$.
\end{solu}

%%%%
\begin{exer}{13}{increasing sequences}{Z:PRIncSeq}
Assume that $\CCC$ is a left-cancellative. For~$\SSSS$ included in~$\CCC$, put $\Div_\SSSS(\hh) = \{\ff \in \CCC \mid \exists \gg \in \SSSS (\ff \gg = \hh)\}$. Show that the restriction of~$\divR$ to~$\SSSS$ is well-founded (that is, admits no infinite descending sequence) if and only if, for every~$\gg$ in~$\SSSS$, every strictly increasing sequence in~$\Div_\SSSS(\gg)$ with respect to left-divisibility is finite.
\end{exer}

\begin{solu}
Assume that $\SSSS$ is not right-Noetherian. Let $\gg_0, \gg_1,...$ be an infinite descending sequence with respect to proper right-divisibility in~$\SSSS$. For each~$\ii$, choose a (necessarily non-invertible) element~$\ff_\ii$ satisfying $\gg_{\ii-1}= \ff_\ii \gg_\ii$. Then we have $\gg_0 = \ff_1 \gg_1 = (\ff_1\ff_2) \gg_2 = ...$, and the sequence $\id\xx$ ($\xx$ the source of~$\gg_0$), $\ff_1, \ff_1\ff_2, ... $ is $\div$-increasing in~$\Div_\SSSS(\gg_0)$.

Conversely, assume that $\gg_0$ lies in~$\SSSS$ and $\hh_1 \div \hh_2 \div ... $ is a strictly increasing sequence in~$\Div_\SSSS(\gg_0)$. Then, for each~$\ii$, there exists a non-invertible element~$\ff_\ii$ satisfying $\hh_\ii \ff_\ii = \hh_{\ii+1}$. On the other hand, as $\hh_\ii$ belongs to~$\Div_\SSSS(\gg_0)$, there exists~$\gg_\ii$ in~$\SSSS$ satisfying $\hh_\ii \gg_\ii = \gg_0$. We find $\gg_0 =\hh_\ii \gg_\ii = \hh_{\ii+1} \gg_{\ii+1} = \hh_\ii \ff_\ii \gg_{\ii+1}$. By left-cancelling~$\hh_\ii$, we deduce $\gg_\ii = \ff_\ii \gg_{\ii+1}$, hence $\gg_{\ii+1}$ is a proper right-divisor of~$\gg_\ii$ for each~$\ii$. So the sequence $\gg_0, \gg_1, ...$ witnesses that $\SSSS$ is not right-Noetherian. 
\end{solu}

%%%
\begin{exer}{15}{left-generating}{Z:PRLeftGenerate}
Assume that $\CCC$ is a left-cancellative category that is right-Noetherian. Say that a subfamily~$\SSSS$ of~$\CCC$ \emph{left-gene\-rates} (\resp \emph{right-generates})~$\CCC$ if every non-invertible element of~$\CCC$ admits at least one non-invertible left-divisor (\resp right-divisor) belonging to~$\SSSS$. \ITEM1 Show that $\CCC$ is right-generated by its atoms. \ITEM2 Show that, if $\SSSS$ is a subfamily of~$\CCC$ that left-generates~$\CCC$, then $\SSSS \cup \CCCi$ generates~$\CCC$. \ITEM3 Conversely, show that, if $\SSSS \cup \CCCi$ generates~$\CCC$ and $\CCCi \SSSS \subseteq \SSSSs$ holds, then $\SSSS$ left-generates~$\CCC$.
\end{exer}

\begin{solu}
\ITEM2 Assume that $\SSSS$ left-generates~$\CCC$. Let $\gg$ be an arbitrary element of~$\CCC$. If $\gg$ is invertible, $\gg$ belongs to~$\CCCi$. Otherwise, by assumption, there exist a non-invertible element~$\gg_1$ in~$\SSSS$ and~$\gg'$ in~$\CCC$ satisfying $\gg = \gg_1 \gg'$. If $\gg'$ is invertible, $\gg$ belongs to~$\SSSS \, \CCCi$. Otherwise, there exist a non-invertible element~$\gg_2$ in~$\SSSS$ and~$\gg''$ satisfying $\gg' = \gg_2 \gg''$, and so on. By Proposition~II.2.28, the sequence $\id\xx$, $\gg_1$, $\gg_1\gg_2, \pdots$, which is increasing with respect to proper left-divisibility and lies in~$\Div(\gg)$, must be finite, yielding $\ell$ and $\gg = \gg_1 \pdots \gg_\ell \ie$ with $\gg_1 \wdots \gg_\ell$ in~$\SSSS$ and $\ie$ in~$\CCCi$.

\ITEM3 Assume that $\SSSS \cup \CCCi$ generates~$\CCC$ and $\CCCi \SSSS \subseteq \SSSSs$ holds. Let $\gg$ be a non-invertible element of~$\CCC$. Let $(\gg_1 \wdots \gg_\pp)$ be a decomposition of~$\gg$ such that $\gg_\ii$ lies in~$\SSSS \cup \CCCi$ for every~$\ii$. As $\gg$ is not invertible, there exists~$\ii$ such that $\gg_\ii$ is not invertible. Assume that $\ii$ is minimal with this property. Then $\gg_1 \pdots \gg_{\ii-1}$ is invertible and, as $\CCCi \SSSS \subseteq \SSSSs$ holds, there exists~$\gg'$ in~$\SSSS \setminus \CCCi$ and~$\ie$ in~$\CCCi$ satisfying $\gg_1 \pdots \gg_\ii = \gg' \ie$. Then $\gg'$ is a non-invertible element of~$\SSSS$ left-dividing~$\gg$. Hence $\SSSS$ left-generates~$\CCC$. 
\end{solu}

%%%
\begin{exer}{20}{equivalence}{Z:PREquiv}
Assume that $(\SSSS, \RRR)$ is a category presentation. Say that an element~$\sss$ of~$\SSSS$ is $\RRR$-right-invertible if $\sss \ww \eqp_\RRR \ew_\xx$ ($\xx$ the source of~$\sss$) holds for some~$\ww$ in~$\SSSS^*$. Show that a category presentation $(\SSSS, \RRR)$ is complete with respect to right-reversing if and only if, for all~$\uu, \vv$ in~$\SSSS^*$, the following are equivalent: \ITEM1 $\uu$ and~$\vv$ are $\RRR$-equivalent (that is, $\uu \eqp_\RRR \vv$ holds), \ITEM2 $\INV\uu \vv \rev_\RRR \vv' \INV{\uu'}$ holds for some $\RRR$-equivalent paths~$\uu', \vv'$ in~$\SSSS^*$ all of which entries are $\RRR$-right-invertible.
\end{exer}

\begin{solu}
Assume that right-reversing is complete for~$(\SSSS, \RRR)$ and $\uu \eqp_\RRR \vv$ holds. Denoting by~$\yy$ the common target of~$\uu$ and~$\vv$, we have $\uu \sepp \eps_\yy \eqp_\RRR \vv \sepp \eps_\yy$. Hence, by definition of completeness, there exist~$\uu', \vv', \ww$ satisfying $\INV\uu \sepp \vv \rev_\RRR \vv' \sepp \INV{\uu'}$, $\eps_\yy \eqp_\RRR\uu' \sepp \ww$, and $\eps_\yy \eqp_\RRR\uu' \sepp \ww$. Hence, in~$\CCC$, we have $\clp{\uu'} \, \clp{\ww} = \id\yy = \clp{\vv'} \, \clp{\ww}$, whence $\clp{\uu'} = \clp{\vv'} = (\clp{\ww})\inv$. So $\uu'$ and $\vv'$ are $\RRR$-equivalent and, as their classes are invertible, they must consist of invertible entries.

Conversely, assume that the condition of~\ITEM2 is satisfied, and that $\uu \sepp \vv'$ and $\vv \sepp \uu'$ are $\RRR$-equivalent. By~\ITEM2, there exist $\RRR$-equivalent $\SSSS$-paths~$\uu_0, \vv_0$ such that $\INV{(\uu \sepp \vv')} \sepp (\vv \sepp \uu')$ is right-reversible to~$\vv_0 \sepp \INV{\uu_0}$ and all entries in~$\uu_0$ and~$\vv_0$ are invertible. By Lemma~II.4.23, the reversing of $\INV{(\uu \sepp \vv')} \sepp (\vv \sepp \uu')$ to~$\vv_0 \sepp \INV{\uu_0}$ splits into four reversings. By construction, all entries in $\uu_1, \uu_2, \vv_1, \vv_2$ are invertible. Put $\ww = \uu''' \sepp \vv_2 \sepp \INV{\uu_2} \sepp \INV{\uu_1}$. By assumption, we have $\vv_1 \sepp \vv_2 \eqp_\RRR \uu_1 \sepp \uu_2$, whence $\INV{\vv_2} \sepp \INV{\vv_1} \eqp_\RRR \INV{\uu_2} \sepp \INV{\uu_1}$. We deduce
$$\uu' \eqp \uu'' \sepp \uu''' \sepp \vv_2 \sepp \INV{\uu_2} \sepp \INV{\uu_1} = \uu'' \sepp \ww, \mbox{\ and\ } \vv' \eqp \vv'' \sepp \uu''' \sepp \vv_2 \sepp \INV{\vv_2} \sepp \INV{\vv_1} \eqp \uu'' \sepp \ww,$$
which means that $(\uu, \vv, \uu', \vv')$ factorizing through right-reversing. Hence right-reversing is complete for~$(\SSSS, \RRR)$.
\end{solu}

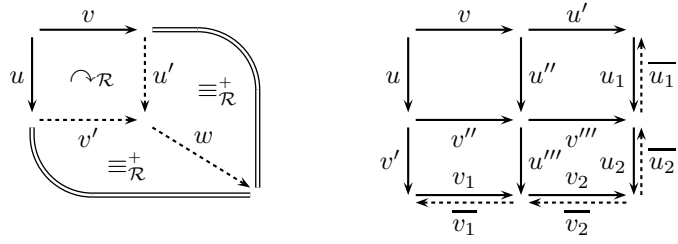
\begin{figure}[htb]\centering
\begin{picture}(80,25)(0,-1)
\pcline{->}(1,22)(14,22)
\taput{$\vv$}
\pcline{->}(0,21)(0,11)
\tlput{$\uu$}
\psline[style=eps](16,22)(22,22)
\psarc[style=eps](22,14){8}{0}{90}
\pcline[style=eps](30,14)(30,1)
\psline[style=eps](0,9)(0,8)
\psarc[style=eps](8,8){8}{180}{270}
\pcline[style=eps](8,0)(29,0)
\pcline[style=exist]{->}(1,10)(14,10)
\tbput{$\vv'$}
\pcline[style=exist]{->}(15,21)(15,11)
\trput{$\uu'$}
\pcline[style=exist]{->}(16,9)(29,1)
\put(21.5,7){$\ww$}
\put(5,15){$\rev_\RRR$}
\put(10,3){$\eqp_\RRR$}
\put(22,13){$\eqp_\RRR$}

\pcline{->}(51,22)(64,22)
\taput{$\vv$}
\pcline{->}(66,22)(79,22)
\taput{$\uu'$}
\pcline{->}(50,21)(50,11)
\tlput{$\uu$}
\pcline{->}(65,21)(65,11)
\trput{$\uu''$}
\pcline{->}(80,21)(80,11)
\tlput{$\uu_1$}
\pcline[style=exist]{<-}(81,21)(81,11)
\trput{$\INV{\uu_1}$}
\pcline{->}(51,10)(64,10)
\tbput{$\vv''$}
\pcline{->}(66,10)(79,10)
\tbput{$\ \vv'''$}
\pcline{->}(50,9)(50,0)
\tlput{$\vv'$}
\pcline{->}(65,9)(65,0)
\trput{$\uu'''$}
\pcline{->}(80,9)(80,0)
\tlput{$\uu_2$}
\pcline[style=exist]{<-}(81,9)(81,0)
\trput{$\INV{\uu_2}$}
\pcline{->}(51,0)(64,0)
\taput{$\vv_1$}
\pcline{->}(66,0)(79,0)
\taput{$\vv_2$}
\pcline[style=exist]{<-}(51,-1)(64,-1)
\tbput{$\INV{\vv_1}$}
\pcline[style=exist]{<-}(66,-1)(79,-1)
\tbput{$\INV{\vv_2}$}
\end{picture}
\caption[]{\sf Solution to Exercise~20}
\label{F:PREquiv}
\end{figure}

%%%%
\begin{exer}{22}{complete \vs cube, complemented case}{Z:PRComplementCube}
Assume that $(\SSSS, \RRR)$ is a right-complemented presentation, associated with the syntactic right-complem\-ent~$\RC$. Show the equivalence of the following three properties: \ITEM1 Right-reversing is complete for~$(\SSSS, \RRR)$; \ITEM2 The map~$\RC^*$ is compatible with $\eqpRR$-equivalence in the sense that, if $\uu' \eqpRR \uu$ and $\vv' \eqpRR \vv$ hold, then $\RC^*(\uu', \vv')$ is defined if and only if $\RC(\uu, \vv)$ is and, if so, they are $\eqpRR$-equivalent; \ITEM3 The $\RC$-cube condition is true for every triple of $\SSSS$-paths. [Hint: For \ITEM1 $\Rightarrow$ \ITEM2, note that, if $\RC^*(\uu, \vv)$ is defined and $\uu' \eqpRR \uu$ and $\vv' \eqpRR \vv$ hold, then $(\uu', \vv,' \RC^*(\uu, \vv), \RC^*(\vv, \uu))$ must be $\rev$-factorable; for \ITEM2 $\Rightarrow$ \ITEM3, compute $\RC^*(\ww, \uu\RC^*(\uu, \vv))$ and $\RC^*(\ww, \vv\RC^*(\vv, \uu))$; for \ITEM3 $\Rightarrow$ \ITEM1, use Lemma~II.4.61.]
\end{exer}

\begin{solu}
\ITEM1$\Rightarrow$\ITEM2 Assume that $\RC^*(\uu, \vv)$ exists and $\uu' \eqpRR \uu$ and $\vv' \eqpRR \vv$ hold. By Corollary~II.4.36, we have $\uu \RC^*(\uu, \vv) \eqpRR \vv \RC^*(\vv, \uu)$, whence $\uu' \RC^*(\uu, \vv) \eqpRR \vv' \RC^*(\vv, \uu)$. By assumption, the quadruple $(\uu', \vv', \RC^*(\uu, \vv), \RC^*(\vv, \uu)$ is $\rev$-factorable: this means that $\RC^*(\uu', \vv')$ and $\RC^*(\vv', \uu')$ are defined and there exists~$\ww$ satisfying 
$$\RC^*(\uu, \vv) \eqpRR \RC^*(\uu', \vv') \ww \mbox{\quad and\quad} \RC^*(\vv, \uu) \eqpRR \RC^*(\vv', \uu') \ww.$$
By a symmetric argument, we obtain the existence of~$\ww'$ satisfying 
$$\RC^*(\uu', \vv') \eqpRR \RC^*(\uu, \vv) \ww'\mbox{\quad and\quad} \RC^*(\vv', \uu') \eqpRR \RC^*(\vv, \uu) \ww'.$$ Merging, we deduce $\RC^*(\uu, \vv) \eqpRR \RC^*(\uu, \vv) \ww' \ww$. By Corollary~II.4.45, the category~$\PRESp\SSSS\RRR$ is left-cancellative, so we deduce $\ew_{\ud} \eqpRR \ww' \ww$, whence $\ww' = \ww = \ew_{\ud}$ and, finally, $\RC^*(\vv', \uu') \eqpRR \RC^*(\uu, \vv)$.

\ITEM2$\Rightarrow$\ITEM3 Assume that $\RC^*$ is compatible with~$\eqpRR$ and $\RCt(\uu, \vv, \ww)$ is defined. By Lemma~II.4.6, we have $\RCt(\uu, \vv, \ww) = \RC^*(\uu (RC^*(\uu, \vv), \ww)$. As $\uu (\RC^*(\uu, \vv) \eqpRR \vv \RC^*(\vv, \uu)$ holds, the compatibility assumption implies that $\RC^*(\vv (RC^*(\vv, \uu), \ww)$ is defined as well, and it is $\eqpRR$-equivalent to the latter. Now, by Lemma~II.4.6 again, the path $\RC^*(\vv (\RC^*(\vv, \uu), \ww)$ is $\RCt(\vv, \uu, \ww)$, so we conclude that $\RCt(\vv, \uu, \ww)$ is defined and it is $\eqpRR$-equivalent to~$\RCt(\uu, \vv, \ww)$, that is, the $\RC$-cube condition is true for~$(\uu, \vv, \ww)$.

\ITEM3$\Rightarrow$\ITEM1 Lemma~II.4.61 gives the result directly. 
\end{solu}

%%%%
\begin{exer}{23}{alternative proof}{Z:PRAlternative}
Assume that $(\SSSS, \RRR)$ is a right-complemented presentation associated with the syntactic right-complement~$\RC$, that $(\SSSS, \RRR)$ is right-Noetherian, and that the $\RC$-cube condition is true on~$\SSSS$. \ITEM1 Show that, for all~$\rr, \sss, \tt$ in~$\SSSS$, the path $\RC^*(\rr, \sss \RC(\sss, \tt))$ is defined if and only if $\RC^*(\rr, \tt \RC(\tt, \sss))$ is and, in this case, the relations $\RC^*(\rr, \sss \RC(\sss, \tt)) \eqpRR \RC^*(\rr, \tt \RC(\tt, \sss))$ and $\RC^*(\sss \RC(\sss, \tt), \rr) \eqpRR \RC^*(\tt \RC(\tt, \sss), \rr)$ are satisfied.
\ITEM2 Show that the map~$\RC^*$ is compatible with~$\eqpRR$, that is, the conjunction of $\uu' \eqpRR \uu$ and $\vv' \eqpRR \vv$ implies that $\RC^*(\uu', \vv')$ exists if and only if $\RC^*(\uu, \vv)$ does and, in this case, they are $\eqpRR$-equivalent. [Hint: Show using on~$\wit^*(\uu \RC^*(\uu, \vv))$, where $\wit^*$ is a right-Noetherianity witness for~$(\SSSS, \RRR)$ that, if $\RC^*(\uu, \vv)$ is defined and we have $\uu' \eqpRR \uu$ and $\vv' \eqpRR \vv$, then $\RC^*(\uu', \vv')$ is defined and we have $\RC^*(\uu', \vv') \eqpRR \RC^*(\uu, \vv)$ and $\RC^*(\vv', \uu') \eqpRR \RC^*(\vv, \uu)$.]
\ITEM3 Apply Exercise~22 to deduce a new proof of Proposition~II.4.16 in the right-Noetherian case.
\end{exer}

\begin{solu}
\ITEM1 Lemma~II.4.6 gives
\begin{gather*}
\RC^*(\sss \RC(\sss, \tt), \rr) = \RC^*(\RC(\sss, \tt), \RC(\sss, \rr)) = \RCt(\sss, \tt, \rr),\\
\RC^*(\rr, \sss \RC(\sss, \tt)) = \RC(\rr, \sss) \RC^*(\RC(\sss, \rr), \RC(\sss, \tt)) = \RC(\rr, \sss) \RCt(\sss, \rr, \tt).
\end{gather*}
As the $\RC$-cube condition is true on~$\{\rr, \sss, \tt\}$, we have $\RCt(\sss, \tt, \rr) \eqpRR \RCt(\tt, \rr, \rr)$, so, using the first equality above and its counterpart exchanging~$\sss$ and~$\tt$, we find
$$\RC^*(\sss \RC(\sss, \tt), \rr) = \RCt(\sss, \tt, \rr) \eqpRR \RCt(\tt, \rr, \rr) = \RC^*(\tt \RC(\tt, \sss), \rr).$$
Similarly, the $\RC$-cube condition gives the relations $\RCt(\sss, \rr, \tt) \eqpRR \RCt(\rr, \sss, \tt)$ and $\RCt(\rr, \tt, \sss) \eqpRR \RCt(\tt, \rr, \sss)$, whence, by Corollary~II.4.36, 
\begin{multline*} 
\RC^*(\rr, \sss \RC(\sss, \tt)) = \RC(\rr, \sss) \RCt(\sss, \rr, \tt) \eqpRR \RC(\rr, \sss) \RCt(\rr, \sss, \tt)\\ \eqpRR \RC(\rr, \tt) \RCt(\rr, \tt, \sss) \eqpRR \RC(\rr, \tt) \RCt(\tt, \rr, \sss) = \RC^*(\rr, \tt \RC(\tt, \sss)).
\end{multline*} 

\ITEM2 The result of~\ITEM1 is the compatiblity in the case $\uu' = \uu = \rr$ and $\vv = \sss \RC(\sss, \tt)$, $\vv' = \tt \RC(\tt, \sss)$, that is, in the basic case of equivalence. We establish the general result using induction on $\wit^*(\uu \RC^*(\uu, \vv))$ and, for a given value~$\alpha$, on the sum~$\dd$ of the combinatorial distances from~$\uu$ to~$\uu'$ and from~$\vv$ to~$\vv'$. For an induction, it is sufficient to consider the case $\dd = 1$, that is, we may assume $\uu' = \uu$ and $\dist(\vv, \vv') = 1$, that is, there exist $\sss, \tt$ in~$\SSSS$ and $\vv_0$, $\vv_1$ in~$\SSSS^*$ satisfying $\vv = \vv_0 \sss \RC(\sss, \tt) \vv_1$ and $\vv'  = \vv_0 \tt \RC(\tt, \sss) \vv_1$. We assume that $\RC^*(\uu, \vv)$ is defined, and our aim is to show that $\RC^*(\uu', \vv')$ is defined as well and we have $\RC^*(\uu', \vv') \eqpRR \RC^*(\uu, \vv)$ and $\RC^*(\vv', \uu') \eqpRR \RC^*(\vv, \uu)$. To this end, we compare the reversing grids below:
$$
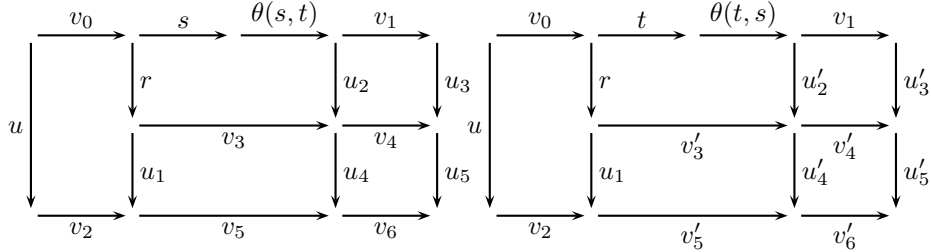
\begin{picture}(56,32)(0,-3)
\psset{xunit=0.9mm}
\pcline{->}(1,24)(14,24)\taput{$\vv_0$}
\pcline{->}(16,24)(29,24)\taput{$\sss$}
\pcline{->}(31,24)(44,24)\taput{$\RC(\sss, \tt)$}
\pcline{->}(46,24)(59,24)\taput{$\vv_1$}
\pcline{->}(16,12)(44,12)\tbput{$\vv_3$}
\pcline{->}(46,12)(59,12)\tbput{$\vv_4$}
\pcline{->}(1,0)(14,0)\tbput{$\vv_2$}
\pcline{->}(16,0)(44,0)\tbput{$\vv_5$}
\pcline{->}(46,0)(59,0)\tbput{$\vv_6$}
\pcline{->}(0,23)(0,1)\tlput{$\uu$}
\pcline{->}(15,23)(15,13)\trput{$\rr$}
\pcline{->}(15,11)(15,1)\trput{$\uu_1$}
\pcline{->}(45,23)(45,13)\trput{$\uu_2$}
\pcline{->}(45,11)(45,1)\trput{$\uu_4$}
\pcline{->}(60,23)(60,13)\trput{$\uu_3$}
\pcline{->}(60,11)(60,1)\trput{$\uu_5$}
\end{picture}
\hspace{5mm}
\begin{picture}(56,32)(0,-3)
\psset{xunit=0.9mm}
\pcline{->}(1,24)(14,24)\taput{$\vv_0$}
\pcline{->}(16,24)(29,24)\taput{$\tt$}
\pcline{->}(31,24)(44,24)\taput{$\RC(\tt, \sss)$}
\pcline{->}(46,24)(59,24)\taput{$\vv_1$}
\pcline{->}(16,12)(44,12)\tbput{$\vv'_3$}
\pcline{->}(46,12)(59,12)\tbput{$\vv'_4$}
\pcline{->}(1,0)(14,0)\tbput{$\vv_2$}
\pcline{->}(16,0)(44,0)\tbput{$\vv'_5$}
\pcline{->}(46,0)(59,0)\tbput{$\vv'_6$}
\pcline{->}(0,23)(0,1)\tlput{$\uu$}
\pcline{->}(15,23)(15,13)\trput{$\rr$}
\pcline{->}(15,11)(15,1)\trput{$\uu_1$}
\pcline{->}(45,23)(45,13)\trput{$\uu'_2$}
\pcline{->}(45,11)(45,1)\trput{$\uu'_4$}
\pcline{->}(60,23)(60,13)\trput{$\uu'_3$}
\pcline{->}(60,11)(60,1)\trput{$\uu'_5$}
\end{picture}$$
By assumption, the left grid exists, and we wish to show that the right grid exists as well and that the corresponding paths are pairwise $\eqpRR$-equivalent. The rectangles on the left ($\uu$ and $\vv_0$) coincide. Next, we find a rectangle as in~\ITEM1, namely $\rr$ and $\sss \RC(\sss, \tt)$ \vs $\rr$ and $\tt \RC(\tt, \sss)$. By~\ITEM1, $\uu'_2$ and $\vv'_3$ exist and we have $\uu'_2 \eqpRR \uu_2$ and $\vv'_3 \eqpRR \vv_3$. Then consider the median bottom rectangles ($\uu_1$ and $\vv'_3$): the point is the inequality
$$\wit^*(\RC^*(\uu_1 \vv_5)) < \wit^*(\rr \uu_1 \vv_5) \le \wit^*(\rr \uu_1 \vv_5 \vv_6) \le \wit^* (\vv_0 \rr \uu_1 \vv_5 \vv_6) = \wit^*(\uu, \vv) \le \alpha.$$ 
As $\RC^*(\uu_1, \vv_3)$ exists and we have $\vv'_3 \eqpRR \vv_3$, the induction hypothesis implies that $\RC^*(\uu_1, \vv'_3)$ exists as well and gives $\uu'_4 \eqpRR \uu_4$ and $\vv'_5 \eqpRR \vv_5$. The argument is the same for the two right squares.

\ITEM3 So, if the $\RC$-cube condition is true on~$\SSSS$, the map~$\RC^*$ is compatible with~$\eqpRR$. By Exercise~22, the latter condition implies that right-reversing is complete for~$(\RRR, \SSSS)$.
\end{solu}

%%%
\begin{figure}[htb]\centering
\begin{picture}(90,27)(0,-2)
\pcline{->}(1,0)(29,0)
\tbput{$\vvt$}
\pcline{->}(1,12)(14,12)
\taput{$\ww$}
\pcline{->}(16,24)(29,24)
\taput{$\vv$}
\pcline{->}(0,11)(0,1)
\tlput{$\uu$}
\pcline{->}(15,23)(15,13)
\tlput{$\ww$}
\pcline{->}(30,23)(30,1)
\trput{$\uut$}
\put(18,8){$\rev$}
\put(40,11){implies}
\pcline{->}(60,23)(60,13)
\tlput{$\uu$}
\pcline{->}(61,24)(74,24)
\taput{$\vv$}
\pcline[style=exist]{->}(75,23)(75,13)
\tlput{$\uu'$}
\pcline[style=exist]{->}(61,12)(74,12)
\taput{$\vv'$}
\pcline[style=exist]{->}(76,11)(89,1)
\put(81,9){$\ww'$}
\pcline(76,24)(83,24)
\psarc(83,17){7}{0}{90}
\pcline{->}(90,17)(90,1)
\put(91,13){$\uut$}
\pcline(60,11)(60,7)
\psarc(67,7){7}{180}{270}
\pcline{->}(67,0)(89,0)
\put(75,-3){$\vvt$}
\put(66,18){$\rev$}
\put(70,6){$\eqp$}
\put(81,16){$\eqp$}
\end{picture}
\caption[]{\sf\smaller Solution to Exercise~24: The cube condition viewed as a factorization property: whenever $\INV\uu \ww \INV\ww \vv$ right-reverses to~$\vvt \INV\uut$, the quadruple $(\uu, \vv, \uut, \vvt)$  is $\rev$-factorable.}
\label{F:PRCubeP}
\end{figure}

\begin{exer}{24}{cube condition}{Z:PRComplCube}
Assume that $(\SSSS, \RRR)$ is a category presentation \ITEM1 Show that the cube condition is true for $(\uu, \vv, \ww)$ if and only if, for all~$\uut, \vvt$ in~$\SSSS^*$ satisfying $\INV\uu \ww \INV\ww \vv \rev_\RRR \vvt \INV{\uut}$, the quadruple~$(\uu, \vv, \uut, \vvt)$  is $\rev$-factorable.  \ITEM2 Show that, if right-reversing is complete for~$(\SSSS, \RRR)$, then the cube condition is true for every triple of $\SSSS$-paths.
\end{exer}

\begin{solu}
\ITEM2 Assume that $\uu, \vv, \ww$ belong to~$\SSSS^*$, and we have $\INV\uu \ww \rev_\RRR \vv_1 \INV{\uu_0}$, $\INV\ww \vv \rev_\RRR \vv_0 \INV{\uu_1}$, and $\INV{\uu_0} \vv_0 \rev_\RRR \vv_2 \INV{\uu_2}$ (so that $\uu, \vv, \ww$ necessarily share the same source). By Proposition~II.4.34, we have $\uu \vv_1 \eqp_\RRR \ww \uu_0$, $\vv \uu_1 \eqp_\RRR \ww \vv_0$, and $\uu_0 \vv_2 \eqp_\RRR \vv_0 \uu_2$, and we deduce $\uu \vv_1 \vv_2 \eqp_\RRR \ww \uu_0 \vv_2 \eqp_\RRR \ww \vv_0 \uu_2 \eqp_\RRR \vv \uu_1 \uu_2$. The assumption that right-reversing is complete then implies that $(\uu, \vv, \uu_1\uu_2, \vv_1\vv_2)$  is $\rev$-factorable, which exactly means that the cube condition is true for~$(\uu, \vv, \ww)$.
\end{solu}

%%%%
\chapter{Chapter~III: Normal decompositions}

%%%%
\section*{Skipped proofs}

\begin{ADprop}{III.1.14}{power}{P:NFPowerGrouping}
If $\SSSS$ is a subfamily of a left-cancellative category~$\CCC$ and $\seqqq{\gg_1}\etc{\gg_\pp}$ is $\SSSS$-greedy, then $\seqq{\gg_1 \pdots \gg_\mm}{\gg_{\mm+1} \pdots \gg_\pp}$ is $\Pow\SSSS\mm$-greedy for $1 \le \mm \le \pp$, that is, 
\begin{equation}
\ADlabel{III.1.15}
\mbox{Each relation $\sss \dive \ff \gg_1 \pdots \gg_\pp$ with $\sss$ in~$\Pow\SSSS\mm$ implies $\sss \dive \ff\gg_1 \pdots \gg_\mm$.}
\end{equation}
\end{ADprop}

\begin{figure}[htb]\centering
\begin{picture}(105,15)(0,0)
\pcline{->}(1,12)(15,12)
\taput{$\sss_1$}
\pcline{->}(16,12)(29,12)
\taput{$\sss_2$}
\pcline[style=etc](31,12)(44,12)
\pcline{->}(46,12)(59,12)
\taput{$\sss_\mm$}
\pcline{->}(1,0)(14,0)
\tbput{$\gg_1$}
\pcline{->}(16,0)(29,0)
\tbput{$\gg_2$}
\pcline[style=etc](31,0)(44,0)
\pcline{->}(46,0)(59,0)
\tbput{$\gg_\mm$}
\pcline{->}(61,0)(74,0)
\tbput{$\gg_{\mm+1}$}
\pcline[style=etc](76,0)(89,0)
\pcline{->}(91,0)(104,0)
\tbput{$\gg_\pp$}
\pcline{->}(0,11)(0,1)
\tlput{$\ff$}
\pcline[style=exist]{->}(15,11)(15,1)
\trput{$\ff_1$}
\pcline[style=exist]{->}(60,11)(60,1)
\psarc[style=thin](14.5,0){2.5}{180}{360}
\psarc[style=thin](59.5,0){2.5}{180}{360}
\psline(61,12)(98,12)
\psarc(98,5){7}{0}{90}
\psline{->}(105,5)(105,1)
\end{picture}
\caption[]{\sf\smaller Inductive proof of Proposition~III.1.14.}
\label{F:NFPower}
\end{figure}
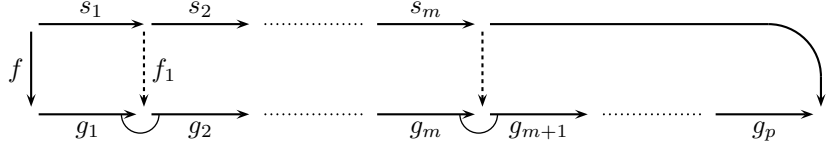

\begin{proof}
(See Figure~\ref{F:NFPower}.) We use induction on~$\mm$. For $\mm = 1$, the result follows from Proposition~III.1.12, and more precisely from~(III.1.13). Assume $\mm \ge 2$. Let $\sss \in \Pow\SSSS\mm$, say $\sss = \sss_1 \pdots \sss_\mm$ with $\sss_1 \wdots \sss_\mm$ in~$\SSSS$, and $\sss \dive \ff \gg_1 \pdots \gg_\pp$. Then we have $\sss_1 \dive \ff \gg_1 \pdots \gg_\pp$, hence, by~(III.1.13), $\sss_1 \dive \ff \gg_1$, say $\ff \gg_1 = \sss_1 \ff_1$, and, therefore, $\sss_2 \pdots \sss_\mm \dive \ff_1 \gg_2 \pdots \gg_\pp$. As $\sss_2 \pdots \sss_\mm$ belongs to~$\Pow\SSSS{\mm-1}$, the induction hypothesis implies $\sss_2 \pdots \sss_\mm \dive \ff_1 \gg_2 \pdots \gg_\mm$, whence $\sss_1 \sss_2 \pdots \sss_\mm \dive \ff \gg_1 \gg_2 \pdots \gg_\mm$, as expected.
\end{proof}

\begin{ADlemm}{III.2.52}{L:NFWeakLcm}
For $\ff, \gg, \ff', \gg'$ in a cancellative category~$\CCC$ satisfying $\ff' \gg = \gg' \ff$, the following conditions are equivalent:

\ITEM1 The elements $\ff'$ and~$\gg'$ are left-disjoint; 

\ITEM2 The element $\ff' \gg$ is a weak left-lcm of~$\ff$ and~$\gg$.
\end{ADlemm}

\begin{proof}
Assume~\ITEM1, and assume that $\ff'' \gg = \gg'' \ff$ is a common left-multiple of~$\ff$ and~$\gg$ such that $\ff'' \gg$ and~$\ff' \gg$ admit a common left-multiple, say $\hh' (\ff'' \gg) = \hh'' (\ff' \gg)$. By assumption we also have $\hh' (\gg'' \ff) = \hh'' (\gg' \ff)$ and, because $\CCC$ is assumed to be right-cancellative, we deduce $\hh' \ff'' = \hh'' \ff' $ and $\hh' \gg'' = \hh'' \gg'$. Thus $(\hh', \hh'')$ witnesses for $(\ff'', \gg'') \EQU (\ff', \gg')$. As $\ff'$ and~$\gg'$ are left-disjoint, we deduce $\hh' \dive \hh''$, that is, there exists~$\hh$ satisfying $\hh'' = \hh' \hh$. We deduce $\hh' (\ff'' \gg) = \hh'' (\ff' \gg) = \hh' \hh (\ff' \gg)$, whence $\ff'' \gg = \hh (\ff' \gg)$ by left-cancelling~$\hh'$. This shows that $\ff'' \gg$ is a left-multiple of $\ff' \gg$, and the latter is a weak left-lcm of~$\ff$ and~$\gg$. So \ITEM1 implies~\ITEM2.

Assume now \ITEM2, and assume that $(\hh', \hh'')$ witnesses for $(\ff'', \gg'') \EQU (\ff', \gg')$, that is, we have $\hh' \ff'' = \hh'' \ff'$ and $\hh' \gg'' = \hh'' \gg'$. We deduce $\hh' \ff'' \gg = \hh'' \ff' \gg = \hh'' \gg' \ff = \hh' \gg'' \ff$, whence $\ff'' \gg = \gg'' \ff$ by left-cancelling~$\hh'$. So $\ff'' \gg$ is a common left-multiple of~$\ff$ and~$\gg$. Moreover, the above equalities show that $\ff'' \gg$ and $\ff' \gg$ admit a common left-multiple. As $\ff' \gg$ is a weak left-lcm of~$\ff$ and~$\gg$, we deduce that $\ff'' \gg$ is a left-multiple of~$\ff' \gg$, that is, there exists~$\hh$ satisfying $\ff'' \gg = \hh \ff' \gg$, hence $\hh'' \ff' \gg = \hh' \ff'' \gg = \hh' \hh \ff' \gg$. Right-cancelling~$\ff' \gg$, we deduce $\hh'' = \hh' \hh$, that is, $\hh' \dive \hh''$. Hence $\ff'$ and~$\gg'$ are left-disjoint, and \ITEM2 implies~\ITEM1.
\end{proof}

\begin{ADprop}{III.2.53}{symmetric normal exist}{P:NFSymExistApp}
If $\SSSS$ is a Garside family in a cancellative category~$\CCC$, the following conditions are equivalent:

\ITEM1 For all~$\ff, \gg$ in~$\CCC$ admitting a common right-multiple, there exists a symmetric $\SSSS$-normal path $\INV\uu \sepp \vv$ satisfying $(\ff, \gg) \EQU (\clp\uu, \clp\vv)$;

\ITEM2 The category~$\CCC$ admits conditional weak left-lcms.
\end{ADprop}

\begin{proof}
Assume~\ITEM1, and let $\ff, \gg$ be two elements of~$\CCC$ that admit a common left-multiple, say $\fft \gg = \ggt \ff$. By~\ITEM1 applied to~$\fft$ and~$\ggt$, there exists a symmetric $\SSSS$-normal path $\INV\uu \sepp \vv$ such that, putting $\ff' = \clp\uu$ and $\gg' = \clp\vv$, we have $(\fft, \ggt) \EQU (\ff', \gg')$, that is, there exist~$\hh', \hht$ satisfying $\hh' \fft = \hht \ff'$ and $\hh' \ggt = \hht \gg'$. By Proposition~III.2.11, $\ff'$ and~$\gg'$ are left-disjoint, so, by Lemma~III.2.52, $\ff' \gg$, which is also $\gg' \ff$, is a weak left-lcm of~$\ff$ and~$\gg$. Hence $\CCC$ admits conditional weak left-lcms, and \ITEM1 implies~\ITEM2.

Conversely, assume~\ITEM2, and let $\ff, \gg$ be two elements of~$\CCC$ that admit a common right-multiple, say $\ff \gg' = \gg \ff'$. By~\ITEM2, there exists a weak left-lcm of~$\ff'$ and~$\gg'$, say $\ff'' \gg' = \gg'' \ff'$, and~$\hh$ satisfying $\ff \gg' = \hh \ff'' \gg'$, hence $\ff = \hh \ff''$ by right-cancelling~$\gg'$ and, similarly, $\gg = \hh \gg''$. By Lemma~III.2.52, the elements~$\ff''$ and~$\gg''$ are left-disjoint. Let $\uu$ be an $\SSSS$-normal decomposition of~$\ff''$ and $\vv$ be an $\SSSS$-normal decomposition of~$\gg''$. By (the trivial part of) Proposition~III.2.11, the first entries of~$\uu$ and~$\vv$ are left-disjoint since $\ff''$ and~$\gg''$ are, so $\INV\uu \sepp \vv$ is symmetric $\SSSS$-normal. Finally, the pair $(\hh, \id\xx)$ (with $\xx$ the source of~$\ff$ and~$\gg$) witnesses for $(\ff, \gg) \EQU (\ff'', \gg'')$. So \ITEM2 implies~\ITEM1.
\end{proof}

\begin{ADlemm}{III.2.55}{L:NFWeakLeftLcm}
A cancellative category~$\CCC$ admits conditional weak left-lcms if and only if $\CCC$ is a strong Garside family in itself.
\end{ADlemm}

\begin{proof}
Assume that $\CCC$ is strong and $\ff \sss = \gg \tt$ holds. Then there exist $\ff', \gg', \hh$ such that $\ff'$ and $\gg'$ are left-disjoint and we have $\ff'\sss = \hh' \tt$, $\ff = \hh \ff'$, and $\gg = \hh \gg'$. By Lemma~III.2.52, $\ff'\gg$ is a weak left-lcm of~$\ff$ and~$\gg$, of which $\ff\gg$ is a left-multiple. So $\CCC$ admits conditional weak left-lcms.

Conversely, assume that $\CCC$ admits conditional weak left-lcms, and $\ff \sss = \gg \tt$ holds. Then there exists a weak left-lcm $\ff'\sss$ of~$\sss$ and~$\tt$ of which $\ff\sss$ is a left-multiple. By Lemma~III.2.52 again, $\ff'$ and~$\gg'$ are left-disjoint, and the condition of Definition~III.2.54 is satisfied. So $\CCC$ is strong.
\end{proof}

\begin{ADprop}{III.2.56}{symmetric normal, short case III}{P:NFSymShort3App}
If $\SSSS$ is a strong Garside family in a cancellative category~$\CCC$ admitting conditional weak  left-lcms, Algorithm~III.2.42 running on a positive--negative $\SSSSs$-path~$\vv \sepp \INV\uu$ such that $\clp\uu$ and~$\clp\vv$ admit a common left-multiple, say $\fft \clp\vv = \ggt \clp\uu$, returns a symmetric $\SSSS$-normal path~$\INV{\uu''} \sepp \vv''$ satisfying $(\ff, \gg) \EQU (\clp{\uu''}, \clp{\vv''})$ and $\clp{\uu'' \sepp \vv} = \clp{\vv'' \sepp \uu}$; moreover there exists~$\hh$ satisfying $\ff = \hh \clp{\uu''}$ and $\gg = \hh \clp{\vv''}$.
\end{ADprop}

\medskip Once again, Proposition~III.2.56 reduces to Proposition~III.2.44 in the case of a left-Ore category as it then says that $\INV\uu \sepp \vv$ is a decomposition of~$\cl{\vv \sepp \INV\uu}$ in~$\Env\CCC$.

\begin{figure}[htb]\centering
\begin{picture}(45,62)(0,0)
\pcline[style=double](46,48)(59,48)
\pcline[style=exist]{->}(16,36)(29,36)
\pcline{->}(31,36)(44,36)
\pcline{->}(16,24)(29,24)
\pcline{->}(31,24)(44,24)
\pcline[style=etc](2,12)(13,12)
\pcline{->}(16,12)(29,12)
\pcline{->}(31,12)(44,12)
\pcline[style=double](46,12)(59,12)
\pcline{->}(1,0)(14,0)
\pcline[style=etc](17,0)(28,0)
\pcline{->}(31,0)(44,0)
\pcline[style=double](46,0)(53,0)
\psarc[style=double](53,7){7}{270}{360}
\pcline[style=double](0,11)(0,1)
\pcline[style=exist]{->}(15,35)(15,25)
\pcline{->}(15,23)(15,13)
\pcline{->}(30,35)(30,25)
\pcline{->}(30,23)(30,13)
\pcline[style=etc](45,46)(45,37)
\pcline{->}(45,35)(45,25)
\pcline{->}(45,23)(45,13)
\pcline[style=double](45,11)(45,1)
\pcline{->}(60,47)(60,37)
\pcline[style=etc](60,34)(60,26)
\pcline{->}(60,23)(60,13)
\pcline[style=double](60,11)(60,7)
\put(13,5){Algorithm~III.1.52}
\put(51,20){\rotatebox{90}{\hbox{Algorithm~III.1.52}}}
\psarc[style=thinexist](15,36){3}{270}{360}
\psarc[style=thin](30,36){3}{270}{360}
\psarc[style=thin](15,24){3}{270}{360}
\psarc[style=thin](30,24){3}{270}{360}
\psarc[style=thin](14.5,12){2.5}{180}{360}
\psarc[style=thin](29.5,12){2.5}{180}{360}
\psarc[style=thin](45,36.5){2.5}{-90}{90}
\psarc[style=thin](45,24.5){2.5}{-90}{90}
\pcline(-14,60)(53,60)
\psarc(53,53){7}{0}{90}
\psline{->}(60,53)(60,49)
\pcline(-15,59)(-15,7)
\psarc(-8,7){7}{180}{270}
\psline{->}(-8,0)(-1,0)
\pscircle(-15,60){0.5}
\put(-18,59.5){$\xx$}
\psarc(38,53){7}{0}{90}
\psline{->}(45,53)(45,49)
\psarc(23,53){7}{0}{90}
\psline{->}(30,53)(30,37)
\psarc(-8,19){7}{180}{270}
\psline{->}(-8,12)(-1,12)
\psarc(-8,31){7}{180}{270}
\psline{->}(-8,24)(14,24)
\psline[style=exist]{->}(-14,59)(14.5,36.5)
\end{picture}
\caption[]{\sf\smaller Proof of Proposition~III.2.56: when the rectangular grid of Figure~III.18 is constructed, there exists at each step an arrow connecting the top--left corner to each local North-West corner of the current diagram.}
\label{F:NFSymShort3App}
\end{figure}
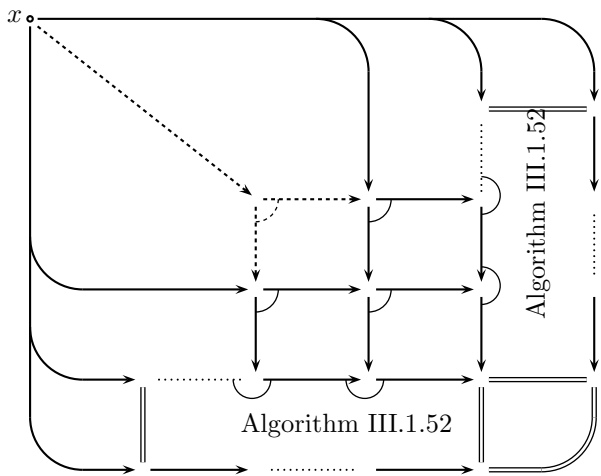

\begin{proof}
The argument is the same as for Proposition~III.2.44, the only difference being that, at each left-reversing step, one has to check the existence of a factorization of the initial equality $\ff \cl\vv = \gg \cl\uu$. The principle is explained in Figure~\ref{F:NFSymShort3App}: the induction hypothesis that is maintained at each step in the construction of the rectangular diagram is that, for every local North-West corner in the current diagram, there exists a factorizing arrow coming from the top-left object~$\xx$. When one more tile is added, the defining property of a strong Garside family guarantees that one can add a new tile in which the left and top arrows represent left-disjoint elements and there exists a factorizing arrow coming from~$\xx$. The rest of the proof is unchanged as, in particular, the third domino rule is still valid in the extended context.
\end{proof}

%%%%
\section*{Solution to selected exercises}

%%%%
\begin{exer}{28}{invertible}{Z:NFIsom}
Assume that $\CCC$ is a left-cancellative category and $\SSSS$ is included in~$\CCC$. Show that, if $\gg_1\pdots \gg_\pp$ belongs to~$\SSSSs$, then $\seqqq{\gg_1}\etc{\gg_\pp}$ being $\SSSS$-greedy implies that $\gg_2 \wdots \gg_\pp$ are invertible.
\end{exer}

\begin{solu}
The element $\gg_1 \pdots \gg_\pp$ lies in~$\SSSSs$ is equal to, hence left-divides, itself. By Proposition~III.1.12, $\gg_1 \sep \gg_2 \pdots \gg_\pp$ is $\SSSS$-greedy, so we deduce that $\gg_1 \pdots \gg_\pp$ left-divides~$\gg_1$, say $\gg_1 = \gg_1 \pdots \gg_\pp \gg'$. Left-cancelling~$\gg_1$, we deduce that $\gg_2 \pdots \gg_\pp \gg'$ is an identity-element, hence $\gg_2 \wdots \gg_\pp$, and $\gg'$ must be invertible.
\end{solu}

%%%%
\begin{exer}{29}{deformation}{Z:NFDeform}
Assume that $\CCC$ is a left-cancellative category. Show that a path $\seqqq{\gg_1}\etc{\gg_\qq}$ is a $\CCCi$-deformation of $\seqqq{\ff_1}\etc{\ff_\pp}$ if and only if $\gg_1 \pdots \gg_\ii \eqir \ff_1\pdots \ff_\ii$ holds for $1 \le \ii \le \max(\pp, \qq)$, the shorter path being extended by identity-elements if needed.
\end{exer}

\begin{solu}
Let $\rr = \max(\pp, \qq)$. Assume that $\ie_0 \wdots \ie_\rr$ are invertible elements witnessing that $\seqqq{\gg_1}\etc{\gg_\rr}$ is a $\CCCi$-deformation of $\seqqq{\ff_1}\etc{\ff_\rr}$. For every~$\ii$, we deduce $\ff_1 \pdots \ff_\ii \ie_\ii = \ie_0 \gg_1 \pdots \gg_\rr$, whence $\gg_1 \pdots \gg_\ii \eqir \ff_1\pdots \ff_\ii$ since $\ie_0$ is an identity-element. 

Conversely, assume $\gg_1 \pdots \gg_\ii \eqir \ff_1\pdots \ff_\ii$ for every~$\ii$, say $\gg_1 \pdots \gg_\ii = \ff_1\pdots \ff_\ii \ie_\ii$ with~$\ie_\ii$ in~$\CCC_\ii$. Set $\ie_0 = \id\xx$ where $\xx$ is the source of~$\ff_1$. Then we have $\ff_1 \ie_1 = \gg_1$ by construction. Assume $\ii \ge 2$. Then we obtain $\ff_1 \pdots \ff_{\ii-1} \ff_\ii \ie_\ii = \gg_1 \pdots \gg_{\ii-1} \gg_\ii = \ff_1 \pdots \ff_{\ii-1} \ie_{\ii-1} \gg_\ii$, whence $\ff_\ii \ie_\ii = \ie_{\ii-1} \gg_\ii$ by left-cancelling $\ff_1 \pdots \ff_{\ii-1}$. So $\seqqq{\gg_1}\etc{\gg_\rr}$ is a $\CCCi$-deformation of $\seqqq{\ff_1}\etc{\ff_\rr}$.
\end{solu}

%%%%
\begin{exer}{33}{left-disjoint}{Z:NFLeftDisjoint2}
Assume that $\CCC$ is a left-cancellative category, $\ff$ and $\gg$ are left-disjoint elements of~$\CCC$, and $\ff$ left-divides~$\gg$. Show that $\ff$ is invertible.
\end{exer}

\begin{solu}
Let $\xx$ be the common source of~$\ff$ and~$\gg$. By assumption we have $\ff \dive \id\xx \ff$ and $\ff \dive \id\xx \gg$. By definition of $\ff$ and $\gg$ being left-disjoint, this implies $\ff \dive \id\xx$, hence $\ff$ must be invertible.
\end{solu}

%%%%
\begin{exer}{34}{normal decomposition}{Z:NFSymNormalExist}
Give a direct argument from deriving Proposition~III.2.20 from Corollary~III.2.50 in the case when $\SSSS$ is strong.
\end{exer}

\begin{solu}
Let $\gg \ff\inv$ be an element of~$\CCC \CCC\inv$. Let $\seqqq{\sss_1}\etc{\sss_\pp}$ be an $\SSSS$-normal decomposition of~$\ff$, and let $\seqqq{\tt_1}\etc{\tt_\qq}$ be an $\SSSS$-normal decomposition of~$\gg$. We prove the existence of an $\SSSS$-normal decomposition for~$\gg \ff\inv$ using an induction on~$\pp$ to construct a rectangular diagram consisting of $\pp$ rows of $\qq$~tiles as in Lemma~III.2.31. As $\seqqq{\tt_1}\etc{\tt_\qq}$ is $\SSSS$-normal, Corollary~III.2.50 inductively implies that the elements of every horizontal line of the diagram make an $\SSSS$-normal path, and so do in particular the elements $\seqqq{\tt'_1}\etc{\tt'_{\qq}}$ of the top line. Similarly, as $\seqqq{\sss_1}\etc{\sss_\pp}$ is $\SSSS$-normal, Corollary~III.2.50 again inductively implies that the elements of every vertical line of the diagram make an $\SSSS$-normal path, and so do in particular the elements $\seqqq{\sss'_1}\etc{\sss'_\pp}$ of the left line. Finally, $\sss'_1$ and $\tt'_1$ are left-disjoint by construction. Hence $\seqqqqqq{\INV{\sss'_{\pp}}}\etc{\INV{\sss'_1}}{\tt'_1}\etc{\tt'_{\qq}}$ is an $\SSSS$-normal decomposition of~$\gg \ff\inv$.
\end{solu}

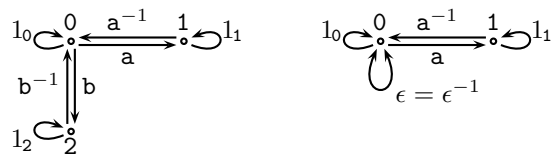
\begin{figure}[htb]
\begin{picture}(20,17)(0,0)
\pscircle(0,12){0.5}
\put(-1,13.5){$\mathtt0$}
\pscircle(15,12){0.5}
\put(14,13.5){$\mathtt1$}
\pscircle(0,0){0.5}
\put(-1,-3){$\mathtt2$}
\pcline{->}(1,11.5)(14,11.5)
\tbput{$\tta$}
\pcline{<-}(1,12.5)(14,12.5)
\taput{$\tta\inv$}
\pcline{->}(0.5,11)(0.5,1)
\trput{$\ttb$}
\pcline{<-}(-0.5,11)(-0.5,1)
\tlput{$\ttb\inv$}
\psbezier{->}(-1,-0.5)(-6,-2.5)(-6,2.5)(-1,0.5)
\put(-8,-2){$\id{\mathtt2}$}
\psbezier{->}(-1,11.5)(-6,9.5)(-6,14.5)(-1,12.5)
\put(-8,12.5){$\id{\mathtt0}$}
\psbezier{->}(16,11.5)(21,9.5)(21,14.5)(16,12.5)
\put(20,12.5){$\id{\mathtt1}$}
\end{picture}
\hspace{20mm}\begin{picture}(20,12)(0,0)
\pscircle(0,12){0.5}\put(-1,13.5){$\mathtt0$}
\pscircle(15,12){0.5}\put(14,13.5){$\mathtt1$}
\pcline{->}(1,11.5)(14,11.5)\tbput{$\tta$}
\pcline{<-}(1,12.5)(14,12.5)\taput{$\tta\inv$}
\psbezier{<->}(-0.5,11)(-3.5,4)(3.5,4)(0.5,11)\put(2,4){$\ie = \ie\inv$}
\psbezier{->}(-1,11.5)(-6,9.5)(-6,14.5)(-1,12.5)\put(-8,12.5){$\id{\mathtt0}$}
\psbezier{->}(16,11.5)(21,9.5)(21,14.5)(16,12.5)\put(20,12.5){$\id{\mathtt1}$}
\end{picture}
\caption[]{Diagrams of the categories of Exercise~35.}
\label{F:Base}
\end{figure}

%%%%
\begin{exer}{35}{Garside base}{Z:NFBase}
\ITEM1 Let $\GGG$ be the category whose diagram is displayed on Figure~\ref{F:Base} left, and let $\SSS = \{\tta, \ttb\}$. Show that $\GGG$ is a groupoid with nine elements, $\SSSS$ is a Garside base in~$\GGG$, the subcategory~$\CCC$ of~$\GGG$ generated by~$\SSSS$ contains no nontrivial invertible element, but $\CCC$ is not an Ore category. Conclusion? \ITEM2 Let \VR(4,0) $\GGG$ be the category whose diagram is displayed on Figure~\ref{F:Base} right, let $\SSSS = \{\ie, \tta\}$, and let $\CCC$ be the subcategory of~$\GGG$ generated by~$\SSSS$. Show that $\GGG$ is a groupoid and every element of~$\GGG$ admits a decomposition that is symmetric $\SSSS$-normal in~$\CCC$. Show that $\ie\tta$ admits a symmetric $\SSSS$-normal decomposition and no $\SSSS$-normal decomposition. Is $\SSSS$ a Garside family in~$\CCC$? Conclusion?
\end{exer}

\begin{solu}
\ITEM1 The nine elements of~$\GGG$ are $\id{\mathtt0}$, $\id{\mathtt1}$, $\id{\mathtt2 }$, $\tta$, $\ttb$, $\tta\inv$, $\ttb\inv$, $\tta\inv\ttb$, and $\ttb\inv\tta$. That $\SSSS$ is a Garside base in~$\GGG$ follows from a direct inspection. The subcategory of~$\GGG$ generated by~$\SSSS$ comprises five elements: $\id{\mathtt0}$, $\id{\mathtt1}$, $\id{\mathtt2 }$, $\tta$, and~$\ttb$, none of which is invertible. Next, $\CCC$ is not an Ore category since the elements~$\tta$ and~$\ttb$ have non common right-multiple although they share the same source. So, in Proposition~III.2.25, the conclusion cannot be strengthened to claim that $\CCC$ necessarily is an Ore category (whereas the assumption of Proposition~III.2.24 cannot be weakened to only assume that $\CCC$ is a left-Ore category).

\ITEM2 The elements of~$\GGG$ are $\id{\mathtt0}$, $\id{\mathtt1}$, $\ie$, $\tta$, $\ie \tta$, $\tta\inv$, and~$\tta\inv \ie$; those of~$\CCC$ are $\id{\mathtt0}$, $\id{\mathtt1}$, $\ie$, $\tta$, $\ie \tta$. Then $\ie\tta$ admits the symmetric $\SSSS$-normal decomposition $\seqq{\INV\ie}\tta$, but admits no $\SSSS$-normal decomposition. Hence $\SSSS$ is not a Garside family in~$\CCC$. So, in Proposition~III.2.25, we cannot simply drop the assumption about nontrivial invertible elements.
\end{solu}

%%%%
\chapter{Chapter~IV: Recognizing Garside families}

%%%%
\section*{Skipped proofs}

\begin{ADlemm}{IV.1.13}{L:GATransfer}
Assume that $\SSSS$ is a subfamily  of a left-cancellative category~$\CCC$. 

\ITEM2 The family~$\SSSS$ is closed under right-comultiple if and only if $\SSSSs$ is. 

\ITEM3 If $\SSSS$ is closed under right-complement and $\CCCi \SSSS \ince \SSSS$  holds, then $\SSSSs$  is closed under right-complement too. 
\end{ADlemm}

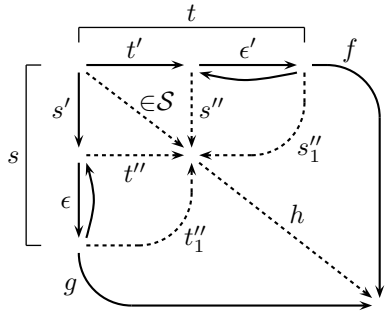
\begin{figure}[htb]\centering
\begin{picture}(47,38)(-7,0)
\pcline{->}(1,32)(14,32)\taput{$\tt'$}
\pcline{->}(16,32)(29,32)\taput{$\ie'$}
\pscurve{<-}(16,31)(22,30)(29,31)
\pcline(31,32)(33,32)
\psarc(33,25){7}{0}{90}
\pcline{->}(40,25)(40,1)\put(35,33){$\ff$}
\pcline{->}(0,31)(0,21)\tlput{$\sss'$}
\pcline{->}(0,19)(0,9)\tlput{$\ie$}
\pscurve{<-}(1,19)(2,14)(1,9)
\psarc(7,7){7}{180}{270}
\pcline{->}(7,0)(39,0)\put(-2,2){$\gg$}
\pcline[style=exist]{->}(15,31)(15,21)\trput{$\sss''$}
\pcline[style=exist]{->}(1,20)(14,20)\tbput{$\tt''$}
\pcline[style=exist]{->}(1,31)(14,21)\put(8,26){${\in}\SSSS$}
\pcline[style=exist]{->}(16,19)(39,1)\put(28,11){$\hh$}
\pcline[style=exist](1,8)(8,8)
\psarc[style=exist](8,15){7}{270}{360}
\psline[style=exist]{->}(15,15)(15,19)\put(14,8){$\tt''_1$}
\pcline[style=exist](30,31)(30,27)
\psarc[style=exist](23,27){7}{270}{360}
\psline[style=exist]{->}(23,20)(16,20)\put(29,20){$\sss''_1$}
\psline[style=thin](0,35)(0,37)
\psline[style=thin](30,35)(30,37)
\pcline[style=thin](0,37)(30,37)\taput{$\tt$}
\psline[style=thin](-5,32)(-7,32)
\psline[style=thin](-5,8)(-7,8)
\pcline[style=thin](-7,32)(-7,8)\tlput{$\sss$}
\end{picture}
\caption[]{\sf\smaller Proof of Lemma~IV.1.13.}
\label{F:GATransfer}
\end{figure}

\begin{proof}
\ITEM2 Assume that $\SSSSs$ is closed under right-comultiple, and we have $\sss  \gg = \tt \ff$ with~$\sss, \tt$ in~$\SSSS$. The assumption implies the  existence of~$\sss', \tt', \hh$ satisfying $\sss \tt' = \tt \sss' \in \SSSSs$,  $\ff = \sss' \hh$, and $\gg = \tt' \hh$. If $\sss \tt'$ is invertible, then  $\sss$ and~$\tt$ must be invertible too, so $\sss$ is a common  right-multiple of~$\sss$ and~$\tt$ lying in~$\SSSS$ of which $\sss \gg$ is a  right-multiple. Assume now that $\sss \tt'$ is not invertible. Then there  exists~$\ie$ in~$\CCCi$ such that $\sss \tt' \ie\inv$ lies in~$\SSSS$. Then we  have $\sss (\tt' \ie\inv) = \tt (\sss' \ie\inv) \in \SSSS$, $\ff = (\sss'  \ie\inv) (\ie \hh)$, and $\gg = (\tt' \ie\inv) (\ie \hh)$.  Hence  $\sss'  \ie\inv$, $\tt' \ie\inv$, and~$\ie \hh$ witness for $\SSSS$ being closed under  right-comultiple.

Conversely, assume that $\SSSS$ is closed under right-comultiple and we have  $\sss \gg = \tt\ff$ with $\sss, \tt \in \SSSSs$. Assume first that $\sss$  or~$\tt$ is invertible, say~$\sss$.  Put $\sss' = \id\yy$ where $\yy$ is the  target of~$\tt$, $\tt' = \sss\inv \tt$,  and $\hh = \ff$. Then we have $\sss \tt' = \tt \sss'$, $\ff = \sss' \hh$, $\gg  = \tt' \hh$, and $\sss'$ and~$\tt'$ are invertible, hence lie in~$\SSSSs$, and,  by assumption, so does~$\sss\tt'$, which is~$\tt$. So $\sss\tt'$ is a common right-multiple of~$\sss$ and~$\tt$ that lies  in~$\SSSSs$ and of which $\sss \gg$ is a right-multiple.  Assume now that neither~$\sss$ nor~$\tt$ is invertible. Then $\sss$ and~$\tt$  lie in~$\SSSS \CCCi$, and there exist~$\sss', \tt'$ in~$\SSSS$ and~$\ie, 
\ie'$ in~$\CCCi$ satisfying $\sss = \sss' \ie$ and~$\tt = \tt' \ie'$,  see  Figure~\ref{F:GATransfer}. Then we have $\sss' (\ie \gg) = \tt' (\ie' \ff)$ with $\sss', \tt'$ in~$\SSSS$. As $\SSSS$ is closed under  right-comultiple, there must exist~$\sss'', \tt''$, and~$\hh$ satisfying
$$\sss' \tt'' = \tt' \sss'' \in \SSSS, \quad \ie' \ff = \sss'' \hh,  \quad\mbox{and}\quad \ie \gg = \tt'' \hh.$$
As $\ie$ and~$\ie'$ are invertible, we can put $\sss''_1 = \ie'{}\inv \sss''$  and $\tt''_1 = \ie\inv \tt''$. Then we have $\ff = \sss''_1 \hh$ and $\gg = \tt''_1  \hh$, and $\sss \tt''_1 = \sss' \tt'' = \tt' \sss'' = \tt \sss''_1 \in \SSSS \ince  \SSSSs$. So, again, $\sss\tt''_1$ is a common right-multiple of~$\sss$ and~$\tt$  that lies in~$\SSSSs$ and of which $\sss \gg$ is a right-multiple. Hence $\SSSSs$ is closed under right-comultiple. 

\ITEM3 Assume now that $\SSSS$ is closed under right-complement, $\CCCi \SSSS  \ince \SSSS$ holds, and we have $\sss, \tt \in \SSSSs$ and $\sss \gg = \tt \ff$.  We follow the same scheme as for~\ITEM2, and keep the same notation. Assume  first that $\sss$ or~$\tt$ is invertible, say $\sss$. Put $\sss'' = \id\yy$ ($\yy$  the target of~$\tt$, $\tt' = \sss\inv \tt$, and $\hh = \ff$. Then we have $\sss  \tt' = \tt \sss''$, $\ff = \sss'' \hh$, and $\gg = \tt' \hh$. Moreover, $\sss''$  belongs to~$\SSSSs$ by definition and $\tt'$, which belongs to~$\CCCi \SSSSs$,  hence to $\CCCi \cup \CCCi \SSSS \CCCi$, belongs to~$\SSSSs$ as $\CCCi \SSSS$ is  included in~$\SSSS$. Assume now that neither $\sss$ nor $\tt$ is invertible.  Then we write $\sss = \sss' \ie$ and $\tt = \tt' \ie'$ with $\sss', \tt'$  in~$\SSSS$ and $\ie, \ie'$ in~$\CCCi$, see  Figure~\ref{F:GATransfer} again.  We have $\sss' (\ie \gg) = \tt' (\ie' \ff)$ so, as $\SSSS$ is closed  under right-complement, there exist~$\sss'', \tt''$ in~$\SSSS$ and~$\hh$  in~$\CCC$ satisfying $\sss' \tt'' = \tt' \sss''$, $\ie \gg = \tt'' \hh$, and $\ie' \ff =  \sss'' \hh$. Put $\sss''_1 = \ie'{}\inv \sss''$ and $\tt''_1 = \ie\inv \tt''$. Then  we have $\sss \tt''_1 = \tt \sss''_1$, $\ff = \sss''_1 \hh$, and $\gg = \tt''_1 \hh$.  Moreover, by construction, $\sss''_1$ and~$\tt''_1$ belong to~$\CCCi \SSSS$, hence,  by assumption, to~$\SSSSs$. Hence $\SSSSs$ is closed  under right-complement. 
\end{proof}

%%%%
\begin{ADlemm}{IV.2.39}{L:GARCClosed}
For every subfamily~$\SSSS$ of a left-cancellative category~$\CCC$ that is right-Noether\-ian and  admits   unique conditional  right-lcms, the following conditions are equivalent:

\ITEM1 The family~$\SSSSs$ is closed under right-complement (in the sense of Definition~IV.1.3);

\ITEM2 The family~$\SSSSs$ is closed under~$\under$, that is, if $\sss$ and~$\tt$ belong to~$\SSSSs$, then so does $\sss \under \tt$ when defined, that is, when $\sss$ and~$\tt$ admit a common right-multiple.
\end{ADlemm}

\begin{proof}
Assume that $\SSSSs$ is closed under right-complement. Let $\sss, \tt$ be elements of~$\SSSSs$ that admit a common right-multiple, and let $\hh$ be the right-lcm of~$\sss$ and~$\tt$. By assumption, we have $\hh = \sss (\sss \under \tt) = \tt (\tt \under \sss)$. As $\SSSSs$ is closed under right-complement, there exist~$\sss', \tt'$ in~$\SSSSs$ and~$\hh'$ satisfying $\sss \tt' = \tt \sss'$, $\tt \under \sss = \sss' \hh'$, and~$\sss \under \tt = \tt' \hh'$, whence $\hh =\nobreak (\sss \tt') \hh'$. By definition of the right-lcm, $\hh'$ must be an identity-element and, therefore, $\sss \under \tt$ and~$\tt \under \sss$ belong to~$\SSSSs$. So $\SSSSs$ is closed under the right-complement operation, and 
\ITEM1 implies~\ITEM2.

Conversely, assume that $\SSSSs$ is closed under the right-complement 
operation. Assume that $\sss, \tt$ belong to~$\SSSSs$ and $\sss \gg = \tt \ff$ holds. Then $\sss \gg$ is a common right-multiple of~$\sss$ and~$\tt$, hence it is a right-multiple of their right-lcm, which is $\sss (\sss \under \tt)$ and~$\tt (\tt \under \sss)$. So there exists~$\hh$ satisfying $\sss \gg = \sss (\sss \under \tt) \hh$. Left-cancelling~$\sss$, we deduce $\gg = (\sss \under \tt) \hh$ and, symmetrically, $\ff = (\tt \under \sss) \hh$. Then $\tt \under \sss$, $\sss \under \tt$, and~$\hh$ witness that the expected instance of closure under right-complement is satisfied. So \ITEM2 implies~\ITEM1.
\end{proof}

%%%%
\section*{Solution to selected exercises}

%%%%
\begin{exer}{38}{multiplication by invertible}{Z:GARIClosed}
Assume that $\CCC$ is a cancellative category, and $\SSSS$ is a subfamily  of~$\CCC$ that is closed under left-divisor and contains at least one element  with source~$\xx$ for each object~$\xx$. Prove $\SSSSs = \SSSS$.
\end{exer}

\begin{solu}
Assume $\gg \in \SSSS$ and $\ie \in \CCCi$. Then we have $\gg = \gg\ie \sepp  \ie\inv$, whence $\gg\ie \dive \gg$, and $\gg\ie \in \SSSS$. So $\SSSSs$ is included in~$\SSSS$.
\end{solu}

%%%%
\begin{exer}{40}{head \vs lcm}{Z:GAHeadVsLcm}
Assume that $\CCC$ is a left-cancellative category, $\SSSS$ is included in~$\CCC$, and $\gg$ belongs to~$\CCCni$. Show that $\sss$ is an $\SSSS$-head of~$\gg$ if and only if it is a right-lcm of~$\Div(\gg) \cap  \SSSS$. 
\end{exer}

\begin{solu}
If $\sss$ is an $\SSSS$-head of~$\gg$, then there exists~$\gg'$ so that $\sss \sep \gg'$ is an $\SSSS$-greedy decomposition of~$\gg$. By Lemma~IV.1.21, we deduce that, for~$\tt$ in~$\SSSS$, the relation $\tt \dive \gg$ implies $\tt \dive \sss$: this means that every  element of~$\Div(\gg) \cap \SSSS$ divides~$\sss$ and, therefore, that $\sss$  is a left-lcm of~$\Div(\gg) \cap \SSSS$. 
\end{solu}

%%%%
\begin{exer}{41}{closed under right-comultiple}{Z:GAComultiple}
Assume that $\CCC$ is a left-cancellat\-ive category, $\SSSS$ is a subfamily of~$\CCC$, and there exists $\HH : \CCCni \to \SSSS$ satisfying~(IV.1.46). Show that $\SSSS$ is closed under right-comultiple. 
\end{exer}

\begin{solu}
Assume $\ff, \gg \in \SSSS$ and $\ff \ggt = \gg \fft$.  If $\ff \ggt$ is  invertible, then everything is obvious. Otherwise, $\HH(\ff \ggt)$ is defined,  and, by~(IV.1.46)\ITEM1,  we have $\HH(\ff \ggt) \dive \ff\ggt$. On the other hand, we have $\ff \in  \SSSS$ and $\ff \dive \ff \ggt$, whence $\ff \dive \HH(\ff \ggt)$  by~(IV.1.46)\ITEM3. A symmetric argument gives $\gg \dive \HH(\ff  \ggt)$ as we have $\ff \ggt = \gg \fft$. So $\HH(\ff \ggt)$, which belongs  to~$\SSSS$ by construction, is a common right-multiple of~$\ff$ and~$\gg$ of  the  expected type.
\end{solu}

%%%%
\begin{exer}{42}{power}{Z:GAPower}
Assume that $\SSSS$ is a Garside family a left-cancellat\-ive category~$\CCC$. Show that, if $\gg_1\sep\etc\sep\gg_\pp$ is an $\SSSS^\mm$-normal decomposition of~$\gg$ and,  for every~$\ii$, the path $\sss_{\ii, 1}\sep\etc\sep\sss_{\ii, \mm}$ is an $\SSSS$-normal decomposition of~$\gg_\ii$, then the path $\sss_{1,1}\sep\etc\sep\sss_{1,\mm}\sep\sss_{2,1}\sep\etc\sep\sss_{2,\mm}\sep\etc\sep\sss_{\pp,1}\sep\etc\sep\sss_{\pp,\mm}$ is an $\SSSS$-normal decomposition of~$\gg$.
\end{exer}

\begin{solu}
By assumption, every~$\sss_{\ii, \jj}$ lies in~$\SSSSs$ and $\sss_{\ii,\jj}\sep\sss_{\ii, \jj+1}$ is $\SSSS$-greedy for all~$\ii, \jj$. So the point is to show that $\sss_{\ii, \mm}\sep\sss_{\ii+1,1}$ is $\SSSS$-greedy. Now assume $\tt \dive \sss_{\ii, \mm}\sss_{\ii+1,1}$ with $\tt$ in~$\SSSS$. Then we deduce $\sss_{\ii,1} \etc \sss_{\ii,\mm-1}\tt \dive \sss_{\ii,1} \etc \sss_{\ii,\mm-1}\sss_{\ii, \mm}\sss_{\ii+1,1}$, whence $\sss_{\ii,1} \etc \sss_{\ii,\mm-1}\tt \dive \gg_\ii\gg_{\ii+1}$. As $\gg_\ii\sep\gg_{\ii+1}$ is $\SSSS^\mm$-greedy, we deduce $\sss_{\ii,1} \etc \sss_{\ii,\mm-1}\tt \dive \gg_\ii$, that is, $\sss_{\ii,1} \etc \sss_{\ii,\mm-1}\tt \dive \sss_{\ii,1} \etc \sss_{\ii,\mm-1}\sss_{\ii, \mm}$, which implies $\tt \dive \sss_{\ii, \mm}$ since $\CCC$ is left-cancellative. As $\SSSS$ is a Garside family, Corollary~IV.1.31 implies that this is enough to conclude that $\sss_{\ii, \mm}\sep\sss_{\ii+1,1}$ is $\SSSS$-greedy.
\end{solu}

%%%%
\begin{exer}{45}{no conditional right-lcm}{Z:GAMcmNoRightLcm}
Let $\MM$ be the monoid generated by  $\tta, \ttb, \tta', \ttb'$ subject to the relations $\tta \ttb {=} \ttb \tta, \tta' \ttb' {=} \ttb' \tta', \tta \tta' {=} \ttb \ttb', \tta' \tta {=} \ttb' \ttb$.
\ITEM1 Show that the cube condition is satisfied on~$\{\tta, \ttb, \tta', \ttb'\}$, and that right- and left-reversing are complete for the above presentation. \ITEM2 Show that $\MM$ is cancellative and admits right-mcms. \ITEM3 Show that $\tta$ and $\ttb$ admit two right-mcms in~$\MM$, but they admit no right-lcm. \ITEM4 Let $\SSS = \{\tta, \ttb, \tta', \ttb', \tta \ttb, \tta' \ttb', \tta \tta', \tta' \tta\}$, see diagram on the side. Show that $\SSS$ is closed under right-mcm, and deduce that $\SSSS$ is a Garside family in~$\MM$. \ITEM5 Show that the (unique) strict $\SSS$-normal decomposition of the element~$\tta^2 \ttb' \tta'{}^2$ is~$\seqqq{\tta \ttb}{\tta' \ttb'}{\ttb'}$.
\end{exer}

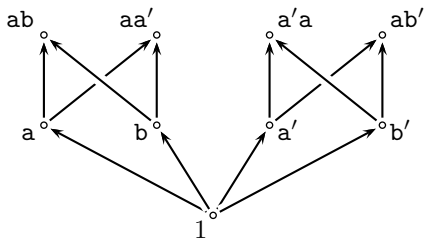
\begin{figure}[htb]
\begin{picture}(45,24)(0,0)
\psset{xunit=0.5mm,yunit=1mm}\SetGrid
\Arrow(30,01)\Arrow(30,21)\Arrow(30,41)\Arrow(30,61)\Arrow(01,02)\Arrow(01,22)\Arrow(21,02)\Arrow(21,22)\Arrow(41,42)\Arrow(41,62)\Arrow(61,42)\Arrow(61,62)
\pscircle[style=thin](45,0){0.5}\put(20,-3){$1$}
\pscircle[style=thin](0,12){0.5}\put(-3,10){$\tta$}
\pscircle[style=thin](30,12){0.5}\put(12,10){$\ttb$}
\pscircle[style=thin](60,12){0.5}\put(31,10){$\tta'$}
\pscircle[style=thin](90,12){0.5}\put(46,10){$\ttb'$}
\pscircle[style=thin](0,24){0.5}\put(-5,25){$\tta\ttb$}
\pscircle[style=thin](30,24){0.5}\put(10,25){$\tta\tta'$}
\pscircle[style=thin](60,24){0.5}\put(31,25){$\tta'\tta$}
\pscircle[style=thin](90,24){0.5}\put(46,25){$\tta\ttb'$}
\end{picture}
\end{figure}

\begin{solu}
\ITEM1 There are several relations of the form $\tta \pdots = \ttb \pdots$, but, as there is no relation $\tta \pdots = \tta'  \pdots$, there is no mixed cube to consider: the only triples to check are $(\tta, \tta, \ttb)$ and the like, and this is easy. 

\ITEM2 The presentation is homogeneous, hence (strongly) Noetherian. By Lem\-ma~II.4.62, right-reversing is complete. Hence $\MM$ is right-cancellative and admits right-mcms. By symmetry of the relations, left-reversing is complete as well, and $\MM$ is right-cancellative. 

\ITEM3 By assumption, $\tta\ttb$ and $\tta\tta'$ are common right-multiples of~$\tta$ and~$\ttb$. Owing to their length, they must be minimal. As right-reversing is complete, every common right-multiple of~$\tta$ and~$\ttb$ is a right-multiple or~$\tta\ttb$ and~$\tta\tta'$. Hence the latter are the only mcms of~$\tta$ and~$\ttb$. Finally neither is a right-multiple of the other, so $\tta$ and~$\ttb$ admit no right-lcm. 

\ITEM4 The set $\SSS$ generates~$\MM$, and it is closed under right-divisor and right-mcm: $\tta\ttb$ and $\tta\tta'$ are common right-multiples of~$\tta$ and~$\ttb$, $\tta' \ttb'$ and~$\tta' \tta$ are the two right-mcms of~$\tta'$ and~$\ttb'$, whereas none of the pairs $\{\tta, \tta'\}$, $\{\tta, \ttb'\}$, $\{\ttb, \ttb'\}$, and~$\{\ttb, \tta'\}$ admits a common right-multiple. Hence, by Corollary~IV.2.26, $\SSS$ is a Garside family in~$\MM$. 

\ITEM5 Push the letters to the left as much as possible.
\end{solu}

%%%%
\begin{exer}{47}{solid}{Z:GAGenerSolid}
Assume that $\CCC$ is a left-cancellative category and $\SSSS$ is a generating subfamily of~$\CCC$. \ITEM1 Show that $\SSSS$ is solid in~$\CCC$ if and only if $\SSSS$ includes~$\Id\CCC$ and it is closed under right-quotient. \ITEM2 Assume moreover that $\SSSS$ is closed under right-divisor. Show that $\SSSS$ includes~$\CCCi \setminus \Id\CCC$, that $\ie \in \SSSS \cap \CCCi$ implies $\ie\inv \in \SSSS$, and that and $\CCCi \SSSS = \SSSS$ holds, but that that $\SSSS$ need not be solid.
\end{exer}

\begin{solu}
\ITEM2 Assume that $\ie$ is a nontrivial invertible element, with target~$\yy$. As  $\SSSS$ generates~$\CCC$, there exists at least one element~$\ie'$ of~$\SSSS$  that right-divides~$\ie$. Then $\id\yy$ right-divides~$\ie'$, which lies  in~$\SSSS$, hence $\id\yy$ must lie in~$\SSSS$. Next, $\ie$ right-divides~$\id\yy$, which lies in~$\SSSS$, hence $\ie$ must lie in~$\SSSS$. Now assume $\ie \in \SSSS \cap \CCCi$. Then either $\ie$ is an identity-element, in which case it coincides with~$\ie\inv$ and the latter lies in~$\SSSS$, or $\ie\inv$ is invertible and is not an identity-element, in which case it belongs to~$\SSSS$ by the above argument. For $\CCCi \SSSS = \SSSS$, the proof is as for Lemma~IV.2.2. If $\yy$ is an object that is the target of no non-identity-element, then $\id\yy$: need not belong to~$\SSSS$. For instance, in the monoid~$\{1\}$, the empty set satisfies the condition. 
\end{solu}

%%%%
\begin{exer}{48}{solid}{Z:GASolid}
Let $\MM$ be the monoid $\PRESp{\tta, \tte}{\tte \tta = \tta, \tte^2 = 1}$.  \ITEM1 Show that every element of~$\MM$ has a unique expression of the  form~$\tta^\pp \tte^\qq$ with $\pp \ge 0$ and $\qq \in \{0,1\}$, and that  $\MM$ is left-cancellative. \ITEM2 Let $\SSS = \{1, \tta, \tte\}$. Show that  $\SSS$ is a solid Garside family in~$\MM$, but that $\SSS = \SSSs$ does not hold.
\end{exer}

\begin{solu}
\ITEM1 Every word in~$\{\tta, \tte\}^*$ is equivalent to a word of the  form~$\tta^\pp$ or~$\tta^\pp \tte$. Conversely, it is impossible that two  distinct words of this form are equivalent, as the number of~$\tta$ is an  invariant, and so is the parity of the number of~$\tte$ that follows the  last~$\tta$. The formulas $\tta \cdot \tta^\pp \tte^\qq = \tta^{\pp+1}  \tte^\qq$ and $\tte \cdot \tta^\pp \tte^\qq = \tta^\pp \tte^\qq$ for $\pp \ge  1$,  plus $\tte \cdot \tte^\qq = \tte^{\qq+1 \mod 2}$ show that, for  every~$\sss$ in~$\{\tta, \tte\}$, the value of~$\gg$  can be recovered  from~$\sss$ and~$\sss \gg$. \ITEM2 $\SSS$ generates~$\MM$. The explicit formula  for the multiplication shows that $\tta$ is right-divisible only by~$1$ and  itself, and so does~$\tte$. Hence $\SSS$ is closed under right-divisor, and it  is solid. For $\pp \ge 1$, define $\HH(\tta^\pp \tte^\qq) = \aa$. Then $\HH$  is a $\SSS$-head function. Hence, by Proposition~IV.2.7\ITEM1,  $\SSS$ is a Garside family in~$\MM$. On the other hand, $\tta \tte$ is an  element of~$\SSSs \setminus \SSS$.
\end{solu}

%%%%
\begin{exer}{49}{not solid}{Z:GANotSolid}
Let $\MM = \PRESp{\tta, \tte}{\tte \tta = \tta \tte, \tte^2 = 1}$, and $\SSS =  \{\tta, \tte\}$. \ITEM1 Show that $\MM$ is left-cancellative. [Hint: $\MM$ is  $\NNNN \times \ZZZZ{/}2\ZZZZ$] \ITEM2 Show that $\SSS$ is a Garside family  of~$\MM$, but $\SSS$ is not solid in~$\MM$.  [Hint: $\tte \tta$  right-divides~$\tta$, but does not belong to~$\SSS$.]
\end{exer}

\begin{solu}
\ITEM1 Every word of~$\{\tta, \tte\}^*$ is equivalent modulo the relations to  a word of the form~$\tta^\pp \tte^\qq$ with $\qq \le 1$. As $(1,\dot0)$ and  $(0, \dot1)$ satisfy in $\NNNN{\times}\ZZZZ{/}2\ZZZZ$ the defining relations  of~$\MM$, there exists a well defined homomorphism~$\FF$ of~$\MM$ to  $\NNNN{\times}\ZZZZ{/}2\ZZZZ$ that satisfies $\FF(\tta) = (1,\dot0)$ and  $\FF(\tte) = (0,\dot1)$. As $(1,\dot0)$ and $(0, \dot1)$ generate  $\NNNN{\times}\ZZZZ{/}2\ZZZZ$, the homomorphism~$\FF$ is surjective. As the  images under~$\FF$ of pairwise distinct words $\tta^\pp \tte^\qq$ are  distinct, $\FF$ is injective, hence it is a isomorphism. Hence  $\NNNN{\times}\ZZZZ{/}2\ZZZZ$ is left-cancellative, so is~$\MM$.  \ITEM2 By definition, $\SSS$ generates~$\MM$. Next, the invertible elements  of~$\MM$ are~$1$ and~$\tte$. The equality $\tte \tta = \tta \tte$ implies  $\Isom\MM \SSS \ince \SSSs$.  Finally, the  three elements of $\Pow\SSS2$, namely~$\tta^2, \tta\tte$, and~$1$, admit  $\SSS$-normal decompositions, namely for instance $\tta \sep \tta$, $\tta \sep \tte$, and $\ew$. So $\SSS$ is a Garside family in~$\MM$. Now $\tta \tte$ does  not lie in~$\SSS$,  but  we have $\tta = \tte (\tta \tte)$, so $\SSS$ is not closed  under right-divisor, hence is not solid.
\end{solu}

%%%%
\begin{exer}{50}{recognizing Garside, right-lcm  solid  case}{Z:GARecGarLcmSolid}
Assume that $\SSSS$ is a solid subfamily  in a  left-cancellative category~$\CCC$ that is right-Noetherian and admits  conditional  right-lcms.  Show that $\SSSS$ is a Garside family in~$\CCC$ if and only if  $\SSSS$ generates~$\CCC$ and  it is weakly closed under right-lcm.
\end{exer}

\begin{solu}
If $\SSSS$ is a (solid or not solid) Garside family in~$\CCC$, then, by Corollary~IV.2.29, $\SSSS$ is weakly closed under right-lcm. Conversely, assume that $\SSSS$ is solid,  generates~$\CCC$,  and is weakly closed under right-lcm. By definition, $\SSSS$ is closed under right-divisor, hence, by Lemma~IV.1.13\ITEM1, so is $\SSSSs$. Then  Corollary~IV.2.29\ITEM3 implies that $\SSSS$ is a Garside family in~$\CCC$.
\end{solu}

%%%%
\begin{exer}{52}{local right-divisibility}{Z:GALocalRightDiv}
Assume that $\CCC$ is a left-cancellative category and $\SSSS$ is a generating subfamily of~$\CCC$ that is closed under right-divisor. \ITEM1 Show that the transitive closure of~$\diveRS$ is the restriction of~$\diveR$ to~$\SSSS$. \ITEM2 Show that the transitive closure of~$\divRS$ is almost the restriction of~$\divR$ to~$\SSSS$, in the following sense: if $\sss \divR \tt$ holds, there exists~$\sss' \eqil \sss$ satisfying $\sss' \divRSs \tt$.  
\end{exer}

\begin{solu}
Let $\diveRSs$ be the transitive closure of~$\diveRS$. As $\sss \diveRS \tt$ implies $\sss \diveR \tt$ trivially and $\diveR$ is transitive, $\sss \diveR^*_{\!\SSSS} \tt$ implies $\sss \diveR \tt$. Conversely, assume $\sss \diveR \tt$, say $\tt = \gg \sss$. As $\SSSS$ generates~$\CCC$, there exist $\sss_1 \wdots \sss_\qq$ in~$\SSSS$ satisfying $\gg = \sss_1 \pdots \sss_\pp$. So we have $\tt = \sss_1 \pdots \sss_\pp \sss$. Put $\tt_\ii = \sss_\ii \pdots \sss_\pp \sss$ for $1 \le \ii \le \pp$. Each element~$\tt_\ii$ right-divides~$\tt$, an element of~$\SSSS$, so it belongs to~$\SSSS$. Now, by construction, we have $\sss = \diveRS \tt_1 \diveRS \pdots \diveRS \tt_\pp = \tt$, whence $\sss  \diveR^*_{\!\SSSS} \tt$.

\ITEM2 Let $\divRSs$ be the transitive closure of~$\divRS$. By Lemma~IV.2.15, $\sss \divRS \tt$ implies~$\sss \divR \tt$, hence $\sss \divRSs \tt$ implies~$\sss \divR \tt$ since $\divR$ is transitive. Conversely, assume $\sss \divR \tt$, say $\tt = \gg \sss$ with $\gg \notin \CCCi$. As $\SSSS$ generates~$\CCC$, there exist $\sss_1 \wdots \sss_\qq$ in~$\SSSS$ satisfying $\gg = \sss_1 \pdots \sss_\pp$. Assume that $\pp$ has been chosen minimal. Then $\sss_1 \wdots \sss_{\pp-1}$ are not invertible: if $\sss_\ii$ is invertible, then $\sss_\ii\sss_{\ii+1}$ right-divides~$\sss_{\ii+1}$, hence belongs to~$\SSSS$, and therefore grouping~$\sss_\ii$ and~$\sss_{\ii+1}$ would provide a shorter decomposition. We have $\tt = \sss_1 \pdots \sss_\pp \sss$. Put $\tt_\ii = \sss_\ii \pdots \sss_\pp \sss$ for $1 \le \ii \le \pp$. Each element~$\tt_\ii$ right-divides~$\tt$, an element of~$\SSSS$, so it belongs to~$\SSSS$. As $\SSSS$ is closed under right-divisor, we have $\CCCi \SSSS \subseteq \SSSS$. So we can assume that $\sss_1 \wdots \sss_{\pp-1}$ are non-invertible. Hence we have $\sss_1 \sss \divRSs \tt$. 
\end{solu}

%%%%
\begin{exer}{53}{local left-divisibility}{Z:GALocalLeftDiv}
Assume that $\SSSS$ is a subfamily of a left-cancellative category~$\CCC$. \ITEM1 Show that $\sss \diveS \tt$ implies $\sss \dive \tt$, and that $\sss \diveS \tt$ is equivalent to $\sss \dive \tt$ whenever $\SSSS$ is closed under right-quotient in~$\CCC$. \ITEM2 Show that, if $\SSSSi = \CCCi \cap \SSSS$ holds, then $\sss \divS \tt$ implies $\sss \div \tt$. \ITEM3 Show that, if $\SSSS$ is closed under right-divisor, then $\diveS$ is the restriction of~$\dive$ to~$\SSSS$ and, if $\SSSSi = \CCCi \cap \SSSS$ holds, $\divS$ is the restriction of~$\div$ to~$\SSSS$.
\end{exer}

\begin{solu}
\ITEM3 First $\sss \diveS \tt$ implies $\sss \dive \tt$. Conversely, assume $\sss, \tt \in \SSSS$ and $\sss \gg' = \tt$. As $\tt$ belongs to~$\SSSS$, the assumption that $\SSSS$ is closed under right-divisor implies that $\gg'$ belongs to~$\SSSS$, hence witnesses for~$\sss \diveS \tt$. Next, by Lemma~IV.2.15, $\sss \divS \tt$ implies $\sss \div \tt$. Conversely, assume $\sss \tt' = \tt$ with $\tt' \notin \CCCi$. As above, $\tt'$ must belong to~$\SSSS$, and it cannot belong to~$\SSSSi$. So $\sss \divS \tt$ holds. 
\end{solu}

%%%%
\begin{exer}{55}{locally right-Noetherian}{Z:GALocRightNoeth}
Assume that $\CCC$ is a left-cancellative category and $\SSSS$ is a subfamily of~$\CCC$. \ITEM1 Prove that $\SSSS$ is locally right-Noetherian if and only if, for  every~$\sss$ in~$\SSSSs$, every $\divS$-increasing sequence in~$\Divl\SSSS(\sss)$ is finite. \ITEM2 Assume that $\SSSS$ is locally right-Noetherian and closed under right-divisor. Show that $\SSSSs$ is locally right-Noetherian. [Hint: For $\XXX \subseteq \SSSSs$ introduce $\XXX' = \{ \sss \in \SSSS \mid \exists \ie, \ie' \in \CCCi \ (\ie \sss \ie' \in \XXX ) \}$, and construct a $\divR$-minimal element in~$\XXX$ from a $\divR$-minimal element in~$\XXX'$.]
\end{exer}

\begin{solu}
\ITEM2 Let $\XXX$ be a nonempty subfamily of~$\SSSSs$. Put 
$$\XXX' = \{ \sss \in \SSSS \mid \exists \ie, \ie' \in \CCCi \ (\ie \sss \ie' \in \XXX ) \}.$$
We have $\XXX \subseteq \XXX'$, whence $\XXX' \not= \emptyset$.  As $\XXX'$ is included in~$\SSSS$,  it contains a $\divR_\SSSS$-minimal element, say~$\sss_0$. By definition, there exists $\ie_0, \ie'_0$ invertible such that $\ie_0 \sss_0 \ie'_0$ lies in~$\XXX$. Put $\tt_0 = \ie_0 \sss_0 \ie'_0$. We claim that $\tt_0$ is $\divR_\SSSSs$-minimal in~$\XXX$. Indeed, assume $\tt_0 = \rr \tt$ with~$\rr$ in~$\SSSSs$ and~$\tt$ in~$\XXX$. We have to prove that $\rr$ lies in~$(\SSSSs)\Isom{}$.  As $\SSSSs$ is solid, it suffices to show that $\rr$ lies in~$\CCCi$.  By definition, $\rr$ belongs to~$\SSSS \CCCi \cup \CCCi$. If $\rr$ belongs to~$\CCCi$, we are done. Otherwise, write $\rr  = \rr' \ie$ with $\rr'$ in~$\SSSS$ and $\ie$ in~$\CCCi$. Then we have $\sss_0 = (\ie_0\inv \rr') (\ie \tt \ie'_0{}\inv)$.  As $\sss_0$ lies in~$\SSSS$ and $\SSSS$ is closed under right-divisor, $\ie \tt \ie'_0{}\inv$ lies in~$\SSSS$.  As $\tt$ belongs to~$\XXX$, we deduce that $\ie \tt \ie'_0{}\inv$ lies in~$\XXX'$. Then, by the choice of~$\sss_0$, the elements~$\ie_0\inv \rr'$, and therefore~$\rr'$ and~$\rr' \ie$, that is, $\rr$, must be invertible in~$\CCC$. Hence $\tt_0$ is $\divR_\SSSSs$-minimal in~$\XXX$, and $\divR_\SSSSs$ is a well-founded relation, that is, $\SSSSs$ is locally right-Noetherian.
\end{solu} 

%%%%
\chapter{Chapter~V: Bounded Garside families}

%%%%
\section*{Skipped proofs}

%%%%
\begin{ADprop}{V.1.59}{right-cancellative II}{P:DERightCancel2}
If  $\SSSS$ is a Garside family in a left-cancellat\-ive category~$\CCC$ that is right-bounded by a target-injective map~$\Delta$, and $\fD$ preserves $\SSSS$-normality and  is surjective on~$\CCCi$,  then the following conditions are equivalent:\\
\indent\ITEM1 The category~$\CCC$ is right-cancellative;\\
\indent\ITEM2 The functor~$\fD$ is injective on~$\CCC$;\\ 
\indent\ITEM3 The functor~$\fD$ is injective on~$\SSSSs$. 
\end{ADprop}

\begin{proof}
By Proposition~V.1.36, \ITEM1 and \ITEM2 are equivalent,  and  \ITEM2 obviously implies~\ITEM3. 

Now  assume that $\fD$ is injective on~$\SSSSs$.  We will prove that it is injective on~$\CCC$. First, we claim that $\fD(\sss') \eqir \fD(\sss)$ implies $\sss \eqir \sss'$ for $\sss, \sss'$ in~$\SSSSs$. Indeed, assume $\fD(\sss') = \fD(\sss) \, \ie$ with~$\ie \in \CCCi$.   If $\ie$ is trivial, the injectivity of~$\fD$ on~$\SSSSs$ implies $\sss' = \sss$. Otherwise, as $\fD$ is surjective on~$\CCCi$, there exists~$\ie'$  in~$\CCCi$ satisfying $\fD(\ie') = \ie$. As $\Delta$ is target-injective, the assumption that $\fD(\sss) \fD(\ie')$ is defined implies that $\sss \ie'$ is defined too: if $\yy$ is the target of~$\sss$ and $\xx'$ is the source of~$\ie'$, we obtain $\fD(\yy) = \fD(\xx')$, whence $\yy = \xx'$ as $\fD$ is injective on objects. We deduce $\fD(\sss') =\nobreak \fD(\sss \ie')$. As $\sss$ belongs to~$\SSSSs$, so does~$\sss \ie'$, and the assumption that $\fD$ is injective on~$\SSSSs$ then implies $\sss' = \sss \ie'$, hence $\sss' \eqir \sss$.

Now we prove using induction on~$\pp \ge 1$ the statement: $\fD(\gg) = \fD(\gg')$ implies $\gg = \gg'$ for all~$\gg, \gg'$ satisfying $\max(\LGG\SSSS\gg, \LGG\SSSS{\gg'}) \le \pp$. For $\pp = 1$, as $\LGG\SSSS\gg \le 1$ implies $\gg \in \SSSSs$, the result follows from  the assumption that $\fD$ is injective on~$\SSSSs$. Assume $\pp \ge 2$, and $\fD(\gg) = \fD(\gg')$ with $\max(\LGG\SSSS\gg, \LGG\SSSS{\gg'}) \le \pp$. Let $\seqqq{\sss_1}\etc{\sss_\pp}$ and $\seqqq{\sss'_1}\etc{\sss'_\pp}$ be $\SSSS$-normal decompositions of~$\gg$ and~$\gg'$. As $\fD$ is a functor and it preserves normality, $\seqqq{\fD(\sss_1)}\etc{\fD(\sss_\pp)}$ and $\seqqq{\fD(\sss'_1)}\etc{\fD(\sss'_\pp)}$ are $\SSSS$-normal decompositions of~$\fD(\gg)$ and~$\fD(\gg')$. By Proposition~III.1.25 (normal unique), $\fD(\gg) = \fD(\gg')$ implies $\fD(\sss_1) \eqir \fD(\sss'_1)$, whence, by the claim above, $\sss'_1 \eqir \sss_1$. Hence $\sss_1$ is an $\SSSS$-head for~$\gg$ and~$\gg'$, and we can write $\gg = \sss_1 \gg_1$, $\gg' = \sss_1 \gg'_1$, with $\max(\LGG\SSSS{\gg_1}, \LGG\SSSS{\gg'_1}) \le \pp - 1$.  Then  $\fD(\gg') = \fD(\gg)$ implies $\fD(\sss_1) \fD(\gg_1) = \fD(\sss_1) \fD(\gg'_1)$, whence $\fD(\gg_1) = \fD(\gg'_1)$. The induction hypothesis implies $\gg_1 = \gg'_1$,  whence  $\gg' = \gg$. So $\fD$ is injective on~$\CCC$, and \ITEM3 implies~\ITEM2.
\end{proof}

%%%%
\begin{ADprop}{V.2.34}{right-cancellativity III}{P:DERightCancel3}
If $\CCC$ is  a left-cancellative category and $\Delta$ is a Garside map of~$\CCC$ that preserves $\Delta$-normality, then $\CCC$ is right-cancellative if and only if $\fD$ is injective on~$\Pow{(\DIV\Delta)}2$.
\end{ADprop}

\begin{proof}
If $\CCC$ is right-cancellative, Proposition~V.1.36 implies that $\fD$ is injective on~$\CCC$, hence a fortiori on~$\Pow{(\DIV\Delta)}2$, so the condition is necessary.

Conversely, assume that $\fD$ is injective on~$\Pow{(\DIV\Delta)}2$. First, $\fD$ must be injective on~$\Id\CCC$, which is included in~$\Pow{(\DIV\Delta)}2$, and, therefore, on~$\Obj(\CCC)$, that is, $\Delta$ must be target-injective. Next, as $\Id\CCC$ is included in~$\DIV\Delta$, the assumption implies that $\fD$ is injective on~$\DIV\Delta$.

We claim that $\fD$ induces a permutation of~$\CCCi$. Indeed, as $\fD$ is a functor, it maps~$\CCCi$ into itself. Now, the assumption that $\Delta$ is a Garside map implies that $\DIV\Delta$ is bounded by~$\Delta$,  hence, by Lemma~V.2.7, $\fD$ induces a surjective map  of~$\DIV\Delta$ into itself, hence a permutation of~$\DIV\Delta$ as it is injective on~$\DIV\Delta$. Assume $\ie \in \CCCi$. Then $\ie$ and~$\ie\inv$ belong to~$\DIV\Delta$, so there  exist~$\sss, \tt$ in~$\DIV\Delta$ satisfying $\fD(\sss) = \ie$  and~$\fD(\tt) = \ie\inv$. As $\Delta$ is target-injective, $\sss \tt$ is  defined. Indeed, let~$\yy$ be the target of~$\sss$ and $\xx'$ be  the source of~$\tt$. As $\fD$ is a functor, $\fD(\yy)$ is the target  of~$\fD(\sss)$, that is, of~$\ie$, whereas $\fD(\xx')$ is the source  of~$\fD(\tt)$, that is, of~$\ie\inv$. Hence we have $\fD(\yy) = \fD(\xx')$  and, therefore, $\yy = \xx'$, that is, $\sss \tt$ is defined. The argument  for~$\tt \sss$ is symmetric. So $\sss \tt$ and $\tt \sss$ belong  to~$\Pow{(\DIV\Delta)}2$. The assumption that $\fD$ is injective  on~$\Pow{(\DIV\Delta)}2$ implies $\sss \tt = \id\xx$ and $\tt \sss = \id\yy$,  where $\xx$ (\resp $\yy$) is the source (\resp target) of~$\sss$. So $\sss$  belongs to~$\CCCi$, and $\fD$ induces a permutation of~$\CCCi$. Then, Proposition~V.1.59 implies that $\CCC$ is  right-cancellative. 

Note that the assumptions of Proposition~V.2.34 can be slightly weakened: the only assumptions used in the proof is that $\fD$ is injective  on~$\DIV\Delta$ and that, for~$\gg$ in~$\Pow{(\DIV\Delta)}2$, the relation  $\fD(\gg) \in \Id\CCC$ implies $\gg \in \Id\CCC$. We do not know whether the latter condition can be skipped.
\end{proof}

%%%%
\begin{ADlemm}{V.2.38}{L:DEGcdLcm}
\ITEM1 A left-cancella\-tive category that is left-Noether\-ian and admits left-gcds admits conditional right-lcms.

\ITEM2 In a cancellative category that admits conditional right-lcms,  any two elements of~$\CCC$ that admit a common left-multiple admit a right-gcd.
\end{ADlemm}

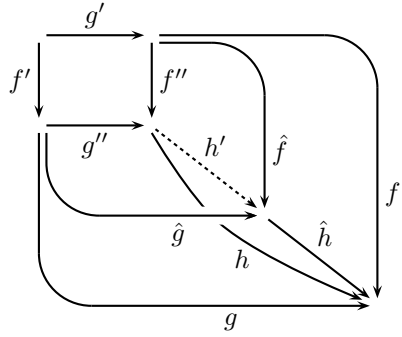
\begin{figure}[htb]\centering
\begin{picture}(45,36)(0,0)
\pcline{->}(1,36)(14,36)\taput{$\gg'$}
\pcline(16,36)(38,36)
\psarc(38,29){7}{0}{90}
\pcline{->}(45,29)(45,1)\trput{$\ff$}
\pcline{->}(0,35)(0,25)\tlput{$\ff'$}
\pcline(0,23)(0,7)
\psarc(7,7){7}{180}{270}
\pcline{->}(7,0)(44,0)\tbput{$\gg$}
\pcline{->}(1,24)(14,24)\tbput{$\gg''$}
\pscurve{->}(15,23)(24,11)(43,0.8)
\put(26,5){$\hh$}
\pcline{->}(15,35)(15,25)\trput{$\ff''$}
\pcline(16,35)(23,35)
\psarc(23,28){7}{0}{90}
\pcline{->}(30,28)(30,13)\trput{$\fft$}
\pcline(1,23)(1,19)
\psarc(8,19){7}{180}{270}
\pcline[border=3pt]{->}(8,12)(29,12)\tbput{$\ggt$}
\pcline{->}(30.5,11.5)(44,1)\put(37,8){$\hht$}
\pcline[style=exist]{->}(15.5,23.5)(29,13)\put(22,20){$\hh'$}
\end{picture}
\caption[]{\sf\small Construction of a right-gcd using a right-lcm.}
\label{F:GERightGcd}
\end{figure}

\begin{proof}
\ITEM1 We first show that every nonempty family of elements of~$\CCC$ sharing the same source admits a left-gcd. Let $\SSSS$ be a nonempty family of elements of~$\CCC$ that share the same source. An obvious induction shows that, in~$\CCC$, every finite nonempty family of elements of~$\CCC$ sharing the same source has a left-gcd. For~$\YY$ a finite nonempty subset of~$\SSSS$, choose a left-gcd~$\gg_\YY$ for~$\YY$, and let~$\SSSS'$ be the family of all elements~$\gg_\YY$. As $\CCC$ is left-Noetherian, there exists an element~$\gg_\XX$ of~$\SSSS'$ that is $\div$-minimal in~$\SSSS'$. We claim that $\gg_\XX$ is a left-gcd for~$\SSSS$. Indeed, let $\gg$ be an arbitrary element of~$\SSSS$. The point is to prove $\gg_\XX \dive \gg$. Now, by construction, we have $\gg_{\XX \cup\{\gg\}} \dive \gg$ and $\gg_{\XX \cup\{\gg\}} \dive \gg_\XX$. As $\gg_\XX$ is $\div$-minimal in~$\SSSS'$, we must have $\gg_{\XX \cup\{\gg\}} \eqir \gg_\XX$, whence $\gg_\XX \dive \gg_{\XX \cup\{\gg\}} \dive \gg$, as expected. So $\SSSS$ has a left-gcd. Then Lemma~II.2.21 guarantees that $\CCC$ admits conditional right-lcms.

\ITEM2 (See Figure~\ref{F:GERightGcd}.) Let $\ff, \gg$ be two elements of~$\CCC$ that admit a common left-multiple, say $\ff' \gg = \gg' \ff$, hence share the same target. The elements~$\ff'$ and~$\gg'$ admit a common right-multiple, namely~$\ff' \gg$, hence they admit a right-lcm, say $\ff' \gg'' = \gg' \ff''$. By definition of a right-lcm, there exists~$\hh$ satisfying $\ff = \ff'' \hh$ and~$\gg = \gg'' \hh$. By construction, $\hh$ is a common right-divisor of~$\ff$ and~$\gg$.

Let $\hht$ be a common right-divisor of~$\ff$ and~$\gg$. So there exist~$\fft, \ggt$ satisfying $\ff = \fft \hht$ and $\gg = \ggt \hht$. Then we have $\ff' \ggt \hht = \ff' \gg = \gg' \ff = \gg' \fft \hht$, whence $\ff' \ggt = \gg' \fft$ by right-cancelling~$\hht$. So $\ff' \ggt$ is a common right-multiple of~$\ff'$ and~$\gg'$, hence it is a right-multiple of their right-lcm, that is, there exists~$\hh'$ satisfying $\ff' \ggt = (\ff' \gg'') \hh'$, whence $\ggt = \gg'' \hh'$ by left-cancelling~$\ff'$. We deduce $\gg'' \hh = \gg = \ggt \hht = \gg'' \hh' \hht$, whence $\hh = \hh' \hht$ by right-cancelling~$\hht$. So $\hht$ right-divides~$\hh$, which shows that $\hh$ is a right-gcd of~$\ff$ and~$\gg$. 
\end{proof}

%%%%
\begin{ADlemm}{V.3.40}{L:DEDivDual}
If $\Delta$ is a Garside map in a cancellative category~$\CCC$, then,  for all~$\sss, \tt$ in~$\DIV\Delta$, we have $\sss \dive \tt \ \Leftrightarrow \ \dual\sss \multeR \dual\tt$ and $\sss \multeR \tt \ \Leftrightarrow \ \dualsym\sss \dive \dualsym\tt$.
\end{ADlemm}

\begin{proof}
As $\dual$ and~$\dualsym$ are bijective, everything is clear. For instance, assume $\sss \dive \tt$. Then we have $\sss \tt' = \tt$ for some~$\tt'$, which necessarily belongs to~$\DIV\Delta$. Calling~$\xx$ the common source of~$\sss$ and~$\tt$, we deduce
$$\Delta(\xx) = \sss \, \dual\sss = \tt \, \dual\tt = \sss \, \tt' \, \dual \tt,$$
whence $\dual\sss = \tt' \, \dual\tt$, and $\dual\sss \multeR \dual\tt$. The other implications are similar.
\end{proof}

%%%%
\begin{ADlemm}{V.3.46}{L:DELcmGcd}
If $\Delta$ is a Garside map in a cancellative category~$\CCC$ and  $\ff, \gg$ are elements of~$\CCC(\xx, \ud)$, then, for every~$\mm$ satisfying  $\mm \ge \max(\SUP\Delta(\ff), \SUP\Delta(\gg))$  and for every~$\hh$ in~$\Div(\DELTA\mm(\xx))$, the following conditions are equivalent:

\ITEM1 The element~$\hh$ is a right-lcm of~$\ff$ and~$\gg$;

\ITEM2 The element~$\DUAL\mm\hh$ is a right-gcd of~$\DUAL\mm\ff$ and~$\DUAL\mm\gg$.
\end{ADlemm}

\begin{proof}
Assume first that $\mm = 1$. As $\DIV\Delta$ is closed under right-comultiple, $\hh$ is a right-lcm of~$\ff$ and~$\gg$ if and only if the three relations 
$$\ff \dive \hh, \ \gg \dive \hh, \mbox{\ and\ } \forall \hh' \in \DIV\Delta\ ((\ff \dive \hh' \mbox{\ and\ } \gg \dive \hh' )  \mbox{\ implies\ } \hh \dive \hh')$$
are satisfied. By the duality formulas of Lemma~V.3.40, this is  equivalent to
$$\dual\ff \multeR \dual\hh, \ \dual\gg \multeR \dual\hh, \mbox{\ and\ } \forall \hh' \in \DIVR\Delta\ ((\dual\ff \multeR \hh' \mbox{\ and\ } \dual\gg  \multeR \hh' ) \mbox{\ implies\ } \dual\hh \multeR \hh'),$$
hence, as $\DIVR\Delta$ coincides with~$\DIV\Delta$, to $\dual\hh$ being a  right-gcd of~$\dual\ff$ and~$\dual\gg$.

Assume now that $\mm$ is arbitrary. By Proposition~V.2.32\ITEM4,  $\DELTA\mm$ is a Garside map in~$\CCC$, and $\DUAL\mm$ is the associated  duality map. Applying the above result to~$\DELTA\mm$ and~$\DUAL\mm$ gives the  expected equivalence.

(Observe that, in the situation of Lemma~V.3.46, $\DELTA\mm(\xx)$ is  a common right-multiple of~$\ff$ and~$\gg$, so every possible right-lcm  of~$\ff$  and~$\gg$ must lie in~$\Div(\DELTA\mm(\xx))$.) 
\end{proof}

%%%%
\section*{Solution to selected exercises}

%%%%
\begin{exer}{59}{preserving $\DIV\Delta$}{Z:DEDIVAuto}
Assume that $\CCC$ is a category, $\Delta$ is a map from~$\Obj(\CCC)$ to~$\CCC$ and $\phi$ is a functor from~$\CCC$ into itself that commutes with~$\Delta$. Show that $\phi$ maps~$\DIV\Delta$ and $\DIVR\Delta$ to themselves.
\end{exer}

\begin{solu}
Assume $\sss \in \DIV\Delta$, say $\sss \dive \Delta(\xx)$. By (the easy direction of) Lemma~II.2.8, this implies $\phi(\sss) \dive \phi(\Delta(\xx))$. By assumption, the latter is~$\Delta(\phi(\xx))$, so $\phi(\sss)$ belongs to~$\DIV\Delta$. The argument is similar for~$\DIVR\Delta$.
\end{solu}

%%%%
\begin{exer}{60}{preserving normality I}{Z:DEPresNormality}
Assume that $\CCC$ is a cancellative category, $\SSSS$ is a Garside family of~$\CCC$, and $\phi$ is a functor from~$\CCC$ to itself. \ITEM1 Show that, if $\phi$ induces a permutation of~$\SSSSs$, then $\phi$ preserves $\SSSS$-normality. \ITEM2 Show that $\phi$ preserves non-invertibility, that is, $\phi(\gg)$ is invertible if and only if $\gg$~is.
\end{exer}

\begin{solu}
\ITEM1 Assume that $\sss_1 \sep \sss_2$ is $\SSSS$-normal. First, by assumption,  $\phi$ maps~$\SSSSs$ to itself, hence $\phi(\sss_1)$ and~$\phi(\sss_2)$ lie  in~$\SSSSs$. Assume that $\sss$ is an element of~$\SSSS$ that satisfies $\sss  \dive \phi(\sss_1) \, \phi(\sss_2)$. By Lemma~II.2.8, we deduce  $\phi\inv(\sss) \dive \sss_1 \sss_2$, whence $\phi\inv(\sss) \dive \sss_1$ as  $\phi\inv$ maps~$\SSSSs$ into itself and $\sss_1 \sep \sss_2$, which is  $\SSSS$-normal by assumption, is also $\SSSSs$-normal by Lemma~III.1.10.  Reapplying~$\phi$, we deduce $\sss \dive \phi(\sss_1)$, and we conclude that $\phi(\sss_1) \sep \phi(\sss_2)$ is $\SSSS$-normal.

\ITEM2 First, $\phi(\gg)$ is invertible whenever $\gg$ is invertible since $\phi$ is a functor. Conversely, assume that $\phi(\gg)$ is invertible. Let $\gg_1 \sep \pdots \sep \gg_\pp$ be an $\SSSSs$-normal decomposition of~$\gg$. As $\phi$ is a functor and preserves normality, $\phi(\gg_1) \sep \etc \sep \phi(\gg_\pp)$ is an $\SSSSs$-normal decomposition of~$\phi(\gg)$. The assumption that $\phi(\gg)$ is invertible implies that each of $\phi(\gg_1) \wdots \phi(\gg_\pp)$ is invertible, and then the assumption that $\phi$ is injective on~$\SSSSs$ implies that $\gg_1, \pdots , \gg_\pp$ are invertible. Hence $\gg$ is invertible.
\end{solu}

%%%%
\begin{exer}{61}{preserving normality II}{Z:DEPresNormal}
Assume that $\CCC$ is a left-cancellative category and $\SSSS$ is a Garside family of~$\CCC$ that is right-bounded by a map~$\Delta$. \ITEM1 Show that $\fD$ preserves normality if and only if there exists an $\SSSS$-head map~$\HH$ satisfying $\HH(\fD(\gg)) \eqir \fD(\HH(\gg))$ for every~$\gg$ in~$\Pow{(\SSSSs)}2$, if and only if, for each $\SSSS$-head map~$\HH$, the above relation is satisfied.
\ITEM2 Show that a sufficient condition for $\fD$ to preserve normality is that $\fD$ preserves left-gcds on~$\SSSSs$, that is, if $\rr, \sss, \tt$ belong to~$\SSSSs$ and $\rr$ is a left-gcd of~$\sss$ and~$\tt$, then $\fD(\rr)$ is a left-gcd of~$\fD(\sss)$ and~$\fD(\tt)$.
\end{exer}

\begin{solu}
\ITEM1 Assume that $\fD$ preserves normality, $\HH$ is a $\SSSS$-head map, and $\seqq{\gg_1}{\gg_2}$ lies in~$\Seq{(\SSSSs)}2$. Put $\gg'_1 = \HH(\gg_1 \gg_2)$, and let $\gg'_2$ satisfy $\gg'_1 \gg'_2 =\nobreak \gg_1 \gg_2$. Then $\gg'_1 \sep \gg'_2$ is $\SSSS$-normal, hence, as $\fD$ preserves normality, so is $\fD(\gg'_1) \sep \fD(\gg'_2)$. Moreover, we have $\fD(\gg'_1)\fD(\gg'_2) = \fD(\gg_1\gg_2)$, hence $\fD(\gg'_1) \sep \fD(\gg'_2)$ is an $\SSSS$-normal decomposition of $\fD(\gg_1\gg_2)$. By uniqueness of the head, we must have $\HH(\fD(\gg_1\gg_2)) \eqir \fD(\gg'_1)$. Conversely, assume that $\HH$ is a $\SSSSs$-head map and $\HH(\fD(\gg)) \eqir \fD(\HH(\gg)))$ holds for every~$\gg$ in~$\Pow{(\SSSSs)}2$. Let $\gg_1 \sep \gg_2$ be an $\SSSS$-normal path. By construction, we have $\gg_1 \eqir \HH(\gg_1\gg_2)$. As $\fD$ is a functor, this implies $\fD(\gg_1) \dive \fD(\HH(\gg_1 \gg_2))$ and $\fD(\HH(\gg_1 \gg_2)) \dive \fD(\gg_1)$, hence $\fD(\gg_1) \eqir \fD(\HH(\gg_1 \gg_2))$. As we have $\fD(\HH(\gg_1\gg_2)) \eqir \HH(\fD(\gg_1\gg_2))$ by assumption, we deduce $\fD(\gg_1) \eqir \HH(\fD(\gg_1\gg_2))$. Hence $\fD(\gg_1)$ is an $\SSSS$-head for~$\fD(\gg_1 \gg_2)$. As we have $\fD(\gg_1 \gg_2) = \fD(\gg_1) \, \fD(\gg_2)$, we deduce that $\fD(\gg_1) \sep \fD(\gg_2)$ is an $\SSSS$-normal decomposition of~$\fD(\gg_1 \gg_2)$. Hence $\fD$ preserves normality. 

\ITEM2 Assume that $\fD$ preserves left-gcds on~$\SSSSs$. Let $\seqq{\gg_1}{\gg_2}$ be an $\SSSS$-normal path. Let~$\xx$ be the source of~$\gg_2$. By Proposition~V.1.53, $\dualD(\gg_1)$ and~$\gg_2$ are left-coprime, that is, $\id\xx$ is a left-gcd of~$\dualD(\gg_1)$ and~$\gg_2$. If the condition of the statement is satisfied, it follows that $\fD(\id\xx)$, which is~$\id{\fD(\xx)}$, is a left-gcd of~$\fD(\dualD(\gg_1))$ and~$\fD(\gg_2)$. By~(V.1.30), we have $\fD(\dualD(\gg_1)) = \dualD(\fD(\gg_1))$. So $\dualD(\fD(\gg_1))$ and~$\fD(\gg_2)$ are left-coprime, hence, by Proposition~V.1.53, $\fD(\gg_1) \sep \fD(\gg_2)$ is $\SSSS$-normal. 
\end{solu}

%%%%
\begin{exer}{62}{normal decomposition}{Z:DENormalCompl}
Assume that $\CCC$ is a left-cancellative category, $\Delta$ is a right-Garside map in~$\CCC$ such that $\fD$ preserves normality, and $\ff, \gg$ are elements of~$\CCC$ such that $\ff\gg$ is defined and $\ff \dive \DELTA\mm(\ud)$ holds, say $\ff\ff' = \DELTA\mm(\ud)$ with $\mm \ge 1$. Show that $\ff'$ and~$\gg$ admit a left-gcd and that, if $\hh$ is such a left-gcd, then concatenating a $\DIV\Delta$-normal decomposition of~$\ff\hh$ and a $\DIV\Delta$-normal decomposition of~$\hh\inv \gg$ yields a $\DIV\Delta$-normal decomposition of~$\ff\gg$. [Hint: First show that concatenating $\ff\hh$ and a $\DIV{\DELTA\mm}$-normal decomposition of~$\hh\inv\gg$ yields a $\DIV{\DELTA\mm}$-normal decomposition of~$\ff\gg$ and apply Exercise~42 in Chapter~IV.]
\end{exer}

\begin{solu}
By Proposition~V.1.58, $\DELTA\mm$ is a right-Garside map in~$\CCC$ and, by definition, $\ff'$ right-divides an element~$\DELTA\mm(\ud)$ so, by Lemma~V.2.36, $\ff'$ and~$\gg$ admit a left-gcd. Assume that $\hh$ is a left-gcd of~$\ff'$ and~$\gg$. Since $\DELTA\mm$ is a right-Garside map, $\gg\hh$ and~$\DELTA\mm(\ud)$, that is, $\ff\ff'$, admit a left-gcd, which (by Exercise~8\ITEM1 in Chapter~II) is of the form~$\ff\hh'$ with~$\hh'$ a left-gcd of~$\ff'$ and~$\gg$. By uniqueness, we have $\hh \eqir \hh'$, whence $\ff\hh \eqir \ff\hh'$. So $\ff\hh$ is a left-gcd of~$\ff\gg$ and~$\DELTA\mm(\ud)$, hence it is a $\DIV{\DELTA\mm}$-head of~$\ff\gg$. Hence concatenating $\ff\hh$ and a $\DIV{\DELTA\mm}$-normal decomposition of~$\hh\inv\gg$ yields a $\DIV{\DELTA\mm}$-normal decomposition of~$\ff\gg$. Then Exercise~42 gives the result.
\end{solu}

%%%%
\begin{exer}{64}{iterated duality}{Z:DEIterDual}
Assume that$\Delta$ is a Garside map in a cancellative category~$\CCC$. \ITEM1 Show that $\DUAL{\mm'}(\gg) = \DUAL\mm(\gg) \, \DELTA{\mm'-\mm}(\fD^\mm(\xx))$ holds for $\mm' \ge \mm$ and~$\gg$ in~$\Div(\DELTA\mm(\xx))$. \ITEM2 Show that $\DUALsym{\mm'}(\gg) = \DELTAsym{\mm'-\mm}(\fD^{-\mm}(\yy)) \, \DUALsym \mm(\gg)$ holds for $\mm' \ge \mm$ and~$\gg$ in~$\Divsym(\DELTAsym\mm(\yy))$.
\end{exer}

\begin{solu}
\ITEM1 By definition, we have
$$\gg \, \DUAL\mm(\gg) \, \DELTA{\mm'-\mm}(\fD^\mm(\xx))  = \DELTA\mm(\xx) \, \DELTA{\mm'-\mm}(\fD^\mm(\xx))  = \DELTA{\mm'} (\xx) = \gg \, \DUAL{\mm'}(\gg),$$
whence the result by left-cancelling~$\gg$. The computation is  symmetric for~\ITEM2.
\end{solu}

%%%%
\chapter{Chapter~VI: Germs}

%%%%
\section*{Skipped proofs}

\noindent(none)

%%%%
\section*{Solutions to selected exercises}

%%%%
\begin{exer}{65}{not embedding}{Z:GENotEmbedding}
Let $\SSSS$ consist of fourteen elements~$1, \tta \wdots \ttn$, all with the same source and target, and $\OP$ be defined by $1 \OP \xx = \xx \OP 1 = \xx$ for each~$\xx$, plus $\tta \OP \ttb = \ttf$, $\ttf \OP \ttc = \ttg$, $\ttd \OP \tte = \tth$, $\ttg \OP \tth = \tti$, $\ttc \OP \ttd = \ttj$, $\ttb \OP \ttj = \ttk$, $\ttk \OP \tte = \ttm$, and $\tta \OP \ttm = \ttn$. \ITEM1 Shows that $\SSSSg$ is a germ. \ITEM2 Show that, in~$\SSSS$, we have $((\tta \OP \ttb) \OP \ttc) \OP (\ttd \OP \tte) = \tti \not= \ttn = \tta \OP ((\ttb \OP (\ttc \OP \ttd)) \OP \tte)$, whereas, in~$\Mon\SSSSg$, we have $\iota\tti = \iota\ttn$. Conclude.
\end{exer}

\begin{solu}
\ITEM1 The only nontrivial triples eligible for~(VI.1.6) are $(\tta, \ttb, \ttj)$, $(\ttf, \ttc, \ttd)$, and $(\ttc, \ttd, \tte)$, for which (VI.1.6) is actually true. \ITEM2 The equalities in~$\SSSSg$ follow from the multiplication table of~$\SSSSg$. On the other hand, in $\Mon\SSSSg$, we find $\iota\tti = ((\iota\tta \, \iota\ttb) \, \iota\ttc) \, (\iota\ttd \, \iota\tte) = \iota\tta \, ((\iota\ttb \, (\iota\ttc \, \iota\ttd)) \, \iota\tte) = \iota\ttn$ applying associativity. Thus $\iota$ is not injective, and $\SSSSs$ does not embed in~$\Mon\SSSSg$.
\end{solu}

%%%%
\begin{exer}{66}{multiplying by invertible elements}{Z:GEMultInv}
\ITEM1 Show that, if $\SSSSg$ is a left-associat\-ive germ, then $\SSSS$ is closed under left-multiplication by invertible elements in~$\Cat\SSSSg$.
\ITEM2 Show that, if $\SSSSg$ is an associative germ, $\sss \OP \tt$ is defined, and $\tt' \eqirS \tt$ holds, then $\sss \OP \tt'$ is defined as well. 
\end{exer}

\begin{solu}
\ITEM1 Assume that $\ie$ admits a left-inverse~$\ie'$---as nothing a priori forces the category~$\Cat\SSSSg$ to be left-cancellative, we have to distinguish between left- and right-inverse---and that $\ie \gg$ is defined for some~$\gg$ lying in~$\SSSS$. Then we have $\gg = \ie' \ie \gg$, so $\ie \gg$ is a right-divisor of an element of~$\SSSS$, hence is an element of~$\SSSS$ as the latter is closed under right-divisor. Applying this to the case when $\gg$ is an identity-element shows that the family of all left-invertible elements is included in~$\SSSS$.
\end{solu}

%%%%
\begin{exer}{67}{atoms}{Z:GEAtoms}
\ITEM1 Show that, if $\SSSSg$  is a left-associative germ, the atoms of~$\Cat\SSSSg$ are the elements of the form~$\tt\ie$ with $\tt$ an atom of~$\SSSSg$ and $\ie$ an invertible element of~$\SSSSg$.

\rightskip30mm
\ITEM2 \VR(3,0)Let~$\SSSSg$ be the germ whose table is shown on the right. Show that the monoid $\Mon\SSSSg$ admits the presentation $\PRESp{\tta, \tte}{\tte \tta = \tta, \tte^2 = 1}$ (see Exercise~48) and that $\tta$ is the only atom of~$\SSSSg$, whereas the atoms of~$\Mon\SSSSg$ are $\tta$ and $\tta\tte$. \ITEM3 Show that $\SSSSg$ is a Garside germ.\par\hfill\begin{picture}(0,0)(22,-15)
\begin{tabular}{c|c@{\HS2}c@{\HS2}c@{\HS2}c@{\HS2}c@{\HS2}c}
$\OP$&$1$&$\tta$&$\tte$\\
\hline
$1$&$1$&$\tta$&$\tte$\\
$\tta$&$\tta$\\
$\tte$&$\tte$&$\tta$&$1$
\end{tabular}
\end{picture}
\end{exer}

\vspace{-3mm}\begin{solu}
\ITEM1 Assume $\sss = \tt \ie$ with $\tt$ an atom of~$\SSSSg$ and $\ie$ an invertible element in~$\SSSS$. Let $\seqqq{\sss_1}\etc{\sss_\pp}$ be a decomposition of~$\sss$ in~$\Cat\SSSSg$. As $\SSSS$ generates~$\Cat\SSSSg$, each element~$\sss_\ii$ can be expressed as a product of elements of~$\SSSS$, leading to a new decomposition $\seqqq{\tt_1}\etc{\tt_\qq}$ of~$\sss$ in terms of elements of~$\SSSS$, whence $\tt = \tt_1 \pdots \tt_\qq \ie\inv$ in~$\Cat\SSSSg$ and, as $\SSSSg$ is left-associative, $\tt = \Pi(\seqqqq{\tt_1}\etc{\tt_\qq}{\ie\inv})$ in~$\SSSSg$, where $\Pi$ is the partial map  of~(VI.1.15).  As $\tt$ is an atom of~$\SSSSg$, at most one of the entries~$\tt_\jj$ is non-invertible in~$\SSSSg$. Hence at most one of the entries~$\sss_\ii$ is non-invertible in~$\Cat\SSSSg$,  and  $\sss$ is an atom in~$\Cat\SSSSg$.

Conversely, assume that $\sss$ is an atom in~$\Cat\SSSSg$. As $\SSSS$ generates~$\Cat\SSSSg$, there exists a decomposition $\seqqq{\sss_1}\etc{\sss_\pp}$ of~$\sss$ into elements of~$\SSSS$. As $\sss$ is an atom in~$\Cat\SSSSg$, at most one entry is not invertible. If every entry is invertible, then $\sss$ is invertible, contradicting the assumption. So exactly one entry, say~$\sss_\ii$, is not invertible. Let $\tt = \sss_1 \pdots \sss_{\ii-1} \sss_\ii$ and $\ie = \sss_{\ii+1} \pdots \sss_\pp$. As $\sss_1 \wdots \sss_{\ii-1}$ are invertible, $\tt$ is a right-divisor of~$\sss$, hence it belongs to~$\SSSS$. Moreover, $\tt$ is an atom of~$\SSSS$, since a decomposition of~$\tt$ with more than one non-invertible entry in~$\SSSS$ would provide a similar decomposition in~$\Cat\SSSSg$, contradicting the assumption. On the other hand, $\ie$ is invertible. So $\sss$, which is~$\tt\ie$, has the expected form.
\end{solu}

%%%%
\begin{exer}{68}{families $\III_{\SSSSg}$ and $\JJJ_{\SSSSg}$}{Z:GEIJFam}
Assume that $\SSSSg$ is a left-associative germ. \ITEM1 Show that a path $\seqq{\sss_1}{\sss_2}$ of~$\Seq\SSSS2$ is $\SSSS$-normal if and only if all elements of~$\JJJ_{\SSSSg}(\sss_1, \sss_2)$ are invertible. \ITEM2 Assuming in addition that $\SSSSg$ is left-cancellative, show that, for $\seqq{\sss_1}{\sss_2}$ in~$\Seq\SSSS2$, the family~$\III_{\SSSSg}(\sss_1, \sss_2)$ admits common right-multiples if and only if $\JJJ_{\SSSSg}(\sss_1, \sss_2)$ does.
\end{exer}

\begin{solu}
\ITEM2 Assume that $\III_{\SSSSg}(\gg_1, \gg_2)$ admits common right-multiples, and let $\hh, \hh'$ belong to~$\JJJ_{\SSSSg}(\gg_1, \gg_2)$. Then $\gg_1 \OP \hh$ and $\gg_2 \OP \hh'$ are defined and belong to~$\III_{\SSSSg}(\gg_1, \gg_2)$, hence, by assumption, $\gg_1 \OP \hh$ and $\gg_1 \OP \hh'$ admit a common right-multiple, say $\gg''$. Then we have $\gg_1 \diveS \gg_1 \OP \hh \diveS \gg$, whence $\gg_1 \diveS \gg''$. So there exists~$\hh''$ in~$\SSSS$ satisfying $\gg'' = \gg_1 \OP \hh''$. By Lemma~VI.1.19, $\gg_1 \OP \hh \diveS \gg_1 \OP \hh''$ implies $\hh \diveS \hh''$, and, similarly, we find $\hh' \diveS \hh''$. So $\hh''$ is a common right-multiple of~$\hh$ and~$\hh'$ in~$\SSSS$. Moreover the assumption that $\gg_1 \OP \hh''$ belongs to~$\III_{\SSSSg}(\gg_1, \gg_2)$ implies that $\hh''$ belongs to~$\JJJ_{\SSSSg}(\gg_1, \gg_2)$. So $\JJJ_{\SSSSg}(\gg_1, \gg_2)$ admits common right-multiples. The converse implication is similar, actually simpler as no cancellation is needed.
\end{solu}

%%%%
\begin{exer}{69}{positive generators}{Z:GEPosGen}
Assume that $\Sigma$ is a family of positive generators in a group~$\GGG$ and $\Sigma$ is closed under inverse, that is, $\gg \in \Sigma$ implies $\gg\inv \in \Sigma$. \ITEM1 Show that $\LT\gg\Sigma = \LT{\gg\inv}\Sigma$ holds for every~$\gg$ in~$\GGG$. \ITEM2 Show that $\ff\inv \Pref\Sigma \gg\inv$ is equivalent to $\ff \Suff\Sigma \gg$.
\end{exer}

\begin{solu}
\ITEM1 An $\SSS$-word~$\ww$ is a minimal length expression for an element~$\gg$ if and only if $\ww\inv$, which is also an $\SSS$-word, is a minimal length expression of~$\gg\inv$. \ITEM2 By definition, $\ff\inv \Pref\Sigma \gg\inv$ is equivalent to $\LT{\ff\inv}\Sigma + \LT{(\ff\inv)\inv \gg\inv}\Sigma = \LT{\gg\inv}\Sigma$, hence, by~\ITEM1, to $\LT\ff\Sigma + \LT{\ff\gg\inv}\Sigma = \LT\gg\Sigma$ and to $\LT\ff\Sigma + \LT{\gg\ff\inv}\Sigma = \LT\gg\Sigma$. The latter is $\ff \Suff\Sigma \gg$.
\end{solu}

%%%%
\begin{exer}{70}{minimal upper bound}{Z:GEMinUpper}
For $\le$ a partial ordering on a family~$\SSSS'$ and $\ff, \gg, \hh$ in~$\SSSS'$, say that $\hh$ is a \emph{minimal upper bound}, or \emph{mub}, for~$\ff$ and~$\gg$, if $\ff \le \hh$ and $\gg \le \hh$ holds, but there exists no~$\hh'$ with $\hh' < \hh$ satisfying $\ff \le \hh'$ and $\gg \le \hh'$. Assume that $\GGG$ is a groupoid, $\Sigma$ positively generates~$\GGG$, and $\HHH$ is a subfamily of~$\GGG$ that is closed under $\Sigma$-suffix. Show that $\Der\HHH\Sigma$ is a Garside germ if and only if, for all $\ff, \gg, \gg', \gg''$ in~$\HHH$ such that $\ff \OPDer \gg$ and $\ff \OPDer \gg'$ are defined and $\gg''$ is a $\Pref\Sigma$-mub of~$\gg$ and~$\gg'$, the product $\ff \OPDer \gg''$ is defined.
\end{exer}

\begin{solu}
By Lemmas~VI.2.60 and~VI.2.62, the germ~$\Der\HHH\Sigma$ is left-associative, cancellative, and Noetherian. Hence, by Proposition~VI.2.44, $\Der\HHH\Sigma$ is a Garside germ if and only if it satisfies (VI.2.43). Now, by Lemma~VI.2.62, for~$\gg, \hh$ in~$\HHH$, the relation $\gg \diveLloc{\Der\HHH\Sigma} \hh$ is equivalent to $\gg \Pref\Sigma \hh$ and, therefore, $\gg''$ is an mcm of~$\gg$ and~$\gg'$ in~$\Der\HHH\Sigma$ if and only it it is a $\Pref\Sigma$-mub of~$\gg$ and~$\gg'$. So the condition is a direct reformulation of(VI.2.43).
\end{solu}

%%%%
\chapter{Chapter~VII: Subcategories}

%%%%
\section*{Skipped proofs}

%%%%
\begin{ADlemm}{VII.1.3}{L:SUInv}
If $\CCCu$ is a subcategory of a left-cancellative category~$\CCC$, we have
\begin{equation}
\ADlabel{VII.1.4}
\CCCui \subseteq \CCCi \cap \CCCu,
\end{equation}
with equality if and only if $\CCCu$ is closed under inverse in~$\CCC$. For every subfamily~$\SSSS$ of~$\CCC$, putting $\SSSSu = \SSSS \cap \CCCu$ and $\SSSSus = \SSSSu\CCCui \cup \CCCui$, we have 
\begin{equation}
\ADlabel{VII.1.5}
\SSSSus \subseteq \SSSSs \cap \CCCu.
\end{equation}
If $\CCC$ has no nontrivial invertible element, then so does~$\CCCu$, and (VII.1.4)--(VII.1.5) are equalities.
\end{ADlemm}

\begin{proof}
If $\ie \ie' = \id\xx$ holds in~$\CCCu$, it holds in~$\CCC$ as well, so (VII.1.4) is clear. For~(VII.1.4) to be an equality means that every invertible element lying in~$\CCCu$ belongs to~$\CCCui$, that is, has an inverse that lies in~$\CCCu$: this means that $\CCCu$ is closed under inverse in~$\CCC$.

Next, assume that $\SSSS$ is included in~$\CCC$, and put $\SSSSu = \SSSS \cap \CCCu$ and $\SSSSus = \SSSSu\CCCui \cap \CCCui$. Then $\CCCui$ is included in~$\CCCi$ and in~$\CCCu$, so we deduce $\SSSSus \subseteq (\SSSS\CCCi \cup \CCCi) \cap \CCCu = \SSSSs \cap \CCCu$.

On the other hand, assume $\CCCi = \Id\CCC$. As an identity-element is its own inverse, it is invertible in every subcategory that contains it, and we obtain $\CCCui = \Id\CCC \cap \CCCu = \CCCi \cap \CCCu$ and, for every~$\SSSS \subseteq \CCC$, as $\SSSSs$ is then $\SSSS \cup \Id\CCC$, we obtain 
$$\SSSSs \cap \CCCu 
= (\SSSS \cap \CCCu) \cup (\Id\CCC \cap \CCCu) 
= \SSSSu \cup \CCCui
= \SSSSu \CCCui \cup \CCCui = \SSSSus.\eqno{\square}$$
\def\qed{\relax}\end{proof}

%%%%
\begin{ADlemm}{VII.1.16}{L:SUClosedFat}
Every subcategory that is closed under left- or under right-divisor in a left-cancellative category~$\CCC$ is closed under~$\eqir$. 
\end{ADlemm}

\begin{proof}
 Assume that $\CCCu$ is a subcategory of~$\CCC$, and we have $\gg \in \CCCu$ and $\gg' \eqir \gg$, say $\gg' = \gg \ie$ with~$\ie$ invertible. If $\CCCu$ is closed under left-divisor, we can write $\gg = \gg' \ie\inv$, and $\gg'$ is a left-divisor of~$\gg$, hence it belongs to~$\CCCu$.

On the other hand, assume that $\CCCu$ is closed under right-divisor. 
Let~$\yy$ be the target of~$\gg$, and $\yy'$ be that of~$\gg'$. The  assumption  that $\gg$ lies in~$\CCCu$ implies that $\yy$ lies in~$\Obj(\CCCu)$, hence $\id\yy$ lies in~$\CCCu$. Now  $\ie\inv$ is a right-divisor of~$\gg$, hence it lies in~$\CCCu$, its  source~$\yy'$ lies  in~$\Obj(\CCCu)$ and $\id{\yy'}$ lies in~$\CCCu$. Now, we have $\id{\yy'} = \ie\inv \ie$, so $\ie$, a right-divisor of~$\id\zz$, must lie in~$\CCCu$.  Hence $\gg'$,  that is, $\gg \ie$, lies in~$\CCCu$. 
\end{proof}

%%%%
\begin{ADlemm}{VII.1.17}{L:SUInvIm}
If $\CCC, \CCC'$ are left-cancellative categories, $\fct : \CCC \to \CCC'$ is a functor, and $\CCCu'$ is a subcategory of~$\CCC'$ that is closed under left-divisor (\resp right-divisor), then so is~$\fct\inv(\CCCu')$. 
\end{ADlemm}

\begin{proof}
Assume that $\CCCu'$ is closed under left-divisor. Put $\CCCu = \fct\inv(\CCCu')$, and assume $\ff \dive \gg \in \CCCu$. By definition, there exists~$\gg'$ satisfying $\ff \gg' = \gg$, which implies $\fct(\ff) \fct(\gg') = \fct(\gg)$, whence $\fct(\ff) \dive \fct(\gg)$. By assumption, $\fct(\gg)$ lies in~$\CCCu'$, hence so does~$\fct(\ff)$ as $\CCCu'$ is closed under left-divisor. Hence $\fct(\ff)$ belongs to~$\CCCu'$, and $\ff$ lies in~$\CCCu$. So $\CCCu$ is closed under left-divisor. The argument when $\CCCu'$ is closed under right-divisor is symmetric.
\end{proof}

%%%%
\begin{ADprop}{VII.2.16}{recognizing compatible}{P:SURecCompat3}
If $\SSSS$ is a Garside family in a left-cancell\-ative category~$\CCC$ and $\CCCu$ is a subcategory of~$\CCC$ that is closed under right-quotient in~$\CCC$, then $\CCCu$ is compatible with~$\SSSS$ if and only if  
\begin{gather}
\ADlabel{VII.2.17}
\parbox{100mm}{The family $\SSSSus$ generates~$\CCCu$, where we put $\SSSSu = \SSSS \cap \CCCu$ and $\SSSSus = \SSSSu\CCCui \cup \CCCui$,}\\
\ADlabel{VII.2.18}
\parbox{100mm}{Every element of~$\Pow{(\SSSSus)}2$ admits an $\SSSS$-normal decomposition with entries \rlap{in~$\SSSSus$.}}
\end{gather}
\end{ADprop}%

\begin{proof}
If $\CCCu$ is compatible with~$\SSSS$, then, by Proposition~VII.2.14, (VII.2.15) holds. This implies in particular that every element of~$\CCCu$ admits a decomposition with entries in~$\SSSSus$, hence that $\SSSSus$ generates~$\CCCu$, which is~(VII.2.17). On the other hand, applying (VII.2.15) to an element of~$\Pow{(\SSSSus)}2$ gives~(VII.2.18).

Conversely, assume that (VII.2.17) and~(VII.2.18) are satisfied. As $\CCCu$ is closed under right-quotient, it is closed under inverse, $\SSSSus \cap \CCCi = \CCCui$ holds, hence so does $\SSSSus(\SSSSus \cap\nobreak \CCCi) \subseteq \SSSSus$. Then Lemma~VII.2.19 is valid for~$\SSSSus$, so every element of the subcategory of~$\CCC$ generated by~$\SSSSus$, hence of~$\CCCu$, admits an $\SSSS$-normal decomposition whose entries lie in~$\SSSSus$. In this case, (VII.2.15) is satisfied, and, by Proposition~VII.2.14, $\CCCu$ is compatible with~$\SSSS$.
\end{proof}

%%%%
\begin{ADlemm}{VII.2.19}{L:SUConstNormal}
Assume that $\SSSS$ is a Garside family in a left-cancellative category~$\CCC$ and $\SSSS'$ is a subfamily of~$\SSSSs$ such that $\SSSS'(\SSSS' \cap \CCCi) \subseteq \SSSS'$ holds and every element of~$\Pow{\SSSS'}2$ admits a $\SSSS$-normal decomposition with entries in~$\SSSS'$. Then every element in the subcategory of~$\CCC$ generated by~$\SSSS'$ admits a $\SSSS$-normal decomposition with entries in~$\SSSS'$. 
\end{ADlemm}

\begin{proof}
First we claim that every element~$\gg$ of~$\Pow{\SSSS'}2$ admits a length two $\SSSS$-normal decomposition with entries in~$\SSSS'$: indeed, as $\SSSS'$ is included in~$\SSSSs$, the $\SSSS$-length of~$\gg$ is at most two, so the possible entries at position~$3$ and beyond in an $\SSSS$-normal decomposition of~$\gg$ must be invertible, and the assumption $\SSSS'(\SSSS' \cap \CCCi) \subseteq \SSSS'$ implies that the latter can be incorporated in the second entry, yielding an $\SSSS$-normal decomposition of length two.

Then the argument is exactly the same as for Proposition~III.1.49 (left multiplication), except that, here, we use two different reference families, namely~$\SSSS'$ for the entries of the decomposition and $\SSSS$ for the greedy property. The point is to prove that, for every~$\pp$, every element of~$\Pow{\SSSS'}\pp$ admits an $\SSSS$-normal decomposition of length~$\pp$ with entries in~$\SSSS'$. We use induction on~$\pp \ge 1$. For $\pp = 1$, the result is obvious as $\SSSS'$ is included in~$\SSSSs$, and, for $\pp = 2$, the result is the assumption. Let $\gg$ belong to~$\Pow{\SSSS'}\pp$ with $\pp \ge 3$. Write $\gg = \sss \gg'$ with $\sss \in \SSSS'$ and $\gg' \in \Pow{\SSSS'}{\pp-1}$. By induction hypothesis, $\gg'$ admits an $\SSSS$-normal decomposition $\seqqq{\sss'_1}\etc{\sss'_{\pp-1}}$ with entries in~$\SSSS'$. Starting with $\sss_0 = \sss$, and applying the assumption $\pp-1$ times, we find an $\SSSS$-normal decomposition $\seqq{\gg_\ii}{\sss_\ii}$ of $\sss_{\ii-1} \gg'_\ii$ with entries in~$\SSSS'$. By the first domino rule (Proposition~III.1.45), $\seqqqq{\sss_1}\etc{\sss_{\pp-1}}{\sss_\pp}$ is an $\SSSS$-normal decomposition of~$\gg$ with entries in~$\SSSS'$.
\end{proof}

%%%%
\begin{ADlemm}{VII.4.7}{L:SUCorrClosed}
If $\SSSS$ is any subfamily of a left-cancellative category~$\CCC$, the identity-functor on~$\CCC$ is correct for inverses (\resp right-comultip\-les, \resp right-complements, \resp right-diamonds) on~$\SSSS$ if and only if $\SSSS$ is closed under inverse (\resp right-comultiple, \resp right-complement, \resp right-diamond) in~$\CCC$.
\end{ADlemm}

\begin{proof}
The identity-functor on~$\CCC$ is correct for inverses on~$\SSSSu$ if and only if, for every~$\sss$ in~$\SSSSu$ that is invertible, $\sss\inv$ belongs to~$\SSSSu$. By definition, this means that $\SSSSu$ is closed under inverse in~$\CCC$.

Next, the identity-functor on~$\CCC$ is correct for right-comultiples on~$\SSSSu$ if and only if, when $\sss, \tt$ lie in~$\SSSSu$ and $\sss \gg = \tt \ff$ holds in~$\CCC$, there exist~$\sss', \tt'$, and~$\hh$ satisfying $\sss \tt' = \tt \sss'$, $\ff = \sss' \hh$, and $\gg = \tt' \hh$, plus $\sss \tt' \in \SSSSu$. By definition, this means that $\SSSSu$ is closed under right-comultiple in~$\CCC$. The result is similar with right-complements and right-diamonds.
\end{proof}

%%%%
\section*{Solutions to selected exercises}

%%%%
\begin{exer}{72}{$\eqir$-closed  subcategory}{Z:SUFat}
\ITEM1 Show that a subcategory~$\CCCu$ of a left-cancellative category~$\CCC$ is  $\eqir$-closed  if and only if, for each~$\xx$ in~$\Obj(\CCCu)$, the families $\CCCi(\xx, \ud)$ and $\CCCi(\ud, \xx)$ are included in~$\CCCu$. \ITEM2 Deduce that $\CCCu$ is  $\eqir$-closed  if and only if $\CCCui$ is a union of connected components of~$\CCCi$.
\end{exer}

\begin{solu}
Assume that $\CCCu$ is $\eqir$-closed and $\xx$ lies in~$\Obj(\CCCu)$. Then $\id\xx$ belongs to~$\CCCu$. If $\ie$ belongs to~$\CCCi(\xx, \ud)$, we have $\ie \eqir \id\xx$, whence $\ie \in \CCCu$. On the other hand, if $\ie$ belongs to~$\CCCi(\yy, \xx)$, then $\ie\inv$ belongs to~$\CCCi(\xx, \yy)$. By the above result, $\ie\inv$ belongs to~$\CCCu$. It follows that $\yy$, the target of~$\ie\inv$, lies in~$\Obj(\CCCu)$, and, therefore, $\ie$ belongs to~$\CCCu$. Conversely, assume that $\CCCi(\xx, \ud)$ is included in~$\CCCu$ for every~$\xx$ in~$\Obj(\CCCu)$, and let $\gg$ belong to~$\CCCu$ and $\gg' \eqir \gg$ hold, say $\gg' = \gg \ie$ with~$\ie$ in~$\CCCi$. Let $\xx$ be the target of~$\gg$. Then $\xx$ belongs to~$\Obj(\CCCu)$, hence $\ie$ belongs to~$\CCCu$, so $\gg \in \CCCu$ implies $\gg' \gg\ie \in \CCCu$.
\end{solu}

%%%%
\begin{exer}{73}{greedy paths}{Z:SUSub1}
Assume that $\CCC$ is a cancellative category, $\SSSS$ is included in~$\CCC$, and $\CCCu$ is a subcategory of~$\CCC$ that is closed under left-quotient. Put $\SSSSu = \SSSS \cap \CCCu$. Show that every $\CCCu$-path that is $\SSSS$-greedy in~$\CCC$ is $\SSSSu$-greedy in~$\CCCu$.
\end{exer}

\begin{solu}
With the notation of the proof of Lemma~VII.2.1, we have also $\ff' = \ff'' \gg_2$ with $\ff'$ and~$\gg_2$ in~$\CCCu$. This implies $\ff'' \in \SSSSu$ as $\CCC$ is right-cancellative and $\CCCu$ is closed under left-quotient in~$\CCC$.
\end{solu}

%%%%
\begin{exer}{74}{compatibility with~$\CCC$}{Z:SUCompatTrivial}
Assume that $\CCC$ is a left-cancellative category. Show that every subcategory of~$\CCC$ that is closed under inverse is compatible with~$\CCC$ viewed as a Garside family in itself.
\end{exer}

\begin{solu}
Let $\CCCu$ be a subcategory of~$\CCC$ that is closed under inverse. Then $\CCC \cap \CCCu = \CCCu$ is a Garside family in~$\CCCu$. A $\CCCu$-path $\seqqq{\gg_1}\etc{\gg_\qq}$ is $\CCC$-normal if and only if $\gg_2 \wdots \gg_\pp$ belong to~$\CCCi$, hence if and only if $\gg_2 \wdots \gg_\pp$ belong to~$\CCCui$, hence if and only if $\seqqq{\gg_1}\etc{\gg_\qq}$ is $\CCCu$-normal.
\end{solu}

%%%%
\begin{exer}{75}{not compatible}{Z:SUNotCompat}
Let $\MM$ be the free Abelian monoid generated by~$\tta$ and~$\ttb$, and let $\NN$ be the submonoid generated by~$\tta$ and~$\tta \ttb$. \ITEM1 Show that $\NN$ is not closed under right-quotient in~$\MM$. \ITEM2 Let $\SSS = \{1, \tta, \ttb, \tta\ttb\}$. Show that $\NN$ is not compatible with~$\SSS$. \ITEM3 Let $\SSS' = \{\tta^\pp \ttb^\ie \mid \pp \ge 0, \ii \in \{0,1\}\}$. Show that $\NN$ not compatible with~$\SSS'$. 
\end{exer}

\begin{solu}
\ITEM1 The elements $\tta\ttb$ and $\tta$ lie in~$\NN$, but the right-quotient~$\ttb$ does not. 

\ITEM2 The family $\SSS \cap \NN$ is not a Garside family in~$\NN$, because $\tta$ and $\tta\ttb$, which belong to~$\SSS \cap \NN$, have no common right-multiple belonging to~$\SSS \cap \NN$. 

\ITEM3 The family $\SSS' \cap \NN$ is a Garside family in~$\MM$, but the $\SSS'$-normal decomposition of~$\tta^2\ttb^2$ is $\tta^2\ttb \sep \ttb$, whereas the $(\SSS' \cap \NN)$-normal decomposition of~$\tta^2 \ttb^2$ in~$\NN$ is $\tta\ttb \sep \tta\ttb$.
\end{solu}

%%%%
\begin{exer}{76}{not closed under right-quotient}{Z:SUNotRQClosed}
\ITEM1 Show that every submonoid~$\mm\NNNN$ of the additive monoid~$\NNNN$ is closed under right-quotient, but that $2\NNNN + 3\NNNN$ of~$\NNNN$ is not. 
\ITEM2 Let $\MM$ be the monoid $\NNNN \rtimes (\ZZZZ/2\ZZZZ)^2$, where the generator~$\tta$ of~$\NNNN$ acts on the generators~$\tte, \ttf$ of $(\ZZZZ/2\ZZZZ)^2$ by $\tta \tte = \ttf \tta$ and $\tta \ttf = \tte \tta$, and let $\NN$ be the submonoid of~$\MM$ generated by~$\tta$ and~$\tte$. Show that $\MM$ is left-cancellative, and its elements admit a unique expression of the form $\tta^\pp \tte ^\ii \ttf^\jj$ with $\pp \ge 0$ and $\ii, \jj \in \{0,1\}$, and that $\NN$ is $\MM \setminus \{\ttf, \tte\ttf\}$. \ITEM3 Show that $\NN$ is not closed under right-quotient in~$\MM$. \ITEM4 Let $\SSS = \{\tta\}$. Show that $\SSS$ is a Garside family in~$\MM$ and determine~$\SSSs$. Show that $\NN$ is compatible with~$\SSS$. [Hint: Show that $\SSS \cap \NN$, which is~$\SSS$, is not a Garside family in~$\NN$.] \ITEM5 Shows that $\SSSs \cap \NN$ is a Garside family in~$\NN$ and $\NN$ is compatible with~$\SSSs$.
\end{exer}

\begin{solu}
\ITEM1 If $\ff$ and $\ff + \gg'$ are multiples of~$\mm$, then $\gg'$ is a multiple of~$\mm$ as well. On the other hand, $2\NNNN + 3\NNNN$ contains $2$ and $3$, but does not contain~$1$ although $3 = 2 + 1$ holds. 

\ITEM3 The set $\Isom\MM$ consists of~$1, \tte, \ttf, \tte\ttf$, whereas $\Isom{\NN}$ consists of~$1$ and~$\tte$, so $\NN$ is closed under inverse in~$\MM$. On the other hand, as $\tta$ and~$\tta\ttf$ belong to~$\NN$ but $\ttf$ does not, $\NN$ is not closed under right-quotient in~$\MM$.

\ITEM4 We have $\SSSs = \{\tta^\pp \tte ^\ii \ttf^\jj \mid \pp, \ii, \jj \le 1\}$ (eight elements). Then $\tta\ttf$, which is $\tte\tta$, belongs to~$\Isom{\NN}\SSS$ but not to~$\SSS\Isom{\NN} \cup \Isom{\NN}$. So $\SSS \cap \NN$, which is~$\SSS$, is not a Garside family in~$\NN$, and the submonoid~$\NN$ is not compatible with~$\SSS$.

\ITEM5 Consider $\SSSs$, which is also a Garside family in~$\MM$. Then $\SSSs \cap \NN$ is $\SSSs \setminus \{\ttf, \tte\ttf\}$ (six elements), and it is equal to $(\SSSs \cap \NN)\Isom{\NN} \cup \Isom{\NN}$. Hence $\SSSs \cap \NN$ generates~$\NN$. Next, $\MM$ and $\NN$ are Noetherian, and both admit right-lcms. Now a direct inspection shows that $\SSSs \cap \NN$ is (weakly) closed under right-lcm and right-divisor in~$\NN$ (here we consider only right-divisors that belong to~$\NN$, so $\ttf$, which is a right-divisor of~$\tta \ttf$ in~$\MM$, is excluded). Hence, by Corollary~IV.2.29 (recognizing Garside, right-lcm case), $\SSSs \cap \NN$ is a Garside family in~$\NN$. Finally, a pair $\tta^{\pp_1} \tte^{\ii_1} \ttf^{\jj_1} \sep \tta^{\pp_2} \tte^{\ii_2} \ttf^{\jj_2}$ is $\SSS$-normal in~$\MM$ if and only if we do not have $\pp_1 = 0$ and $\pp_2 = 1$ and the same condition characterizes $(\SSSs \cap \NN)$-normal pairs in~$\NN$. Hence $\NN$ is compatible with~$\SSSs$. \end{solu}

%%%%
\begin{exer}{77}{not closed under divisor}{Z:SUNotClosedDiv}
Let $\MM = \PRESp{\tta, \ttb}{\tta\ttb = \ttb\tta, \tta^2 = \ttb^2}$, and let $\NN$ be the submonoid of~$\MM$ generated by~$\tta^2$ and~$\tta\ttb$. Show that $\NN$ is compatible with every Garside family~$\SSS$ of~$\MM$, but that $\MM$ is not closed under left- and right-divisor.
\end{exer}

\begin{solu}
As seen in Example~IV.2.34 (no proper Garside), the only Garside family in~$\MM$ is $\MM$ itself. 
\end{solu}

%%%%
\begin{exer}{78}{head implies closed}{Z:SUHeadComult}
Assume that $\CCC$ is a left-cancellative category, $\SSSS$ is a subfamily of~$\CCC$ that is closed under right-comultiple in~$\CCC$, and $\CCCu$ is a subcategory of~$\CCC$. Put $\SSSSu = \SSSS \cap \CCCu$. \ITEM1 Show that, if every element of~$\SSSS$ admits a $\CCCu$-head that lies in~$\SSSSu$, then $\SSSSu$ is closed under right-comultiple in~$\CCC$. \ITEM2 Show that, if, moreover, $\SSSSu$ is closed under right-quotient in~$\CCC$, then $\SSSSu$ is closed under right-diamond in~$\CCC$.
\end{exer}

\begin{figure}[htb]\centering
\begin{picture}(45,36)(0,0)
\pcline{->}(1,36)(14,36)\tbput{$\tt$}
\pcline{->}(0,35)(0,25)\tlput{$\sss$}
\pcline(16,36)(38,36)\psarc(38,29){7}{0}{90}
\pcline{->}(45,29)(45,1)\trput{$\ff$}
\pcline(0,23)(0,7)\psarc(7,7){7}{180}{270}
\pcline{->}(7,0)(44,0)\taput{$\gg$}
\pcline(16,35)(23,35)\psarc(23,28){7}{0}{90}
\pcline{->}(30,28)(30,13)\trput{$\sss'$}
\pcline(1,23)(1,19)\psarc(8,19){7}{180}{270}
\pcline{->}(8,12)(29,12)\tbput{$\tt'$}
\pcline{->}(30.5,11.5)(44,1)
\pcline[style=exist]{->}(1,24)(14,24)\put(8,21){$\tt_1$}
\pcline[style=exist]{->}(15,35)(15,25)\put(15.5,28){$\sss_1$}
\pcline{->}(0.5,35.5)(14,25)\put(6,27.5){$\rr$}
\pcline{->}(15.5,23.5)(29,13)
\end{picture}
\caption[]{\sf\smaller Solution to Exercise~79}
\label{F:SUHeadComult}
\end{figure}
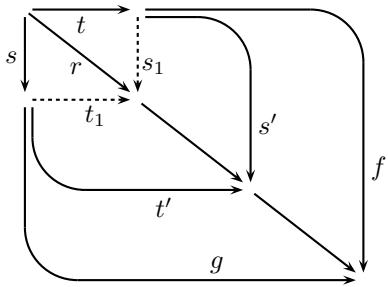

\begin{solu}
(See Figure~\ref{F:SUHeadComult}.) \ITEM1 Assume $\sss \gg = \tt \ff$ with $\sss, \tt \in \SSSSu$. As $\SSSS$ is closed under right-comultiple, there exist~$\sss', \tt'$ satisfying $\sss \tt' = \tt \sss' \dive \sss \gg$ with $\sss \tt' \in \SSSS$. Let $\rr$ be a $\CCCu$-head of~$\sss \tt'$ that lies in~$\SSSSu$. As $\sss$ lies in~$\CCCu$ and $\sss \tt'$ lies in~$\SSSS$, the relation $\sss \dive \sss\tt'$ implies $\sss \dive \rr$, so we have $\rr = \sss \tt''$ for some~$\tt''$ in~$\CCCu$. Similarly, we have $\rr = \tt \sss_1$ for some~$\sss_1$ in~$\CCCu$, and, therefore, $\rr$ witnesses that $\SSSSu$ is closed under \VR(0,1.8) right-comultiple in~$\CCC$. \ITEM2 If, in addition, $\SSSSu$ is closed under right-quotient in~$\CCCu$, then $\sss_1$ and~$\tt_1$ must lie in~$\SSSSu$ since $\sss, \tt$, and~$\rr$ do, and $\SSSSu$ is closed under right-diamond in~$\CCC$.
\end{solu}

%%%%
\begin{exer}{79}{head on generating family}{Z:SUHeadGener}
Assume that $\CCC$ is a left-cancellative category that is right-Noetherian, $\CCCu$ is a subcategory of~$\CCC$ that is closed under inverse, and $\SSSS$ is a subfamily of~$\CCC$ such that every element of~$\SSSS$ admits a $\CCCu$-head that lies in~$\SSSS$. Assume moreover that $\SSSS$ is closed under right-comultiple and that $\SSSS \cap \CCCu$ generates~$\CCCu$ and is closed under right-quotient in~$\CCC$. Show that $\CCCu$ is a head-subcategory of~$\CCC$. [Hint: Apply Exercise~78.]
\end{exer}

\begin{solu}
Exercise~78 implies that $\SSSS \cap \CCCu$ is closed under right-diamond in~$\CCC$. Then Proposition~IV.1.15 (factorization grid) implies that the subcategory of~$\CCC$ generated by~$\SSSS \cap \CCCu$, which is $\CCCu$ by assumption, is closed under right-diamond in~$\CCC$. By Lemma~VII.1.8, it follows that $\CCCu$ is closed under right-quotient in~$\CCC$. Then, as $\CCC$ is right-Noetherian, Proposition~VII.1.21 implies that $\CCCu$ is a head-subcategory of~$\CCC$.
\end{solu}

%%%%
\begin{exer}{80}{transitivity of compatibility}{Z:SUTransit}
Assume that $\SSSS$ is a Garside family in a left-cancellative category~$\CCC$, $\CCCu$ is a subcategory of~$\CCC$ that is compatible with~$\SSSS$, and $\CCC_2$ is a subcategory of~$\CCCu$ that is compatible with $\SSSSu = \SSSS \cap \CCCu$. Show that $\CCC_2$ is compatible with~$\SSSS$.
\end{exer}

\begin{solu} 
Put $\SSSS_2 = \SSSSu \cap \CCC_2$. By assumption, we have $\Isom{\CCC_2} = \CCCui\cap\CCC_2 = (\CCCi\cap\CCCu) \cap\CCC_2 = \CCCi\cap\CCC_2$. Therefore, $\CCC_2$ in closed under inverse in $\CCC$. Next, as $\CCCu$ is compatible with~$\SSSS$, the family~$\SSSSu$ is a Garside family in~$\CCCu$. Now, as $\CCC_2$ is compatible with~$\SSSSu$ in~$\CCCu$, the family $\SSSS_2$, which is $\SSSSu \cap \CCC_2$, is a Garside family in~$\CCC_2$. Finally, as $\SSSS_2 = \SSSSu \cap \CCC_2$ holds, an $\CCC_2$ path is $\SSSS_2$-normal in $\CCC_2$ if and only if it is $\SSSSu$-normal in $\CCCu$, hence if and only if it is $\SSSS$-normal in $\CCC$. Hence $\CCC_2$ is compatible with~$\SSSS$.
\end{solu}

%%%%
\begin{exer}{81}{transitivity of head-subcategory}{Z:SUHeadTransitivity}
Assume that $\CCC$ is a left-cancellative category, $\CCCu$ is a head-subcategory of $\CCC$, and $\CCC_2$ is a subcategory of~$\CCCu$. Show that $\CCC_2$ is a head-subcategory of~$\CCC$ if and only if it is a head-subcategory of~$\CCCu$. 
\end{exer}

\begin{solu} 
Assume that $\CCC_2$ is a head-subcategory of~$\CCC$. If $\ie$ belongs to~$\CCCui \cap \CCC_2$, then it belongs to~$\CCCi \cap \CCC_2$, hence $\ie\inv$ belongs to~$\CCC_2$, so $\CCC_2$ is closed under inverse in~$\CCCu$. Moreover, every element of~$\CCCu$ admits a $\CCC_2$-head since every element of~$\CCC$ does. Conversely assume that $\CCC_2$ is a head-subcategory of~$\CCCu$. If $\ie$ belongs to~$\CCCi \cap \CCC_2$, then it belongs to~$\CCCi \cap \CCCu$, hence to~$\CCCui$ since $\CCCu$ is a head-subcategory of~$\CCC$. So $\ie$ belongs to~$\CCCui \cap \CCC_2$, hence to~$\CCC_2i$ since 
$\CCC_2$ is a head-subcategory of~$\CCCu$. So $\CCC_2$ is closed under inverse in~$\CCC$. Next assume $\gg \in \CCC$. Let $\gg'$ be a $\CCCu$-head of~$\gg$ in~$\CCC$, and $\gg''$ be a $\CCC_2$-head of~$\gg'$ in~$\CCCu$. Assume $\hh \in \CCC_2$ and $\hh \dive \gg$. As $\hh$ belongs to~$\CCCu$, we must have $\hh \dive \gg'$. Then, as $\hh$ belongs to~$\CCC_2$, we must have $\hh \dive \gg''$. So $\gg''$ is a $\CCC_2$-head of~$\gg$, and $\CCC_2$ is a head-subcategory of~$\CCC$. 
\end{solu}

%%%%
\begin{exer}{82}{recognizing compatible IV}{Z:SUCompat4}
Assume that $\SSSS$ is a Garside family in a left-cancellative category~$\CCC$ and $\CCCu$ is a subcategory of~$\CCC$ that is closed under right-quotient in~$\CCC$. Show that $\CCCu$ is compatible with~$\SSSS$ if and only if, putting $\SSSSu = \SSSS \cap \CCCu$, \ITEM1 the family~$\SSSSu$ is a Garside family in~$\CCCu$, and \ITEM2 a $\CCCu$-path is strictly $\SSSSu$-normal in~$\CCCu$ if and only if it is strictly $\SSSS$-normal in~$\CCC$.
\end{exer}

\begin{solu}
Assume that $\CCCu$ is compatible with~$\SSSS$. By definition, $\SSSSu$ is a Garside family in~$\CCCu$, so \ITEM1 is satisfied. Next, assume that $\seqqq{\gg_1}\etc{\gg_\qq}$ is a strict $\SSSS$-normal $\CCCu$-path. By assumption, $\seqqq{\gg_1}\etc{\gg_\qq}$ is an $\SSSSu$-normal path in~$\CCCu$. Moreover, by definition, $\gg_1 \wdots \gg_\qq$ are not invertible in~$\CCC$, hence they are not invertible either in~$\CCCu$, and $\gg_1 \wdots \gg_{\qq-1}$ belong to~$\SSSS$ and to~$\CCCu$, hence to~$\SSSSu$. So $\seqqq{\gg_1}\etc{\gg_\qq}$ is strictly $\SSSSu$-normal in~$\CCCu$. Conversely, assume that $\seqqq{\gg_1}\etc{\gg_\qq}$ is strictly $\SSSSu$-normal in~$\CCCu$. As $\CCCu$ is compatible with~$\SSSS$, the path $\seqqq{\gg_1}\etc{\gg_\qq}$ is $\SSSS$-normal in~$\CCC$. Moreover, the assumption that $\gg_1 \wdots \gg_{\qq-1}$ belong to~$\SSSSu$ implies that they belong to~$\SSSS$. Finally, the assumption that no~$\gg_\ii$ is invertible in~$\CCCu$ implies that they are not invertible in~$\CCC$ either. So \ITEM2 is satisfied. 

Conversely, assume that $\CCCu$ satisfies~\ITEM1 and~\ITEM2. Assume that $\seqq{\gg_1}{\gg_2}$ is an $\SSSSu$-path. Then $\gg_1 \gg_2$ is invertible in~$\CCCu$ if and only if it is invertible in~$\CCC$. In this case, $\seqq{\gg_1}{\gg_2}$ is both $\SSSSu$-normal in~$\CCCu$ and $\SSSS$-normal in~$\CCC$. Otherwise, if $\gg_1 \gg_2$ belongs to~$\SSSSu \CCCui$, then $\seqq{\gg_1}{\gg_2}$ is $\SSSSu$-normal in~$\CCCu$ if and only if $\gg_2$ is invertible in~$\CCCu$ and, similarly, it is $\SSSS$-normal in~$\CCC$ if and only if $\gg_2$ is invertible in~$\CCC$, so both conditions are equivalent. Finally, if $\gg_1 \gg_2$ has $\SSSSu$-length~$2$, it admits a strict $\SSSSu$-normal decomposition in~$\CCCu$, say $\seqq{\hh_1}{\hh_2}$. Then $\seqq{\gg_1}{\gg_2}$ is $\SSSSu$-normal in~$\CCCu$ if and only if there exist~$\ie$ invertible in~$\CCCu$ satisfying $\gg_1 = \hh_1 \ie$ and $\hh_2 = \ie \gg_2$; similarly, $\seqq{\gg_1}{\gg_2}$ is $\SSSS$-normal in~$\CCC$ if and only if there exist~$\ie$ invertible in~$\CCC$ satisfying $\gg_1 = \hh_1 \ie$ and $\hh_2 = \ie \gg_2$: both conditions are equivalent again. So an $\SSSSu$-path is $\SSSSu$-normal in~$\CCCu$ if and only if it is $\SSSS$-normal in~$\CCC$, and, by definition, $\CCCu$ is compatible with~$\SSSS$. 
\end{solu}

%%%%
\begin{exer}{83}{inverse image}{Z:SUInvIm}
Assume that $\CCC, \CCC'$ are left-cancellative categories, $\phi$ is a functor from~$\CCC$ to~$\CCC'$, and $\CCCu'$ is a subcategory of~$\CCC'$ that is closed under left- and right-divisor. Show that the subcategory~$\phi\inv(\CCCu')$ is compatible with every Garside family of~$\CCC$. \ITEM2 Let $\BP{}$ be the Artin--Tits monoid of type B as defined in Example~VII.4.21. Show that the map~$\phi$ defined by~$\phi(\sig0) = 1$ and $\phi(\sig\ii) = 0$ for $\ii \ge 1$ extends into a homomorphism of~$\BP{}$ to~$\NNNN$, and that the submonoid $\NN = \{\gg\in\MM \mid\phi(\gg) = 0\}$ of~$\BP{}$ is compatible with every Garside family of~$\BP{}$. 
\end{exer}

\begin{solu}
\ITEM1 By Proposition~VII.1.18, the subcategory~$\phi\inv(\CCC'_1)$ is closed under right-quotient and under~$\eqir$. Then apply Proposition~VII.2.21. \ITEM2 Use~\ITEM1
\end{solu}

%%%%
\begin{exer}{84}{intersection}{Z:SUIntersection}
Assume that $\SSSS$ is a Garside family in a left-cancellative category~$\CCC$. \ITEM1 Let $\FFF$ be the family of all subcategories of~$\CCC$ that are closed under right-quotient, compatible with~$\SSSS$, and $\eqir$-closed. Show that every intersection of  elements of~$\FFF$ belongs to~$\FFF$. \ITEM2 Same question when ``$\eqir$-closed'' is replaced with ``including~$\CCCi$''. \ITEM3 Same question when $\CCC$ contains no nontrivial invertible element and ``$\eqir$-closed'' is skipped. 
\end{exer}

\begin{solu}
\ITEM1 Let $(\CCC_\ii)_{\ii \in \II}$ be a family of $\eqir$-closed subcategories of~$\CCC$ that are compatible with~$\SSSS$. First, an intersection of subcategories is a subcategory. Next, an intersection of $\eqir$-closed families is $\eqir$-closed and, similarly, an intersection of families that are closed under right-quotient is closed under right-quotient. Finally, let $\gg$ belong to~$\bigcap \CCC_\ii$. Then $\gg$ admits an $\SSSS$-normal decomposition $\sss_1 \sep \etc \sep \sss_\pp$. By Lemma~VII.2.20, the latter is $(\SSSS \cap \CCC_\ii)$-normal for every~$\ii$, hence it its entries lie in every subfamily~$\CCC_\ii$, hence in their intersection. By Proposition~VII.2.21 we deduce that $\bigcap\CCC_\ii$ is compatible with~$\SSSS$. 

\ITEM2 An intersection of families that include~$\CCCi$ includes~$\CCCi$. On the other hand, a subcategory that includes~$\CCCi$ is $\eqir$-closed and we apply~\ITEM1. 

\ITEM3 In this case, every subcategory is $\eqir$-closed, and we apply~\ITEM1.
\end{solu}

%%%%
\begin{exer}{85}{fixed points}{Z:SUWhole}
Assume that $\CCC$ is a left-cancellative category and $\phi: \CCC\to \CCC$ is a functor. Show that the fixed point subcategory~$\CCC^\phi$ is compatible with~$\CCC$ viewed as a Garside family in itself.
\end{exer}

\begin{solu}
First, assume that $\ie$ belongs to~$\CCC(\xx, \yy) \cap \CCC^\phi$. Necessarily we have $\phi(\xx) = \xx$. As $\phi$ is a functor, we find $\ie \ie\inv = \id\xx = \phi(\id\xx) = \phi(\ie \ie\inv) = \phi(\ie) \phi(\ie\inv) = \ie \phi(\ie\inv)$, whence $\ie\inv = \phi(\ie\inv)$, so $\CCCi \cap \CCC^\phi \subseteq \Isom{(\CCC^\phi)}$ holds. Next, we have $\sh{\CCC} = \CCC$, so $\sh{\CCC} \cap \CCC^\phi \subseteq \SSSS^\phi$ trivially holds, and $\CCC^\phi$ and~$\CCC$ satisfy (the counterpart of)~(VII.2.11). On the other hand, every element~$\gg$ of~$\CCC$ admits the $\CCC$-normal decomposition~$\seq\gg$, so, in particular, every element of~$\CCC^\phi$ has a $\CCC$-normal decomposition whose entries lie in $\CCC^\phi$. So (the counterpart of)~(VII.2.12) is satisfied and, by Proposition~VII.2.10, $\CCC^\phi$ is compatible with~$\CCC$.
\end{solu}

%%%%
\begin{exer}{86}{connection between closure properties}{Z:SUClosureConn}
Assume that $\SSSS$ is a subfamily in a left-cancellative category~$\CCC$ and $\SSSSu$ is a subfamily of~$\SSSS$. \ITEM1 Show that, if $\SSSSu$ is closed under product, inverse, and right-complement in~$\SSSS$, then $\SSSSu$ is closed under right-quotient in~$\SSSS$.
\ITEM2 Assume that $\SSSSu$ is closed under product and right-complement in~$\SSSS$. Show that, if $\SSSS$ is closed under left-divisor in~$\CCC$, then $\SSSSu$ is closed under right-comultiple in~$\SSSS$. Show that, if $\SSSS$ is closed under right-diamond in~$\CCC$, then $\Sub\SSSSu$ is closed under right-comultiple in~$\SSSS$.
\ITEM3 Show that, if $\SSSSu$ is closed under identity and product in~$\SSSS$, then $\SSSSu$ is closed under inverse and right-diamond in~$\SSSS$ if and only if $\SSSSu$ is closed under right-quotient and right-comultiple in~$\SSSS$. 
\end{exer}

\begin{solu}
\ITEM1 The argument is the same as for Lemma~VII.1.8. Assume $\sss \gg = \tt$ with $\sss, \tt$ in~$\SSSSu$ and $\gg$ in~$\SSSS$. As $\SSSSu$ is closed under right-complement in~$\SSSS$, there exist $\sss', \tt'$ in~$\SSSSu$ and $\rr$ in~$\SSSS$ satisfying $\sss \tt' = \tt \sss'$, $\id\yy = \sss' \rr$, and $\gg = \tt' \rr$, where $\yy$ is the target of~$\tt$. The second equality implies that $\sss'$ is invertible, and the assumption that $\SSSSu$ is closed under inverse in~$\SSSS$ then implies that $\sss'{}\inv$, that is, $\rr$, lies in~$\SSSSu$.  Hence  $\tt' \rr$, that is~$\gg$, belongs to~$\Seq\SSSSu2 \cap \SSSS$, hence to~$\SSSSu$ as the latter is closed under product in~$\SSSS$. So $\SSSSu$ is closed under right-quotient in~$\SSSS$.

\ITEM3 Assume that $\SSSSu$ is closed under inverse and right-diamond in~$\SSSS$. Then, by definition, $\SSSSu$ is closed under right-comultiple and right-complement, and, by Exercice~89 (transfer of closure), it is closed under right-quotient. 

Conversely, assume that $\SSSSu$ is closed under right-quotient and right-comultiple in~$\SSSS$. First, the assumption that $\SSSSu$ is closed under right-quotient trivially implies that $\SSSSu$ is closed under inverse. Next, by the same argument as in Lemma~IV.1.8, the assumption that $\SSSSu$ is closed under right-comultiple and right-quotient implies that it is closed under right-diamond. 
\end{solu}

%%%%
\begin{exer}{87}{subgerm}{Z:SUSubgerm}
Assume that $\SSSSg$ is a left-cancellative germ and $\SSSSgu$ is a subgerm of~$\SSSS$ such that the relation~$\diveSu$ is the restriction to~$\SSSSu$ of the relation~$\diveS$. Show that $\SSSSu$ is closed under right-quotient in~$\SSSS$. 
\end{exer}

\begin{solu}
Assume $\ff \OP \hh = \gg$ with $\ff, \gg \in \SSSSu$. By definition, $\ff \diveS \gg$ holds, hence so does $\ff \diveSu \gg$. This means that there exists~$\hh_1$ in~$\SSSSu$ satisfying $\ff \OP \hh_1 = \gg$. The assumption that $\SSSS$ is left-cancellative implies $\hh_1 = \hh$, whence $\hh \in \SSSS$. So $\SSSSu$ is closed under right-quotient in~$\SSSS$.
\end{solu}

%%%%
\begin{exer}{88}{transitivity of closure}{Z:SURelClosConverse}
Assume that $\SSSSgu$ is a subgerm of a Garside germ~$\SSSSg$, the subcategory~$\Sub\SSSSu$ is closed under right-quotient and right-diamond in $\Cat{\SSSSg}$, and $\SSSSu$ is closed under inverse and right-complement (\resp right-diamond) in~$\Sub\SSSSu$. Show that $\SSSSu$ is closed under inverse and right-comple\-ment (\resp right-diamond) in~$\SSSS$.
\end{exer}

\begin{solu}
By Lemma~VII.1.7, the assumption that $\Sub\SSSSu$ is closed under right-quotient implies that it is closed under inverse, and Lemma~VII.3.7, which is valid since the assumption that $\SSSSg$ is a Garside germ implies that $\SSSS$ is closed under right-quotient in~$\Cat\SSSSg$, then implies that $\SSSSu$ is closed under inverse in~$\SSSS$. Next, by Lemma~VII.3.8 applied with~$\Sub\SSSSu$ in place of~$\SSSS$, the assumption that $\SSSSu$ is closed under right-complement (\resp right-diamond) in~$\Sub\SSSSu$ implies that $\SSSSu$ is closed under right-complement (\resp right-diamond) in~$\Cat\SSSSg$. Now, as $\SSSS$ is a solid Garside family in~$\Cat\SSSSg$, it is closed under right-divisor in~$\Cat\SSSSg$, hence a fortiori under right-quotient. Applying once more Lemma~VII.3.8, we deduce from the fact that $\SSSSu$ is closed under right-complement (\resp right-diamond) in~$\Cat\SSSSg$ that $\SSSSu$ is closed under right-complement (\resp right-diamond) in~$\SSSS$. 
\end{solu}

%%%%
\begin{exer}{89}{transfer of closure}{Z:SUTransRQ}
Assume that $\CCC$ is a left-cancellative category, $\SSSS$ is a subfamily of~$\CCC$, and $\SSSSu$ is a subfamily of~$\SSSS$ that is closed under identity and product. \ITEM1 Show that, if $\Sub\SSSSu$ is closed under right-quotient in~$\CCC$, then $\SSSSu$ is closed under right-quotient in~$\SSSS$. \ITEM2 Show that, if, moreover, $\SSSS$ is closed under right-quotient in~$\CCC$, then $\SSSSu$ is closed under right-quotient in~$\Sub\SSSSu$.
\end{exer}

\begin{solu}
\ITEM1 Assume that $\tt$ belongs to~$\SSSS$ and $\sss$ and $\sss \tt$ belong to~$\SSSSu$. Then $\sss$ and $\sss \tt$ belong to~$\Sub\SSSSu$, so, as the latter is assumed to be closed under right-quotient, $\tt$ must belong to~$\Sub\SSSSu$, hence to $\Sub\SSSSu \cap \SSSS$. By Proposition~VII.3.3, the latter is~$\SSSSu$. So $\SSSSu$ is closed under right-quotient in~$\SSSS$.

\ITEM2 Assume that $\SSSS$ is closed under right-quotient in~$\CCC$. Assume that $\tt$ belongs to~$\Sub\SSSSu$ and $\sss$ and $\sss \tt$ belong to~$\SSSSu$. Then $\sss$ and $\sss \tt$ belong to~$\SSSS$, so the assumption implies that $\tt$ belong to~$\SSSS$, hence to $\Sub\SSSSu \cap \SSSS$, which is~$\SSSSu$ by Proposition~VII.3.3. So $\SSSSu$ is closed under right-quotient in~$\Sub\SSSSu$.
\end{solu}

%%%%
\begin{exer}{90}{braid subgerm}{Z:SUSubgerm1}
Let $\SSSSg$ be the six-element Garside germ associated with the divisors of~$\Delta_3$ in the braid monoid~$\BP3$. \ITEM1 Describe the subgerm~$\SSSSgu$ of~$\SSSSg$ generated by~$\sig1$ and~$\sig2$. Compare $\Mon\SSSSgu$ and~$\Sub\SSSSu$ (describe them explicitly). \ITEM2 Same questions with $\sig1$ and~$\sig2\sig1$. Is $\SSSSu$ closed under right-quotient in~$\SSSSg$ in this case? 
\end{exer}

\begin{solu}
\ITEM1 The subgerm of~$\SSSSg$ generated by~$\sig1$ and~$\sig2$ is the closure of~$\{\sig1, \sig2\}$ under identity and product in~$\SSSSg$: so it contains~$1$, and, as $\sig1 \OP \sig2$ is defined in~$\SSSSg$, it contains $\sig1 \OP \sig2$, that is, $\sig1\sig2$. Then it contains $\sig1\sig2 \OP \sig1$, which is~$\Delta_3$, etc. Finally, one finds $\SSSSgu = \SSSSg$, whence $\Mon\SSSSgu = \Sub\SSSSu = \BP3$. 

\ITEM2 The closure of~$\{\sig1, \sig2\sig1\}$ under identity and product is $\{1, \sig1, \sig2\sig1, \Delta_3\}$. In the associated germ~$\SSSSgu$, the only nontrivial product is $\sig1 \OPu \sig2\sig1 = \Delta_3$, so $\Mon\SSSSgu$ is $\PRESp{\tta, \ttb, \Delta_3}{\tta\ttb = \Delta_3}$, that is, a free monoid based on~$\tta$ and~$\ttb$. On the other hand, in~$\Sub\SSSSu$, the relation $(\sig2\sig1)^3 = \Delta_3^2$ holds, corresponding to $\ttb^3 = (\tta\ttb)^2$, which fails in~$\Mon\SSSSgu$: so $\Mon\SSSSgu$ is not isomorphic to~$\Sub\SSSSu$. Here $\Sub\SSSSu$ is not closed under right-quotient in~$\SSSSg$, as $\Delta_3$ and~$\sig2\sig1$ lie in~$\Sub\SSSSu$, but we have $\Delta_3 = \sig2\sig1 \OP \sig2$ in~$\SSSSg$ and $\sig2 \notin \SSSSu$. 
\end{solu}

%%%%
\begin{exer}{91}{$\eqir$-closed}{Z:SUEqirClosed}
Show that, if $\SSSSgu$ is a subgerm of an associative germ~$\SSSSg$, then $\Sub\SSSSu$ is $\eqir$-closed in~$\Cat\SSSSg$ if and only if $\SSSSgu$ is $\eqirS$-closed in~$\SSSS$.
\end{exer}

\begin{solu}
Assume that $\SSSSgu$ is  an $\eqirS$-closed  subgerm of~$\SSSSg$. Assume that $\gg$ is an element of~$\Sub\SSSSu$ and $\gg' \eqir \gg$ holds in~$\Cat\SSSSg$. By definition, $\gg$ admits an $\SSSSu$-decomposition, say $\seqqq{\sss_1}\etc{\sss_\pp}$, and, in~$\Cat\SSSSg$, we then have $\gg' = \sss_1 \pdots \sss_\pp \ie$. As $\SSSS$ is closed under left-divisor in~$\Cat\SSSSg$ since $\SSSSg$ is right-associative, the assumption that $\sss_\pp$ lies in~$\SSSS$ implies that its left-divisor~$\sss_\pp \ie$ also lies in~$\SSSS$, that is, $\sss_\pp \OP \ie$ is defined in~$\SSSSg$. The assumption that $\SSSSgu$ is  $\eqirS$-closed  in~$\SSSSg$ then implies that $\sss_\pp \OP \ie$ belongs to~$\SSSSu$, which implies that $\gg'$, which is $\sss_1 \pdots \sss_{\pp-1} (\sss_\pp\ie)$, lies in~$\Sub\SSSSu$.
Conversely, assume that $\Sub\SSSSu$ is  an $\eqir$-closed  subcategory of~$\Cat\SSSSg$. Assume that $\sss$ belongs to~$\SSSSu$ and $\sss' \eqirS \sss$ holds in~$\SSSSg$. This means that there exists~$\ie$ in~$\SSSSi$ satisfying $\sss' = \sss \OP \ie$. Then, in~$\Cat\SSSSg$, we have $\sss' = \sss\ie$, whence $\sss' \eqir \sss$. The assumption that $\Sub\SSSSu$ is  $\eqir$-closed  implies $\sss' \in \Sub\SSSSu$, whence $\sss' \in \Sub\SSSSu \cap \SSSS$. As $\SSSSgu$ is a subgerm of~$\SSSSg$, the latter family is~$\SSSSu$. 
\end{solu}

%%%%
\begin{exer}{92}{correct \vs mcms}{Z:SUCorrectMcm}
Assume that $\CCC$ and~$\CCC'$ are left-cancellative categories and $\SSSS$ is included in~$\CCC$. Assume moreover that $\CCC$ and $\CCC'$ admit right-mcms and $\SSSS$ is closed under right-mcm. Show that a functor~$\phi$ from~$\CCC$ to~$\CCC'$ is correct for right-comultiples on~$\SSSS$ if and only if, for all~$\sss, \tt$ in~$\SSSS$, every right-mcm of~$\phi(\sss)$ and~$\phi(\tt)$ is $\eqir$-equivalent to the image under~$\phi$ of a right-mcm of~$\sss$ and~$\tt$.
\end{exer}

\begin{solu}
Assume that $\phi$ is correct for right-comultiples on~$\SSSS$, that $\sss, \tt$ lie in~$\SSSS$, and that $\hh$ is a right-mcm of~$\phi(\sss)$ and~$\phi(\tt)$. By definition, there exists a common right-multiple~$\rr$ of~$\sss$ and~$\tt$ that lies in~$\SSSS$ and satisfies $\phi(\rr) \dive \hh$. As $\CCC$ admits right-mcms, there exists a right-mcm~$\rr'$ of~$\sss$ and~$\tt$ satisfying $\rr' \dive \rr$ and, as $\SSSS$ is closed under right-mcm, $\rr'$ lies in~$\SSSS$. Then $\phi(\rr') \dive \hh$ holds, and $\phi(\rr')$ is a common right-multiple of~$\phi(\sss)$ and~$\phi(\tt)$. As $\rr$ is minimal, we deduce $\phi(\rr') \eqir \hh$. So the condition is necessary.
Conversely, assume that $\phi$ preserves mcms in the sense of the statement. Assume that $\sss, \tt$ lie in~$\SSSS$, and $\hh$ is a common right-multiple of~$\phi(\sss)$ and~$\phi(\tt)$. As $\CCC'$ admits right-mcms, $\hh$ is a right-multiple of some right-mcm of~$\phi(\sss)$ and~$\phi(\tt)$, hence, by~\ITEM2, of some element~$\phi(\rr)$ where $\rr$ is a right-mcm of~$\sss$ and~$\tt$. By assumption, $\rr$ belongs to~$\SSSS$, and there exist~$\sss', \tt'$ in~$\CCC$ satisfying $\sss \tt' = \tt \sss' = \rr$. Hence $\phi$ is correct for right-comultiples on~$\SSSS$.
\end{solu}

%%%%
\chapter{Chapter~VIII: Conjugacy}

%%%%
\section*{Skipped proofs}

\noindent(none)

%%%%
\section*{Solution to selected exercises}

%%%%
\begin{exer}{97}{quasi-distance}{Z:COQuasiDist}
\ITEM1 In Context~VIII.3.7, show that $\CAN\Delta(\gg) = \SUP\Delta(\rep\gg)$ holds for every~$\gg$ in~$\GGG$. \ITEM2 For $\gg, \gg'$ in~$\GGG$, define $\dist(\gg, \gg')$ to be $\infty$ if $\gg, \gg'$ do not share the same source, and to be $\CAN\Delta(\gg\inv \gg')$ otherwise. Show that $\dist$ is a quasi-distance on~$\GGG$ that is compatible with~$\eqD$.
\end{exer}

\begin{solu}
\ITEM2 The canonical length is invariant under left- and right-multiplication by~$\Delta$ and, therefore, $\dist$ takes constant values on $\eqD$-classes. 
\end{solu}

%%%%
\chapter{Chapter~IX: Braids}

%%%%
\section*{Skipped proofs}

\noindent(none)

%%%%
\section*{Solution to selected exercises}

%%%%
\begin{exer}{102}{smallest Garside, right-angled type}{Z:BRRightAngled}
Assume that $\BP{}$ is a right-angled Artin--Tits monoid, that is, $\BP{}$ is associated with a Coxeter system~$(\WW, \Gen)$ satisfying $\mm_{\sss, \tt} \in \{2, \infty\}$ for all~$\sss, \tt$ in~$\Gen$. \ITEM1 For~$\II \subseteq \Gen$, denote by~$\Delta_\II$ the right-lcm (here the product) of the elements of~$\II$, when it exists (that is, when the elements pairwise commute). Show that the divisors of~$\Delta_\II$ are the elements~$\Delta_\JJ$ with $\JJ \subseteq \II$. \ITEM2 Deduce that the smallest Garside family in~$\BP{}$ is finite and consists of the elements~$\Delta_\II$ for $\II \subseteq \Gen$. 
\end{exer}

\begin{solu}
\ITEM1 Clearly $\JJ \subseteq \II$ implies $\Delta_\JJ \dive \Delta_\II$ and $\Delta_\JJ \diveR \Delta_\II$, since the elements of~$\II$ commute. Conversely, the point is to see that, if $\uu$ is a word in~$\SSS$ that involves at least one letter not in~$\II$, then (the class of)~$\uu$ cannot right-divide~$\Delta_\II$: this is so since, when left-reversing~$\Delta_\II \INV\uu$, a negative letter~$\INV\sss$ can only vanish when it is adjacent to the positive letter~$\sss$, and it, if $\uu$ is a word in~$\II$ in which some letter is repeated twice, then the second occurrence cannot vanish as we are arguing in the free abelian monoid based on~$\II$. 

\ITEM2 Let $\SSS$ be the family of all elements~$\Delta_\II$ with $\II$ a set of pairwise commuting atoms in~$\BP{}$. As the smallest Garside family of~$\BP{}$ includes~$\II$ and is closed under right-lcm, it must include~$\SSS$. On the other hand, $\SSS$ includes~$\Gen$, it is closed under right-lcm by definition, and it is closed under right-divisor by~\ITEM1. Corollary~IV.2.29 (recognizing Garside, right-lcm case) implies that $\SSS$ is a Garside family in~$\BP{}$.
\end{solu}

%%%%
\begin{exer}{103}{smallest Garside, large type}{Z:BRLarge}
Assume that $\BP{}$ is an Artin--Tits monoid of large type, $\BP{}$ is associated with a Coxeter system~$(\WW, \Gen)$ satisfying $\mm_{\sss, \tt} \ge 3$ for all~$\sss, \tt$ in~$\Gen$. Put 

$\Gen_1 = \{\sss {\in} \Gen \mid \forall \rr{\in} \Gen \, (\mm_{\rr, \sss} {=} \infty)\}$,  

$\Gen_2 = \{(\sss, \tt) {\in} \Gen^2 \mid \mm_{\sss, \tt} <\nobreak \infty \mbox{\ and\ }\forall \rr{\in} \Gen\, (\mm_{\rr, \sss} {+} \mm_{\rr, \tt} {=} \infty\}$, 

$\Gen_3 = \{(\rr, \sss, \tt) {\in} \Gen^3 \mid \mm_{\rr, \sss} {+} \mm_{\rr, \tt} {+} \mm_{\sss, \tt} {<} \infty\}$, 

\noindent and $\EE = \Gen_1 \cup \{\Delta_{\sss, \tt} \mid (\sss, \tt) \in \Gen_2\} \cup \{\rr\Delta_{\sss, \tt} \mid (\rr, \sss, \tt) \in \Gen_3\}$ (we write $\Delta_{\sss, \tt}$ for the right-lcm of~$\sss$ and~$\tt$ when it exists). \ITEM1 Explicitly describe the elements of the closure~$\SSS$ of~$\EE$ under right-divisor, and deduce that $\SSS$ is finite. \ITEM2 Show that $\SSS$ is closed under right-lcm and deduce that $\SSS$ is a Garside family in~$\BP{}$. \ITEM3 Show that $\SSS$ is the smallest Garside family in~$\BP{}$. \ITEM4 Show that, if $\Gen$ has $\nn$~elements and $\mm_{\sss, \tt} \not= \infty$ holds for all~$\sss, \tt$, then $\EE$ has $3{\nn\choose3}$ elements. \ITEM5 Show that, if $\Gen$ has $\nn$~elements and $\mm_{\sss, \tt} = \mm$ holds for all~$\sss, \tt$, then $\SSS$ has $(\nn + 2\mm-5){\nn\choose2} + \nn + 1$ elements. Apply to $\nn = \mm = 3$.
\end{exer}

\begin{solu}
\ITEM1 For $\sss$ in~$\Gen_1$, the only right-divisors of~$\sss$ are~$1$ and~$\sss$; for $(\sss, \tt)$ in~$\Gen_2$, the right-divisors of~$\Delta_{\sss, \tt}$ are the elements~$\Prod\sss\tt\kk$ and~$\Prod\tt\sss\kk$ with $\kk \le \mm_{\sss, \tt}$: this is so because right-divisors are detected using left-reversing (Section~II.4) and it is clear that, if $\uu$ is a word involving a letter different from~$\sss$ and~$\tt$ or if $\uu$ involves two letters~$\sss$ or two letters~$\tt$, or if $\uu$ is of the form $\Prod\sss\tt\kk$ and~$\Prod\tt\sss\kk$ with $\kk > \mm_{\sss, \tt}$, then left-reversing $\Delta_{\sss, \tt} \uu\inv$ cannot result in a positive word. Finally, for $(\rr, \sss, \tt)$ in~$\Gen_3$, the right-divisors of~$\rr\Delta_{\sss, \tt}$ are $\rr\Delta_{\sss, \tt}$ and the elements~$\Prod\sss\tt\kk$ and~$\Prod\tt\sss\kk$ with $\kk \le \mm_{\sss, \tt}$: again, this is so because left-reversing $\Delta_{\sss, \tt} \uu\inv$ can result in a positive word only in the cases described above, and it can result in~$\rr\inv$ concatenated with a positive word only if $\uu$ is~$\rr\Delta_{\sss, \tt}$ itself.

\ITEM2 We have to analyze when two elements of~$\SSS$ may admit a common right-multiple. As there are three cases for each factor, nine cases are to be considered. Many of them are trivial. First, all cases involving an atom of~$\Gen_1$ are trivial. Then, in the case of two elements $\Prod\sss\tt\kk$ and $\Prod{\sss'}{\tt'}{\kk'}$, the case $\{\sss, \tt\} = \{\sss', \tt'\}$ is obvious. Otherwise, for $\sss' = \sss$ and $\tt' \not= \tt$, the only cases when a common right-multiple may exist are the trivial cases $\kk = 1$ or $\kk' = 1$, plus the case $\kk = \kk' = 2$, in which case there exists a common right-multiple for $\mm_{\tt, \tt'} < \infty$, and the right-lcm is then $\sss \Delta_{\tt, \tt'}$, an element of~$\SSS$. Finally, for $\card\{\sss, \tt, \sss', \tt'\} = 4$, no common right-multiple may exist for $\kk, \kk' \ge 2$, and the remaining cases are treated as above. The cases of an element $\rr\Delta_{\sss, \tt}$ and an element $\Prod{\sss'}{\tt'}{\kk'}$ and of two elements $\rr\Delta_{\sss, \tt}$ and $\rr'\Delta_{\sss', \tt'}$ are treated in the same way: once again, the point is that common right-multiples may exist only in the obvious cases. Finally, $\SSS$ is closed under right-lcm. As, by definition, $\SSS$ includes~$\Gen$  and is closed under right-divisor, Corollary~IV.2.29 (recognizing Garside, right-lcm case) implies that $\SSS$ is a Garside family in~$\BP{}$.

\ITEM3 Conversely, every Garside family~$\SSS'$ of~$\BP{}$ containing~$1$ must include~$\SSS$. Indeed, $\SSS'$ must include~$\Gen$, hence all elements~$\Delta_{\sss, \tt}$ since it is closed under right-lcm. Moreover, for $(\rr, \sss, \tt)$ in~$\Gen_3$, the family~$\SSS'$ must contain $\rr\sss$ since it is closed under right-divisor and $\rr\sss$ right-divides~$\Delta_{\rr, \sss}$; it contains $\rr\tt$ for a similar reason and, therefore, it contains the right-lcm $\rr\Delta_{\sss, \tt}$ of~$\rr\sss$ and~$\rr\tt$. So $\SSS'$ includes~$\EE$, hence~$\SSS$, and $\SSS$ is the smallest Garside family containing~$1$ in~$\BP{}$.

\ITEM4 When no coefficient~$\mm_{\sss, \tt}$ is~$\infty$, $\Gen_1$ and~$\Gen_2$ are empty and $\EE$ consists of all elements~$\rr\Delta_{\sss, \tt}$ with $\rr, \sss, \tt$ pairwise distinct atoms in~$\BP{}$. As $\mm_{\sss, \tt}$ does not matter, there exist $3{\nn\choose 3}$ such elements for $\Gen$ of size~$\nn$.

\ITEM5 By the description of~\ITEM1, $\SSS$ comprises $1$, plus the $\nn$ atoms, plus, for each of the $\nn\choose2$ pairs~$(\sss, \tt)$, the $2\mm-3$ right-divisors of~$\Delta_{\sss, \tt}$ of length~$\ge 2$, plus, for each of the $\nn\choose3$ triples~$(\rr, \sss, \tt)$, the $3$ elements $\rr \Delta_{\sss, \tt}, \sss\Delta_{\tt, \rr}, \tt\Delta_{\rr, \sss}$, whence the formula. The case $\nn = \mm = 3$ corresponds to the affine type~$\wAt$, and confirms that the smallest Garside family containing~$1$ has $16$~elements, as seen in Reference Structure~9, page~111.
\end{solu}

%%%%
\chapter{Chapter~X: Deligne--Lusztig varieties}

%%%%
\section*{Skipped proofs}

\noindent(none)

%%%%
\section*{Solution to selected exercises}

\noindent(none)

%%%%
\chapter{Chapter~XI: Left self-distributivity}

%%%%
\section*{Skipped proofs}

\begin{ADlemm}{XI.1.8}{L:LDCancel}
Assume that $\LDAct$ is a partial action of a monoid~$\MM$ on a set~$\XX$. 

\ITEM1 If the monoid~$\MM$ is left-cancellative, then so is the category~$\LDCAT\MM\LDAct$.

\ITEM2 Conversely, if $\LDAct$ is proper and $\LDCAT\MM\LDAct$ is left-cancellative, then so is~$\MM$.
\end{ADlemm}

\begin{proof}
\ITEM1 Assume that $\MM$ is left-cancellative and, in~$\LDCAT\MM\LDAct$, the equality $\NT(\xx, \gg, \yy) \cdot \NT(\yy, \hh, \zz) = \NT(\xx, \gg, \yy) \cdot \NT(\yy, \hh', \zz')$ holds. This implies $\NT(\xx, \gg \hh, \zz) = \NT(\xx, \gg \hh', \zz')$, whence, by definition, $\gg \hh = \gg \hh'$ and~$\zz = \zz'$, and $\hh = \hh'$ since $\MM$ is left-cancellative. We deduce $\NT(\yy, \hh, \zz) = \NT(\yy, \hh', \zz')$, and $\LDCAT\MM\LDAct$ is left-cancellative. 

\ITEM2 Conversely, assume that $\LDAct$ is proper and $\LDCAT\MM\LDAct$ is left-cancellative. Assume $\gg \hh = \gg \hh'$. As $\LDAct$ is proper, there exists~$\xx$ in~$\XX$ such that both $\xx \act \gg\hh$ and $\xx \act \gg \hh'$ are defined. Let $\zz = \xx \act \gg \hh$. By~(XI.1.3), $\xx \act \gg$ must be defined and, putting $\yy = \xx \act \gg$, we have $\zz = \yy \act \hh = \yy \act \hh'$. Then, in~$\LDCAT\MM\LDAct$, we have $\NT(\xx, \gg, \yy) \cdot \NT(\yy, \hh, \zz) = \NT(\xx, \gg, \yy)\cdot \NT(\yy, \hh', \zz')$. As $\LDCAT\MM\LDAct$ is left-cancellative, we deduce $\NT(\yy, \hh, \zz) = \NT(\yy, \hh', \zz')$, whence $\hh = \hh'$. Hence $\MM$ is left-cancellative.
\end{proof}

%%%%
\begin{ADlemm}{XI.1.9}{L:LDDiv}
Assume that $\LDAct$ is a partial action of a monoid~$\MM$ on a set~$\XX$. 

\ITEM1 Assume that $\xx \act \gg$ is defined. Then $\NT(\yy, \hh, \ud) \dive \NT(\xx, \gg, \ud)$ holds in~$\LDCAT\MM\LDAct$ if and only if we have $\yy = \xx$ and $\hh \dive \gg$ in~$\MM$. 

\ITEM2 Assume that $\xx \act \ff$ and $\xx \act \gg$ are defined. Then $\NT(\xx, \hh, \ud)$ is a left-gcd of $\NT(\xx, \ff, \ud)$ and $\NT(\xx, \gg, \ud)$ in~$\LDCAT\MM\LDAct$ if and only if $\hh$ is a left-gcd of~$\ff$ and~$\gg$ in~$\MM$.
\end{ADlemm}

\begin{proof}
\ITEM1 Assume $\NT(\yy, \hh, \yy\act\hh) \cdot \NT(\xx', \gg', \ud) = (\xx, \gg, \ud)$. Then we have $\yy = \xx$ and $\hh\gg' =\nobreak \gg$, hence $\hh \dive \gg$. Conversely, assume $\hh \dive \gg$, say $\hh \gg' = \gg$. By~(XI.1.3), the assumption that $\xx \act \gg$ is defined implies that $\xx \act \hh$ is defined, say $\xx \act \hh = \xx'$. Then we have $\NT(\xx, \hh, \xx') \cdot \NT(\xx', \gg', \ud) = \NT(\xx, \gg, \ud)$, whence $\NT(\xx, \hh, \ud) \dive \NT(\xx, \gg, \ud)$ in~$\LDCAT\MM\LDAct$.

\ITEM2 Assume that $\NT(\xx, \hh, \ud)$ is a left-gcd of $\NT(\xx, \ff, \ud)$ and $\NT(\xx, \gg, \ud)$ in~$\LDCAT\MM\LDAct$. By~\ITEM1, $\hh$ left-divides~$\ff$ and~$\gg$ in~$\MM$. Let $\hh'$ be a common left-divisor of~$\ff$ and~$\gg$ in~$\MM$. By~(XI.1.3), the assumption that $\xx \act \ff$ is defined implies that $\xx \act \hh'$ is defined. Then $\NT(\xx, \hh, \ud)$ is a common left-divisor of $\NT(\xx, \ff, \ud)$ and $\NT(\xx, \gg, \ud)$ in~$\LDCAT\MM\LDAct$, hence it is a left-divisor of $\NT(\xx, \hh, \ud)$. By~\ITEM1, this implies $\hh' \dive \hh$. So $\hh$ is a left-gcd of~$\ff$ and~$\gg$. 

Conversely, assume that $\hh$ is a left-gcd of~$\ff$ and~$\gg$. By~\ITEM1, $\NT(\xx, \hh, \ud)$ left-divides $\NT(\xx, \ff, \ud)$ and $\NT(\xx, \gg, \ud)$ in~$\LDCAT\MM\LDAct$. Now, consider an arbitrary common left-divisor of $\NF(\xx, \ff, \ud)$ and $\NT(\xx, \gg, \ud)$. By~\ITEM1, it must be of the form $\NT(\xx, \hh', \ud)$ with $\hh'$ left-dividing~$\ff$ and~$\gg$. Then $\hh$ left-divides~$\hh$ in~$\MM$, and $\NT(\xx, \hh', \ud)$ left-divides $\NT(\xx, \hh, \ud)$ in~$\LDCAT\MM\LDAct$.
\end{proof}

%%%%
\begin{ADlemm}{XI.4.16}{L:LDProjCompl2}
Assume that $\MM$ is a left-cancellative monoid, $\LDAct$ is a proper partial action of~$\MM$ on some set~$\XX$, and $(\Delta_\xx)_{\xx \in \XX}$ is a right-Garside map on~$\MM$ with respect to~$\LDAct$. Assume moreover that $\SSS$ is a family of atoms that generate~$\MM$. Now assume that $\pi: \MM \to \und\MM$ is a surjective homomorphism, $\LDpibis : \XX \to \und\XX$ is a surjective map, and, for all~$\xx$ in~$\XX$ and $\gg$ in~$\MM$,
\begin{gather}
\ADlabel{XI.4.17}
\parbox{100mm}{The value of~$\LDpibis(\xx \act \gg)$ only depends on~$\LDpibis(\xx)$ and~$\pi(\gg)$;}\\ 
\ADlabel{XI.4.18}
\parbox{100mm}{The value of~$\pi(\Delta_\xx)$ only depends on~$\LDpibis(\xx)$.} 
\end{gather}
Assume finally that $\LDsec: \und\SSS \to \SSS$ is a section of~$\pi$ such that, for~$\xx$ in~$\XX$, $\und\sss$ in~$\und\SSS$, and $\und\ww$ in~$\und\SSS{}^*$,
\begin{gather}
\ADlabel{XI.4.19}
\parbox{100mm}{If $\LDpibis(\xx) \act \und\sss$ is defined, then so is $\xx \act \LDsec(\und\sss)$.}\\
\ADlabel{XI.4.20}
\parbox{100mm}{The relation  $\cl{\und\ww} \dive \Delta_{\LDpibis(\xx)}$ implies $\cl{\LDsec^*(\und\ww)} \dive \Delta_\xx$.}
\end{gather}
Define a partial action~$\und\LDAct$ of~$\und\MM$ on~$\und\XX$ by $\LDpibis(\xx) \act \pi(\gg) = \LDpibis(\xx \act \gg)$, and, for~$\und\xx$ in~$\und\XX$, let $\Delta_{\und\xx}$ be the common value of~$\pi(\Delta_\xx)$ for $\xx$ satisfying $\LDpibis(\xx) = \und\xx$. Then $\und\LDAct$ is a proper partial action of~$\und\MM$ on~$\und\XX$ and $(\Delta_{\und\xx})_{\und\xx \in \und\XX}$ is a right-Garside sequence on~$\und\MM$.
\end{ADlemm}

\begin{proof}
First, (XI.4.17) and (XI.4.18) guarantee that the definitions of~$\und\xx \act \und\gg$ and~$\Delta_{\und\xx}$ are unambiguous. Next $\und\LDAct$ is a partial action of~$\und\MM$ on~$\und\XX$. Indeed, if $\und\gg$ and~$\und\hh$ lie in~$\und\MM$, and if $\gg, \hh$ satisfy $\pi(\gg) = \und\gg$ and $\pi(\hh) = \und\hh$, then we have $\pi(\gg \hh) = \und\gg \und\hh$ and, for~$\und\xx$ in~$\und\XX$ satisfying $\und\xx = \LDpibis(\xx)$, we can write
$$(\und\xx \act \und\gg) \act \und\hh = \LDpibis(\xx \act \gg) \act \pi(\hh) = \LDpibis((\xx \act \gg) \act \hh) = \LDpibis(\xx \act \gg \hh) = \und\xx \act \und\gg \und\hh,$$
equality meaning as usual that the involved expressions are simultaneously defined and, in this case, they have the same value.

The partial action~$\und\LDAct$ is proper. Indeed, assume that $\und\gg{}_1 \wdots \und\gg{}_\mm$ are elements of~$\und\MM$. As $\pi$ is surjective, there exists for every~$\ii$ an element~$\gg_\ii$ of~$\MM$ that satisfies $\pi(\gg_\ii) = \und\gg{}_\ii$. As $\LDAct$ is proper, there exists~$\xx$ in~$\XX$ such that $\xx \act \gg_\ii$ is defined for each~$\ii$. Then $\LDpibis(\xx) \act \und\gg{}_\ii$ is defined as well. 

We now check that $(\Delta_{\und\xx}){}_{\und\xx \in \und\XX}$ is a right-Garside sequence on~$\und\MM$. First, let $\und\xx$ belong to~$\und\XX$. There exists~$\xx$ in~$\XX$ satisfying $\LDpibis(\xx) = \und\xx$. By assumption, $\xx \act \Delta_\xx$ is defined, hence $\und\xx \act \pi(\Delta_\xx)$ is defined as well, and it is $\LDpibis(\xx \act \Delta_\xx)$. So the sequence $(\Delta_{\und\xx})_{\und\xx \in \und\XX}$ satisfies~(XI.1.11). 

Next, the assumption that $\MM$ is generated by $\bigcup_{\xx \in \XX} \Div(\Delta_\xx)$ implies that $\und\MM$ is generated by $\bigcup_{\xx \in \XX} \pi(\Div(\Delta_\xx))$, which is $\bigcup_{\und\xx \in \und\XX} \Div(\Delta_{\und\xx})$. So, the sequence $(\Delta_{\und\xx})_{\und\xx \in \und\XX}$ satisfies~(XI.1.12). 

Now, assume $\und\sss \in \und\SSS$ and $\und\sss \dive \Delta_{\und\xx}$. Let $\sss$ satisfy $\pi(\gg) = \und\gg$ and $\xx$ satisfy $\LDpibis(\xx) = \und\xx$. By assumption, $\LDpibis(\xx) \act \und\sss$ is defined, hence, by~(XI.4.19), $\xx \act \LDsec(\und\sss)$ is defined. As $\bigcup_\xx \Div(\Delta_\xx)$ generates~$\MM$, the element~$\LDsec(\und\sss)$ is left-divisible by some non-invertible element~$\sss$ that lies in~$\Div(\Delta_\yy)$ for some~$\yy$. As $\SSS$ consists of atoms, every element~$\sss$ of~$\SSS$ belongs to~$\bigcup_\yy \Div(\Delta_\yy)$ and, therefore, by~(XI.1.15), $\sss \dive \Delta_\xx$ must hold for every~$\xx$ such that $\xx \act \sss$ is defined. Applying this to~$\LDsec(\und\sss)$, we deduce $\LDsec(\und\sss) \dive \Delta_\xx$, whence, by~(XI.1.12), $\Delta_\xx \dive \LDsec(\und\sss) \Delta_{\xx \act \LDsec(\und\sss)}$. Applying~$\pi$, we deduce $\Delta_{\und\xx} \dive \und \sss \Delta_{\und\xx \act \und\sss}$. Then an easy induction on the length of~$\und\gg$ shows that $\und\gg \dive \Delta_{\und\xx}$ implies $\Delta_{\und\xx} \dive \und\gg \Delta_{\und\xx \act \und\gg}$ for every~$\und\gg$ in~$\und\MM$. So, the sequence $(\Delta_{\und\xx})_{\und\xx \in \und\XX}$ satisfies~(XI.1.13). 

As of~(XI.1.14), it is automatically satisfied by the sequence $(\Delta_{\und\xx})_{\und\xx \in \und\XX}$ since any two elements of~$\und\MM$ are supposed to admit a left-gcd.

Finally, assume that $\und\gg$ is an element of~$\und\MM$ satisfying $\und\gg \dive \Delta_{\und\xx}$ and $\und\yy \act \und\gg$ is defined. Let~$\xx, \yy$ in~$\XX$ satisfy $\LDpibis(\xx) = \und\xx$, $\LDpibis(\yy) = \und\yy$. Let~$\und\ww$ be an $\und\SSS$-word representing~$\und\gg$ in~$\und\MM$. By construction, $\cl{\und\ww}$ left-divides~$\Delta_{\LDpibis(\xx)}$ hence, by~(XI.4.20), $\cl{\LDsec^*(\und\ww)}$ left-divides~$\Delta_\xx$. By~(XI.4.19) and an induction on the length of~$\und\ww$, the assumption that $\und\yy \act \cl{\und\ww}$ is defined implies that $\yy \act \cl{\LDsec^*(\und\ww)}$ is defined. As the sequence $(\Delta_\xx)_{\xx \in \XX}$ satisfies~(XI.1.15), we deduce that $\cl{\LDsec^*(\und\ww)}$ left-divides~$\Delta_\yy$. Applying~$\pi$, we conclude that $\cl{\und\ww}$, that is, $\und\gg$, left-divides~$\Delta_{\und\yy}$. Hence, the sequence $(\Delta_{\und\xx})_{\und\xx \in \und\XX}$ satisfies~(XI.1.15) and, therefore, it is a right-Garside sequence in~$\und\MM$.
\end{proof}

%%%%
\section*{Solution to selected exercises}

%%%%
\begin{exer}{104}{skeleton}{Z:LDDomain}
Say that a set of addresses is an \emph{antichain} if it does not contain two addresses, one is a prefix of the other; an antichain is called \emph{maximal} if it is properly included in no antichain. \ITEM1 Show that a finite maximal antichain is a family~$\{\alpha_1 \wdots \alpha_\nn\}$ such that every long enough address admits as a prefix exactly one of the addresses~$\alpha_\ii$. 
\ITEM2 Show that, for every $\LDO{}$-word~$\ww$, there exists a unique finite maximal antichain~$\AA_{\ww}$ such that $\term \act \ww$ is defined if and only if the skeleton of~$\term$ includes~$\AA_{\ww}$.
\end{exer}

\begin{solu}
\ITEM2 Use induction on the length of~$\ww$. For $\ww$ of length one, say $\ww = \LDO\alpha$, the result is true with $\term_{\ww}^- = \term_{\ww}^+ = \xx_1 = \{\alpha 10\}$. Assume now $\ww = \LDO\alpha \ww'$. It follows from the definition of the action of~$\LDO\alpha$ that, for every address~$\gamma$, there exists a well defined address~$\gamma'$ such that, for every term~$\term$ such that $\term \act \LDO\alpha$ is defined, the skeleton of~$\term \act \LDO\alpha$ contains~$\gamma$ if and only if the skeleton of~$\term$ contains~$\gamma'$. Put $\gamma' = \LDO\alpha\inv(\gamma)$. Then we claim that the result holds with $\AA_{\ww} = \{\alpha10\} \cup \LDO\alpha\inv(\AA_{\ww'})$. Indeed, $\term \act \ww$ is defined if and only if $\term \act \LDO\alpha$ and $(\term \act \LDO\alpha)\act \ww'$ are defined, hence, by induction hypothesis, if and only if the skeleton of~$\term$ contains~$\alpha 10$ and the skeleton of~$\term \act \LDO\alpha$ includes~$\AA_{\ww'}$.
\end{solu}

%%%%
\begin{exer}{105}{preservation}{Z:LDLcm}
Assume that $\MM$ is a left-cancellative monoid and $\LDAct$ is a partial action of~$\MM$ on a set~$\XX$. \ITEM1 Show that, if $\MM$ admits right lcms (\resp conditional right-lcms), then so does the category~$\CFMX$. 
\ITEM2 Show that, if $\MM$ is right-Noetherian, then so is~$\CFMX$.
\end{exer}

\begin{solu}
By definition, $(\xx, \ff, \yy) \dive (\xx', \ff', \yy')$ in~$\CFMX$ implies $\xx' = \xx$ and $\ff \dive \ff'$ in~$\MM$. So the assumption that $\MM$ is right-Noetherian implies that $\CFMX$ is right-Noetherian as well. Assume that $(\xx, \ff, \yy)$ and $(\xx, \gg, \zz)$ admit a common right-multiple in~$\CFMX$, say $(\xx, \ff, \yy)(\yy, \gg', \xx') = (\xx, \gg, \zz)(\zz, \ff', \xx')$. Then $\ff\gg' = \gg\ff'$ holds in~$\MM$. As $\MM$ is a right-lcm monoid, $\ff$ and $\gg$ admit a right-lcm~$\hh$, and we have $\ff \dive \cc$, $\gg \dive \hh$, and $\hh \dive \ff\gg'$. By assumption, $\xx \act \ff\gg'$ is defined, hence so is $\xx \act \hh$, and it is obvious to check that $(\xx, \hh, \xx\act\hh)$ is a right-lcm of~$(\xx, \ff, \yy)$ and~$(\xx,\gg,\zz)$ in~$\CFMX$. 
\end{solu}

%%%%
\begin{exer}{107}{Noetherianity}{Z:LDMass}
Assume that $\LDAct$ is a proper partial action of a monoid~$\MM$ on a set~$\XX$ and there exists a map $\mu: \XX \to \NNNN$ such that $\mu(\xx\act\gg) > \mu(\xx)$ holds whenever $\gg$ is not invertible. Show that $\MM$ is Noetherian and every element of~$\MM$ has a finite height.
\end{exer}

\begin{solu}
Let $\gg$ be a non-invertible element of~$\MM$, and let $\seqqq{\gg_1}\etc{\gg_\pp}$ be any decomposition of~$\gg$. At the expense of possibly gathering entries, we may assume that each element~$\gg_\ii$ is non-invertible. As $\LDAct$ is proper, there exists~$\xx$ in~$\XX$ such that $\xx \act \gg$ is defined. Then, by~(XI.1.3), $\xx \act \gg_1 \pdots \gg_\ii$ is defined for every~$\ii$. We then obtain $\mu(\xx) < \mu(\xx \act \gg_1) < \mu(\xx \act \gg_1\gg_2) < \pdots < \mu(\xx \act \gg)$, whence $\pp \le \mu(\xx \act \gg) - \mu(\xx)$. So $\gg$ has a finite height bounded above by $\mu(\xx \act \gg) - \mu(\xx)$. In particular, by Proposition~II.2.47 (height), $\LDM$ is Noetherian.
\end{solu}

%%%%
\begin{exer}{109}{common multiple}{Z:LDCommon}
Assume that $\MM$ is a left-cancellative monoid, $\LDAct$ is a partial action of~$\MM$ on~$\XX$, and $(\Delta_\xx)_{\xx\in\XX}$ is a right-Garside sequence in~$\MM$ with respect to~$\LDAct$. Show that, for every~$\xx$ in~$\XX$, any two elements of~$\LDDef\xx$ admit a common right-multiple that lies in~$\LDDef\xx$.
\end{exer}

\begin{solu}
Let $\ff, \gg$ belong to~$\LDDef\xx$. Then $(\xx, \ff, \xx \act \ff)$ and $(\xx, \gg, \xx\act\gg)$ are two elements of~$\CFMX$ sharing the same source. By Proposition~XI.3.25, $\CFMX$ possesses a right-Garside map, hence any two elements of~$\CFMX$ with the same source admit a common right-multiple. So do in particular $(\xx, \ff, \xx \act \ff)$ and $(\xx, \gg, \xx\act\gg)$. A common right-multiple must be of the form $(\xx, \hh, \xx \act \hh)$ where $\hh$ is a common right-multiple of~$\ff$ and~$\gg$ in~$\MM$. Then $\hh$ lies in~$\LDDef\xx$.
\end{solu}

%%%%
\chapter{Chapter~XII: Ordered groups}

%%%%
\section*{Skipped proofs}

\noindent(none)

%%%%
\section*{Solution to selected exercises}

%%%%
\begin{exer}{111}{braid ordering}{Z:ORBraidOrder}
Show that  $1 \lD \sig1\sig2 \lDe(\sig1\sig2\sig1)^{2p} \sigg2\qq$ holds in~$\BB_3$ for $\pp > 0$.
\end{exer}

\begin{solu}
Put $\gg = (\sig1\sig2\sig1)^{2p} \sigg2\qq$ and $\gg' = \sig1\sig2$. We find $\gg'{}\inv \gg = \sig1 (\sig1\sig2\sig1)^{2p-1} \sigg2\qq$. Hence $\gg' \lD \gg$ is true.
\end{solu}

%%%%
\begin{exer}{112}{limit of conjugates}{Z:ORLimit}
Assuming that $\lD$ is a limit of its conjugates in~$\BR3$, show the same result in~$\BR\nn$. [Hint: Use the subgroup of~$\BR\nn$ generated by~$\sig{\nn-2}$ and~$\sig{\nn-1}$, which is isomorphic to~$\BR3$.]
\end{exer}

\begin{solu}
Let $\PP_\nn$ be the positive cone of~$\lD$ on~$\BR\nn$, and let $\HH$ be the subgroup of~$\BR\nn$ generated by~$\sig{\nn-2}$ and~$\sig{\nn-1}$. Then $H$ is isomorphic to $B_3$ and the $\lD$-ordering of~$\BR\nn$ restricted to~$\HH$ corresponds with the $\lD$-ordering of $B_3$.  Moreover, the positive cone for the $\lD$-ordering of~$H$ is $P_\nn \cap H$. Let $S$ be a finite subset of~$P_\nn$. By assumption, there exists~$\ff$ in~$H$ such that $\ff\inv \hh \ff$ is lies in~$\PP_\nn$ for every~$\hh$ in~$S \cap H$, and $\ff (P_\nn \cap H)\ff\inv$ is distinct from~$P_\nn \cap H$. Assume now $\hh \in S \setminus H$. By assumption, $\hh$ is $\sig\ii$-positive for some $\ii < n-2$. Its conjugate~$\ff\inv \hh \ff$ is $\sig\ii$-positive as well, hence it lies in~$\PP_\nn$. We deduce $\ff\inv \SSS \ff \subseteq \PP_\nn$, hence $\SSS \subseteq \ff P_\nn \ff\inv$.  Finally, $\ff P_\nn \ff\inv$ and $P_\nn$ are distinct, because their intersections with $\HH$ are distinct.
\end{solu}

%%%%
\begin{exer}{113}{closure of conjugates}{Z:ORClosure}
Let $\PP_\nn$ be the positive cone of the ordering~$\lD$ on~$\BR\nn$ considered in Example~XII.1.23. Show that the closure of the conjugates of~$\PP_\nn$ in~$\LO{\BR\nn}$ is a Cantor set. 
\end{exer}

\begin{solu}
Let $\ZZ_\nn$ be the family of all conjugates of~$\PP_\nn$. We saw in Example~XII.1.23 that $\PP_\nn$ is a limit of its conjugates, hence it is a limit point in~$\ZZ_\nn$. Hence so is every conjugate of~$\PP_\nn$, that is, every point of~$\ZZ_\nn$ is a limit point in~$\ZZ_\nn$. By continuity, every point in the closure of~$\ZZ_\nn$ is a limit point. Next, the closure of~$\ZZ_\nn$ is a closed subspace in a totally disconnected nonempty compact metric space, hence it is itself totally a disconnected, compact, and metric space. As it is nonempty and every point is a limit point, the closure of~$\ZZ_\nn$ is homeomorphic to the Cantor set.
\end{solu}

%%%%
\begin{exer}{114}{space $\LO{\BR\infty}$}{Z:ORInfiniteBraid}
Show that every point in the space $\LO{\BR\infty}$ is a limit of its conjugates and that $\LO{\BB_\infty}$ is homeomorphic to the Cantor set (contrary to the spaces~$\LO{\BB_\nn}$ for finite~$\nn$).
\end{exer}

\begin{solu}
Consider an arbitrary positive cone $P$ for a left-invariant ordering of~$B_\infty$ and suppose $S$ is a finite subset of~$P$. We will show there is a positive cone $\sig\ii P \siginv\ii$ in~$B_\infty$ that also includes~$S$ and is distinct from~$P$. Choose $\nn$ such that $S$ is included in~$\BB\nn$. Then, for each~$\ii > \nn$, every braid in~$S$ commutes with~$\sig\ii$, so we have $S = \sig\ii S \siginv\ii \subseteq \sig\ii P \siginv\ii$. On the other hand, there exists~$\ii > \nn$ such that the sets~$P$ and $\sig\ii P\siginv\ii$ are different. For otherwise, using $\shift$ for the shift endomorphism that maps~$\sig\ii$ to~$\sig{\ii+1}$ for every~$\ii$, consider the subgroup $\shift^\nn(B_\infty)$, which is isomorphic to~$B_\infty$. The sets $P \cap \shift^\nn(B_\infty)$ and $( \sig\ii P\siginv\ii) \cap \shift^\nn(B_\infty)$ are positive cones for orderings of $\shift^\nn(B_\infty)$.  If $P = \sig\ii P\siginv\ii$ is true for each $\ii >\nn$, then the cone $P \cap \shift^\nn(B_\infty)$ of~$\shift^\nn(B_\infty)$ is invariant under conjugation by all elements of~$\shift^\nn(B_\infty)$.  This would imply that  $\shift^\nn(B_\infty)$ and, therefore, $\BB\infty$ are bi-orderable, which is not true.
\end{solu}

%%%%
\begin{exer}{116}{non-terminating reversing}{Z:ORNotTermin}
Assume that $\Pres\SSS\RR$ is a triangular presentation. \ITEM1 Show that, if a relation of~$\RRh$ has the form $\sss = \ww$ with $\lg\ww > 1$ and $\ww$ finishing with~$\sss$, then the monoid $\PRESp\SSS\RR$ is not of right-$O$-type. \ITEM2 Let $\Pres\SSS\RRh$ be the maximal right-triangular deduced from~$\Pres\SSS\RR$. Show that, if a relation of~$\RRh$ has the form $\sss = \ww$ with $\ww$ beginning with $(\uu\vv)^\rn \uu \sss$ with~$\rn \ge\nobreak 1$, $\uu$~nonempty, and $\vv$ such that $\vv\inv \sss$ reverses to a word beginning with~$\sss$, hence in particular if $\vv$ is empty or it can be decomposed as $\uu_1 \wdots \uu_\mm$ where $\uu_\kk\sss$ is a prefix of~$\ww$ for every~$\kk$, then $\sss\inv \uu \sss$ cannot be terminating, and deduce that$\PRESp\SSS\RR$ is not of right-$O$-type. \ITEM3 Show that a relation $\tta = \ttb\tta\ttb\tta\ttb^3\tta^2...$ is impossible in a right-triangular presentation for a monoid of right-$O$-type.
\end{exer}

\begin{solu}
\ITEM1 If $\RR$ contains a relation $\sss = \uu \sss$ with $\uu$ nonempty, $\sss = \clp\uu \sss$ holds in~$\PRESp\SSS\RR$, whereas $1 = \clp\uu$ fails. So $\PRESp\SSS\RR$ is not right-cancellative. \ITEM2 Write the involved relation $\sss = (\uu \vv)^\rn \uu \sss \ww_1$ with $\vv\inv \sss \revRh \sss \ww_2$. We find 
\begin{align*}
\sss\inv \uu \sss \ &\revRh\  \ww_1\inv \sss\inv (\vv \uu)^{-(\rn-1)} \uu\inv \vv\inv \sss\\
&\revRh\  \ww_1\inv \sss\inv (\vv \uu)^{-(\rn-1)} \uu\inv \sss \ww_2\ \revRh\\  
&\ww_1\inv \sss\inv (\vv \uu)^{-(\rn-1)} (\vv\uu)^{\rn-1} \vv \uu \sss \ww_1 \ww_2\\
&\revRh\  \ww_1\inv \sss\inv \vv \uu \sss \ww_1 \ww_2\ \revRh\  \ww_1\inv \ww_2\inv \cdot \sss\inv  \uu \sss \cdot \ww_1 \ww_2.
\end{align*}
We deduce that $\sss\inv \uu \sss \ \revRh \ (\ww_1\inv \ww_2\inv)^\nn\cdot \sss\inv  \uu \sss \cdot (\ww_1 \ww_2)^\nn$ holds for every~$\nn$ and, therefore, it is impossible that $\sss\inv \uu \sss$ leads in finitely many steps to a positive--negative word. Proposition~II.4.51 (completeness), right-reversing is complete for $(\SSS, \RRh)$, the elements $\sss$ and $\clp\uu \sss$ admit no common right-multiple in~$\PRESp\SSS\RR$, and the latter cannot be of right-$O$-type.  

\ITEM3 The right-hand side of the relation can be written as $(\ttb\tta)\ttb\tta\ttb^2(\ttb\tta)\tta...$, eligible for\ITEM2 with $\uu = \ttb\tta$ and $\vv = \ttb\tta\ttb\cdot\ttb$, a product of two words~$\uu_1, \uu_2$ such that $\uu_\ii\tta$ is a prefix of the right-hand term of the relation.
\end{solu}

%%%%
\begin{exer}{117}{roots of Garside element}{Z:ORDominBis}
Assume that $\MM$ is a left-cancellative monoid generated by a set~$\SSS$ \ITEM1 Show that, for $\delta, \gg$ in a left-cancellative monoid~$\MM$, a necessary and sufficient condition for $\delta$ to right-dominate~$\gg$ is that there exist~$\mm \ge 1$ satisfying $(*)\forall\kk{\ge}0\ (\,\gg\delta^{\kk\mm + \mm-1} \mathbin{\dive} \delta^{\kk\mm +1}\,)$. \ITEM2 Assume that $\delta^\mm$ is a right-Garside element in~$\MM$. Show that $\delta$ right-dominates every element~$\gg$ that satisfies $\gg \delta^{\mm-1} \dive \delta$.
\end{exer}

\begin{solu}
\ITEM1 If $\delta$ right-dominates~$\gg$, then, by definition, $(*)$ holds with $\mm = 1$. Conversely, assume $(*)$. Let $\nn$ be a nonnegative integer. Let $\kk$ be maximal with $\kk\mm \le \nn$. Then we have $\nn \le \kk\mm + \mm - 1$, and $(*)$ implies $\gg\delta^\nn \dive \gg\delta^{\kk\mm + \mm - 1} \dive \delta^{\kk\mm + 1} \dive \delta^{\nn+1}$, so $\delta$ right-dominates~$\sss$. 

\ITEM2 Put $\Delta = \delta^\mm$, and let $\phi$ be the (necessarily unique) endomorphism of~$\MM$ witnessing that $\Delta$ is right-quasi-central. First, we have $\delta \Delta =\nobreak \delta^{\mm+1} = \Delta \phi(\delta)$, whence $\phi(\delta) = \delta$ since $\MM$ is left-cancellative. Next, we claim that $\gg \dive \hh$ implies $\phi(\gg) \dive \phi(\hh)$. Indeed, by definition, $\gg \dive \hh$ implies the existence of~$\hh'$ satisfying $\gg \hh' = \hh$, whence $\phi(\gg) \phi(\hh') = \phi(\hh)$ since $\phi$ is an endomorphism. This shows that $\phi(\gg) \dive \phi(\hh)$ is satisfied. So, in particular, and owing to the above equality, $\gg \dive \delta$ implies $\phi(\gg) \dive \delta$. Now assume $\gg \delta^{\mm-1} \dive \delta$. Then, for every~$\kk$, we find 
$$\gg \delta^{\kk\mm + \mm-1} = \gg \delta^{\mm-1} \Delta^\kk = \Delta^\kk \phi^\kk(\gg\delta^{\mm-1}) \dive \Delta^\kk \phi^\kk(\delta) = \Delta^\kk \delta = \delta^{\kk\mm + 1},$$
and we conclude that $\delta$ right-dominates~$\gg$ by~\ITEM1.
\end{solu} 

%%%%
\begin{exer}{119}{right-ceiling}{Z:ORCentralCeiling}
Assume that $\MM$ is a cancellative monoid of right-$O$-type, and that $\sss_\ell \pdots \sss_1$ is a right-top $\SSS$-word in~$\MM$ such that $\clp{\sss_\ell \pdots \sss_1}$ is central in~$\MM$. Show that $\sss_\ii = \sss_1$ must hold for every~$\ii$, and deduce that ${}^\infty\sss_1$ is the right-$\SSS$-ceiling in~$\MM$.
\end{exer}

\begin{solu}
Let $\Delta = \sss_\ell \pdots \sss_1$. First, we have $\gg \dive \Delta$ for every~$\gg$ in~$\SSS^\ell$, so that $\Delta$ right-dominates~$\SSS^\ell$. Hence, by Lemma~XII.3.9, the right-$\SSS$-ceiling is periodic with period $\sss_\ell \pdots \sss_1$. Now consider its length $\ell+1$ final fragment $\sss_1 \sss_\ell \pdots \sss_1$. Then, in $\MM$, we have $\sss_1 \sss_\ell \pdots \sss_1 = \sss_\ell \pdots \sss_1 \sss_1$, so $\sss_1 \sss_\ell \pdots \sss_1$ and $\sss_\ell \pdots \sss_1 \sss_1$ are two right-top $\SSS$-words of length~$\ell + 1$. By uniqueness of the right-$\SSS$-ceiling, these words must coincide, which is possible only for $\sss_1 = ... = \sss_\ell$. 
\end{solu}

%%%%
\begin{exer}{123}{no triangular presentation]}{Z:ORNoTriangular}
Assume that $\MM$ is a monoid of right-$O$-type that is generated by~$\tta, \ttb, \ttc$ with $\tta \mult \ttb \mult \ttc$ and $\ttb, \ttc$ satisfying some relation $\ttb = \ttc \vv$ with no~$\tta$ in~$\vv$. \ITEM1 Prove that, unless $\MM$ is generated by~$\ttb$ and~$\ttc$, there is no way to complete $\ttb = \ttc\vv$ with a relation $\tta = \ttb \uu$ so as to obtain a presentation of~$\MM$. \ITEM2 Deduce that no right-triangular presentation made of $\ttb = \ttc \ttb \ttc$ (Klein bottle relation) or $\ttb = \ttc \ttb^2 \ttc$ (Dubrovina--Dubrovin braid relation) plus a relation of the form $\tta = \ttb ...$ may define a monoid of right-$O$-type. 
\end{exer}

\begin{solu}
For a contradiction, assume that $(\tta, \ttb, \ttc ; \tta = \ttb \uu, \ttb = \ttc \vv)$ is a presentation of~$\MM$. If there is no~$\tta$ in~$\uu$, the assumption that $\tta = \ttb \uu$ is valid in~$\MM$ implies that $\tta$ belongs to the submonoid generated by~$\ttb$ and~$\ttc$, so $\MM$ must be generated by~$\ttb$ and~$\ttc$.Assume now that there is at least one~$\tta$ in~$\uu$. As $\tta$ does not occur in $\ttb = \ttc \vv$, a word containing~$\tta$ cannot be equivalent to a word not containing~$\tta$. This implies that $\tta$ is preponderant in~$\{\tta, \ttb, \ttc\}$. Indeed, assume that $\gg, \hh$ belong to the submonoid of~$\MM$ generated by~$\ttb$ and~$\ttc$. By the above remark, $\hh \tta \gg' = \gg$ is impossible, hence so is $\hh \tta \dive \gg$. As, by assumption, $\MM$ is of right-$O$-type, we deduce $\gg \dive \hh \tta$. Then Proposition~XII.3.14 gives the result.
\end{solu}

%%%%
\begin{exer}{124}{Birman--Ko--Lee generators}{Z:ORBKL}
Put $\bb_{\ii, \jj} = \aa_{\ii, \jj}^{(-1)^{\ii+1}}$ in the braid group~$\BB_\nn$. Show that, for every~$\nn$, the monoid~$\BDD\nn$ is generated by the elements~$\bb_{\ii, \jj}$.
\end{exer}

\begin{solu}
By Lemma~XII.3.13, an element of~$\BB_\nn$ belongs to~$\BDD\nn$ if and only if it is either $\sig\ii$-positive for some odd~$\ii$ or $\sig\ii$-negative for some even~$\ii$. The braid relations imply $\aa_{\ii, \jj} = \sig{\jj-1} \pdots \sig{\ii+1} \sig\ii \siginv{\ii + 1} \pdots \siginv{\jj-1}$ for $\ii < \jj$, hence $\aa_{\ii, \jj}$ is $\sig\ii$-positive, and $\bb_{\ii, \jj}$ is $\sig\ii$-positive for odd~$\ii$ and $\sig\ii$-negative for even~$\ii$. Therefore, $\bb_{\ii, \jj}$ belongs to~$\BDD\nn$ for all~$\ii, \jj$. Conversely, in~$B_\nn$, we have 
\begin{multline*}
\HS5 \sig\ii \pdots \sig{\nn-1} = (\sig\ii \pdots \sig{\nn-2} \sig{\nn-1} \siginv{\nn-2} \pdots \siginv\ii)\\
(\sig\ii \pdots \sig{\nn-3} \sig{\nn-2} \siginv{\nn-3} \pdots \siginv\ii) \pdots (\sig\ii \sig{\ii+1} \siginv\ii) (\sig\ii),\HS5
\end{multline*}
whence $\sss_\ii = (\sig\ii \pdots \sig{\nn-1})^{(-1){\ii+1}} = \bb_{\ii, \nn} \bb_{\ii, \nn-1} \pdots \bb_{\ii, \ii+1}$ for odd~$\ii$, and $ = \bb_{\ii, \ii+1} \pdots \bb_{\ii, \nn-1} \bb_{\ii, \nn}$ for even~$\ii$. Hence $\sss_\ii$ belongs to the submonoid of~$B_\nn$ generated by the~$\bb_{\ii, \jj}$'s and, finally, $\BDD\nn$ coincides with the latter. 
\end{solu}

%%%%
\chapter{Chapter~XIII: Set-theoretic solutions of YBE}

%%%%
\section*{Skipped proofs}

\noindent(none)

%%%%
\section*{Solution to selected exercises}

%%%%
\begin{exer}{125}{bijective RC-quasigroup}{Z:YBBij}
Assume that $(\XX, \opR)$ is an RC-quasigroup. Let $\double : \XX \to \XX$ and $\Double : \XX\times\XX \to \XX\times\XX$ be defined by $\double(\aa) = \aa\opR\aa$ and $\Double(\aa, \bb) = (\aa\opR\bb, \bb\opR\aa)$. Show that $\Double$ is injective (\resp bijective) if and only if $\double$ is. 
\end{exer}

\begin{solu}
If $\Double$ is injective, then so is $\double$ since $\double(\aa) = \double(\aa')$ implies $\Double(\aa, \aa) = \Double(\aa', \aa')$. On the other hand, if $\Double$ is surjective, then so is $\double$: for every~$\cc$, there exists~$(\aa, \bb)$ satisfying $\Double(\aa, \bb) = (\cc, \cc)$. As seen in the proof of Lemma~XII.2.23\ITEM3, this implies $\aa = \bb$, whence $\double(\aa) = \cc$. For the converse, we first compute $(*)$ $\double(\bb\opR\aa) = (\bb\opR\aa) \opR (\bb\opR\aa) = (\aa\opR\bb) \opR (\aa\opR\aa) = (\aa\opR\bb) \opR \double(\aa)$. Assume that $\double$ is injective and $\Double(\aa, \bb) = \Double(\aa', \bb') = (\cc, \dd)$ holds. By~$(*)$, we have $\cc \opR\double(\aa) = \double(\dd) = \cc \opR \double(\aa')$. As the left-translation by~$\cc$ and~$\double$ are injective, we deduce $\double(\aa) = \double(\aa')$, whence $\aa = \aa'$, and a symmetric argument gives $\bb = \bb'$, so $\Double$ is injective. Finally, assume that $\double$ is bijective and $\cc, \dd$ belong to~$\XX$. As the left-translation by~$\cc$ and~$\double$ are surjective, we can find~$\aa$ satisfying $\cc \opR\double(\aa) = \double(\dd)$, whence, by~$(*)$, $\double(\bb\opR\aa) = \double(\dd)$ and, similarly, we can find~$\bb$ satisfying $\double(\aa\opR\bb) = \double(\cc)$. As $\double$ is injective, we deduce $\aa\opR\bb = \cc$ and $\bb\opR\aa = \dd$, whence $\Double(\aa, \bb) = (\cc, \dd)$. So $\Double$ is surjective, hence bijective.
\end{solu}

%%%%
\begin{exer}{126}{right-complement}{Z:YBRCompl}
Assume that $(\XX, \opR)$ is an RC-quasigr\-oup and $\MM$ is the associated structure monoid. \ITEM1 Show that, for every element~$\ff$ in~$\MM \cap \XX^\pp$, the function from~$\XX$ to~$\XX \cup\nobreak \{1\}$ that maps~$\tt$ to~$\ff\under\tt$ takes pairwise distinct values in~$\XX$ plus at most $\pp$~times the value~$1$. 
\ITEM2 Deduce that, for $\II$ a finite subset of~$\XX$ with cardinal~$\nn$, the right-lcm~$\Delta_\II$ of~$\II$ lies in~$\XX^\nn$.
\end{exer}

\begin{solu}
\ITEM1 We use induction on~$\pp$. For $\pp = 0$, that is, for $\ff = 1$, we have $\ff\under\tt = \tt$, which takes pairwise distinct values in~$\XX$. For $\pp = 1$, that is, for $\ff$ in~$\XX$, we have $\ff\under\ff = 1$ and $\ff\under\tt = \ff\opR\tt$ for $\tt\neq\ff$. As $\tt \mapsto \ff\opR\ff$ injective, the expected result is true. Assume $\pp \ge 2$ and write $\ff = \gg\sss$ with $\gg$ in~$\XX^{\pp-1}$. Then, by the formula for an iterated right-complement, we have $\ff \under \tt = \sss \under (\gg \under \tt)$, whence $\ff\under\tt = 1$ if $\gg\under\tt$ lies in~$\{1, \sss\}$ and $\ff\under\tt = \sss\opR(\gg\under\tt)$ otherwise. Then the result follows from the induction hypothesis. 

\ITEM2 Use induction on~$\nn \ge 0$. For $\nn = 0$, we have $\Delta_\II = 1$. For $\nn = 1$, say $\II = \{\sss\}$ with~$\sss$ in~$\XX$, we have $\Delta_\II = \sss$, an element of~$\XX$. For  $\nn = 2$, say $\II = \{\sss, \tt\}$ with $\sss\not=\tt$, we have $\Delta_\II = \sss(\sss\opR\tt)$, so the induction starts. Assume $\nn \ge 3$ and $\II = \{\sss_1 \wdots \sss_\nn\}$. Let $\JJ = \{\sss_1 \wdots \sss_{\nn-2}\}$. By induction hypothesis, $\Delta_\JJ$ belongs to~$\XX^{\nn-2}$, and $\Delta_{\JJ \cup \{\sss_{\nn-1}\}}$, which is $\Delta_\JJ (\Delta_\JJ \under\sss_{\nn-1})$, belongs to~$\XX^{\nn-1}$. This implies $\Delta_\JJ \under\sss_{\nn-1} \neq1$, so $\Delta_\JJ \under\sss_{\nn-1}$ must lie in~$\XX$. For the same reason, $\Delta_\JJ \under\sss_\nn$ lies in~$\XX$. Moreover, as $\sss_{\nn-1}$ and $\sss_\nn$ are distinct, we have $\Delta_\JJ \under\sss_{\nn-1} \neq \Delta_\JJ \under\sss_\nn$ by~\ITEM1. Now, as $\II$ is $\JJ \cup \{\sss_{\nn-1}\} \cup \{\sss_\nn\}$, Proposition~II.2.12 (iterated lcm) gives $\Delta_\II = \Delta_{\JJ \cup \{\sss_{\nn-1}\}} \cdot (\Delta_\JJ \under \sss_{\nn-1})\under(\Delta_\JJ \under \sss_\nn)$ with a commutative diagram as below, and the question is to know whether the last term may be~$1$.
$$\begin{picture}(95,16)
\pcline{->}(1,12)(29,12)
\taput{$\Delta_\JJ$}
\pcline{->}(31,12)(59,12)
\taput{$\Delta_\JJ \under \sss_{\nn-1} \neq1$}
\pcline{->}(1,0)(29,0)
\pcline{->}(31,0)(59,0)
\pcline{->}(0,11)(0,1)
\tlput{$\sss_\nn$}
\pcline{->}(30,11)(30,1)
\trput{$\Delta_\JJ \under \sss_\nn \neq1$}
\pcline{->}(60,11)(60,1)
\trput{$(\Delta_\JJ \under \sss_{\nn-1}) \under (\Delta_\JJ \under \sss_\nn) \in \XX?$}
\end{picture}$$
Now we saw above that $\Delta_\JJ \under \sss_{\nn-1}$ and $\Delta_\JJ \under \sss_\nn$ are distinct elements of~$\XX$, so the element $(\Delta_\JJ \under \sss_{\nn-1})\under(\Delta_\JJ \under \sss_\nn)$ is equal to $(\Delta_\JJ \under \sss_{\nn-1})\opR(\Delta_\JJ \under \sss_\nn)$, an element of~$\XX$. Hence $\Delta_\II$ belongs to~$\XX^\nn$, as expected.
\end{solu}

%%%%
\begin{exer}{128}{$I$-structure}{Z:YBIStructure}
Assume that $(\XX, \opR)$ is a bijective RC-quasigroup. \ITEM1 Show (by a direct argument) that the map~$\opRu$ from~$\XX^* \times \XX$ to~$\XX$ defined by $1 \opRu \tt = \tt$, $\sss \opRu \tt = \sss\opR\tt$ for $\sss$ in~$\XX$, and $(\uu\sep\vv) \opRu \tt = \vv \opRu (\uu \opRu \tt)$ induces a well-defined map from~$\MM \times \XX$ to~$\XX$. \ITEM2 Show that the map~$\nu$ from~$\XX^*$ to~$\MM$ defined by $\nu(1) = 1$, $\nu(\sss) = \sss$, and $\nu(\ww\sss) = \nu(\ww) \cdot (\nu(\ww)\opRu\sss)$ for~$\sss$ in~$\XX$ induces a well-defined map from~$\NNNN^{(\XX)}$ to~$\MM$.
\end{exer}

\begin{solu}
\ITEM1 Owing to the presentation of~$\MM$, it suffices to check that the translations associated with $\rr\sep(\rr\opR\sss)$ and $\sss\sep(\sss\opR\rr)$ coincide. Now, using the RC-law, we find $(\rr\sep(\rr\opR\sss)) \opRu \tt = (\rr\opR\sss)\opR(\rr\opR\tt) = (\sss\opR\rr)\opR(\sss\opR\tt) = (\sss\sep(\sss\opR\tt)) \opRu \tt$. 

\ITEM2 As $\NNNN^{(\XX)}$ is the quotient of~$\XX^*$ by the equivalence relation generated by all pairs $(\uu\sep\sss\sep\tt\sep\vv, \uu\sep\tt\sep\sss\sep\vv)$ with $\sss, \tt$ in~$\XX$, it suffices to check that the images of such words coincide. Now we find
\begin{align*}
\nu(\uu\sep\sss\sep\tt) 
&= \nu(\uu\sep\sss) \cdot (\nu(\uu\sep\sss)\opRu\tt) 
= \nu(\uu) \cdot (\nu(\uu)\opRu \sss) \cdot (\nu(\uu)\sep(\nu(\uu)\opRu\sss)\opRu\tt)\\
&= \nu(\uu) \cdot (\nu(\uu)\opRu \sss) \cdot ((\nu(\uu)\opRu\sss)\opRu(\nu(\uu)\opRu\tt))\\
&= \nu(\uu) \cdot (\nu(\uu)\opRu \sss) \cdot ((\nu(\uu)\opRu\sss)\opR(\nu(\uu)\opRu\tt))
\end{align*}
and $\nu(\uu\sep\tt\sep\sss) = \nu(\uu) \cdot (\nu(\uu)\opRu \tt) \cdot ((\nu(\uu)\opRu\tt)\opR(\nu(\uu)\opRu\sss))$, whence the equality. Then multiplying by~$\vv$ on the right gives the same result.
\end{solu}

%%%%
\begin{exer}{129}{parabolic submonoid}{Z:YBParab}
\ITEM1 Assume that $(\XX,\sol)$ is a finite involutive nondegenerate set-theoretic solution of YBE and $\MM$ is the associated structure monoid. Show that a submonoid~$\MMu$ of~$\MM$ is parabolic if and only if there exists a (unique) subset~$\II$ of~$\XX$ satisfying $\sol(\II\times\II) = \II\times\II$ such that $\MMu$ is the submonoid of~$\MM$ generated by~$\II$. \ITEM2 Assume that $(\XX, \opR)$ is a   finite   RC-system and $\MM$ is the associated structure monoid. Then a submonoid~$\MMu$ of~$\MM$ is parabolic if and only if there exists a (unique) subset~$\II$ of~$\XX$ such that $\II$ is closed under~$\opR$ and $\MMu$ is the submonoid generated by~$\II$. \ITEM3 Show that, if $(\XX, \opR)$ is an infinite RC-quasigroup, there may exist subsets~$\II$ of~$\XX$ that are closed under~$\opR$ but the induced RC-system $(\II, \opR)$ is not an RC-quasigroup. 
\end{exer}

\begin{solu}
\ITEM1 Assume that $\MMu$ is a parabolic submonoid of $\MM$. Set $\II =\nobreak \MMu \cap \XX$. Since $\XX$ generates $\MM$ and $\MMu$ is closed under factors, $\MMu$ is generated by $\II$. Now, assume that $\aa,\bb$ belong to $\II$ and $(\cc,\dd) = \sol(\aa,\bb)$ holds. Then $\aa\bb$ belongs to~$\MMu$, and $\aa\bb = \cc\dd$ holds in~$\MM$. As $\MMu$ is closed under factor, we deduce that $\cc$ and $\dd$ lie in~$\II$, whence $\sol(\II\times \II) \subseteq \II\times \II$. As $\sol$ is involutive, the latter inclusion must be an equality. 
Conversely assume $\II$ is a subset of $\XX$ satisfying $\sol(\II\times\II) = \II\times\II$ and $\MMu$ is the submonoid of~$\MM$ generated by $\II$. The monoid is (right)-Noetherian, so, in order to prove that $\MMu$ is a parabolic submonoid, and owing to Proposition~VII.1.32 (parabolic subcategory), it is sufficient to establish that $\MMu$ is closed under factor and right-comultiple, that is, in a context where right-lcms exist, that $\MMu$ is closed under factor and right-lcm. Now, an obvious induction shows that every $\XX$-word that is equivalent to an $\II$-word is itself an $\II$-word, implying that $\MMu$ is closed under factor. On the other hand, let $\sol_\II$ be the restriction of~$\sol$ to~$\II \times \II$. The assumption that $(\XX, \sol)$ is a set-theoretic solution of YBE implies that so is $(\II, \sol_\II)$, and the assumption that $(\XX, \sol)$ is involutive implies that so is $(\II, \sol_\II)$. Finally, the assumption that $(\XX, \sol)$ is nondegenerate implies that so is $(\II, \sol_\II)$: with the notation of Definition~XIII.1.8, the left-translations associated with~$\sol_1$ and right-translations associated with~$\sol_2$ are injective on~$\XX$, hence so are their restrictions to~$\II$, and therefore, the latter are bijective since $\II$ is finite. Now, let $\aa, \bb$ belong to~$\II$. As $(\II, \sol_\II)$ is an involutive nondegenerate solution of YBE, there exist~$\cc, \dd$ in~$\II$ satisfying $\aa\cc = \bb\dd$ in~$\MM$, hence in~$\MMu$, and $\MMu$ is closed under right-lcm in~$\MM$. Hence $\MMu$ is a parabolic submonoid of~$\MM$. 

\ITEM2 Translate the result of~\ITEM1 in terms of the operation~$\opR$.

\ITEM3 Consider $\XX = \ZZZZ$ and $\xx \opR \yy = \yy+1$. Then the restriction of~$\opR$ to~$\NNNN$ does not give an RC-quasigroup.
\end{solu}

%%%%
\chapter{Chapter~XIV: More examples}

%%%%
\section*{Skipped proofs}

\noindent(none)

%%%%
\section*{Solution to selected exercises}

\noindent(none)

\vspace{10mm}

\end{document}